\RequirePackage{fix-cm}
\documentclass[smallextended]{svjour3}       
\smartqed  

\let\proof\relax
\let\endproof\relax

\usepackage[pdfencoding=auto]{hyperref}

\usepackage{graphicx}
\usepackage{bm}
\usepackage{booktabs} 
\usepackage[title]{appendix}

\usepackage{bbm}
\usepackage{mathptmx}
\usepackage{graphicx} 
\usepackage{amssymb}

\usepackage{cite}
\usepackage{natbib}
\usepackage{amsmath,mathtools}

\DeclareMathOperator*{\argmax}{argmax}
\DeclareMathOperator*{\argmin}{argmin}
\DeclareMathOperator*{\argsup}{argsup}

\let\Gamma\varGamma
\let\Delta\varDelta
\let\Theta\varTheta
\let\Lambda\varLambda
\let\Xi\varXi
\let\Pi\varPi
\let\Sigma\varSigma
\let\Upsilon\varUpsilon
\let\Phi\varPhi
\let\Psi\varPsi
\let\Omega\varOmega

\makeatletter
\newcommand*{\rom}[1]{\expandafter\@slowromancap\romannumeral #1@}
\makeatother

\usepackage{algorithm}
\usepackage{algorithmic}
\usepackage{amsmath}  

\usepackage{xcolor}

\usepackage{amsthm}

\newtheorem{assumption}{Assumption}

\begin{document}

\title{On Incentive Compatibility in Dynamic Mechanism Design With Exit Option in a Markovian Environment
}


\author{Tao Zhang         \and
        Quanyan Zhu 
}


\institute{Tao Zhang \and Quanyan Zhu \at
              Department of Electrical and Computer Engineering\\
              New York University \\
              \email{\{tz636, qz494\}@nyu.edu}           
}

\date{Received: date / Accepted: date}

\maketitle

\begin{abstract}
This paper studies dynamic mechanism design in a Markovian environment and analyzes a direct mechanism model of a principal-agent framework in which the agent is allowed to exit at any period. 
We consider that the agent's private information, referred to as state, evolves over time.
The agent makes decisions of whether to stop or continue and what to report at each period.
The principal, on the other hand, chooses decision rules consisting of an allocation rule and a set of payment rules to maximize her ex-ante expected payoff.
In order to influence the agent's stopping decision, one of the terminal payment rules is posted-price, i.e., it depends only on the realized stopping time of the agent.
This work focuses on the theoretical design regime of the dynamic mechanism design when the agent makes coupled decisions of reporting and stopping.
A dynamic incentive compatibility constraint is introduced to guarantee the robustness of the mechanism to the agent's strategic manipulation.
A sufficient condition for dynamic incentive compatibility is obtained by constructing the payment rules in terms of a set of functions parameterized by the allocation rule.
The payment rules are then pinned down up to a constant in terms of the allocation rule by deriving a first-order condition.
We show cases of relaxations of the principal's mechanism design problem and provide an approach to evaluate the loss of robustness of the dynamic incentive compatibility when the problem solving is relaxed due to analytical intractability.
A case study is used to illustrate the theoretical results.

\keywords{dynamic mechanism design \and principal-agent problem \and optimal stopping}
\end{abstract}

\section{Introduction}\label{sec:intro}

Mechanism design theory provides a theoretical foundation for designing games that can induce desired outcomes.
The players of the game have private information that is not publicly observable. 
Hence, the mechanism designer's collective decisions have to rely on the players to reveal their private information.
This information asymmetry is an important feature of mechanism design problems.
%
%
The revelation principle allows the mechanism designer to focus on a class of incentive-compatible direct mechanisms to replicate equilibrium outcomes of indirect mechanisms. 
In the celebrated work by Vickery (\cite{vickrey1961counterspeculation}), it has been shown that the seller receives the same expected revenue independent of the mechanism within a large class of auctions.
Vickrey-Clark-Groves (VCG) mechanism is an example of truthful mechanism to achieve a social-optimal solution.
This work investigates mechanism design problems in a dynamic environment, in which a player, a.k.a., the agent, sends a sequence of messages based on the gathered information to the designer.
The designer, a.k.a., the principal, chooses dynamic rules of encounter of the agent to maximize the profit based on the messages. 
%
%
Our model allows the agent to decide how to report his private information to the principal and whether to stop the mechanism immediately or continue to the future at the same time.
These two coupled decision makings depend on the agent's dynamic private information that endogenously depends on the past outcomes of the mechanism.
%
%
The model covers different economic scenarios when the agent establishes an agreement with the principal to, for example, dynamically purchase private or public goods when his valuation stochastically changes over time or frequently consume experience goods when his preference is refined after every usage of the goods, while the agent is allowed to terminate this agreement at any period of time.
Unlike the mechanism design with deadline or a solid commitment period, the agent in our model owns the right to stop.
This additional freedom reduces the agent's risk in the lone-term relationship with the principal due to the uncertainty of the dynamic environment.
However, agent's such freedom complicates the principal's characterization of the incentive compatibility in the dynamic environment.
In this work, we aim to settle the design regimes of such dynamic mechanism by characterizing the allocation and the payment rules and elaborate the incentive compatibility of the mechanism when the agent with time-varying private information couples his decision making of how to reveal his private information with his optimal stopping decision.

Many real-world problems are fundamentally dynamic in nature.
Research of dynamic mechanism design has studied many applications in optimal auctions (e.g., \cite{lin2010dynamic, esHo2007optimal}), screening (e.g., \cite{courty2000sequential,deb2015dynamic,akan2015revenue}), optimal taxation (e.g., \cite{makris2015incentives,findeisen2016education}), contract design (e.g., \cite{williams2011persistent,zhang2009dynamic}), matching market (e.g., \cite{akbarpour2020thickness,anderson2010dynamic}), to name a few.
In dynamic mechanism problems, there are mechanisms without private information. 
For example, in airline revenue management problems, an airline makes decisions about seat pricing on a flight by taking into account the time-varying inventory and the time-evolution of the customer base.
In this paper, however, we consider an information-asymmetric dynamic environment in which the agent privately possesses information that evolves over time. 
The time evolution of the private information may be caused by external factors, the past observations, as well as the decisions from the principal, as when the agent employs learning-by-doing regimes.
For example, in repeated sponsored search auctions, the advertisers privately learn about the profitability of clicks on their ads based on evaluations of the past ads as well as observations from market analysis.
In this work, we consider a dynamic environment, when the agent's private information changes endogenously due to the outcomes of the past decision makings.
Evidence in many economic scenarios has shown that people's past decisions often play a significant role in shaping their future preferences (\cite{Zhang2012Endogenous}). 
Business models have taken into account the belief that customers' preferences can be changed endogenously by using the products or services to form the habit of the customers and encourage long-term shopping sprees.
For example, many products and services, such as newspapers, online video streaming, software as a service, food delivery, provide free trial offers to attract new customers to commit long-term subscriptions.

Optimal stopping theory studies the timing decisions under conditions of uncertainty and has been successfully adopted in applications of economics, finance, and engineering.
Examples include gambling problems (e.g., \cite{he2019optimal,dubins1977countably}), option tradings (e.g., \cite{lundgren2010optimal,lamberton2009optimal,ano2009optimal}), and quick detection problems (e.g., \cite{tartakovsky2008asymptotically,li2014quickest,peng2020sequential}).
This paper studies a class of general dynamic decision-making models, in which the agent has the right to stop the mechanism at any period based on his current observations and the anticipations of the future. 
The agent adopts a stopping rule and is allowed to realize the stopping time in any period before the mechanism terminates naturally upon reaching the final period.
In contract design problems, for example, allowing the stopping decision is made as a clause of a dynamic contract that specifies the agent's right to terminate the agreement at any period.
However, our optimal stopping setting is fundamentally different from the contract intra-deal renegotiation (i.e., renegotiation within the life of the contract) and the contract break.
In general, a renegotiation occurs due to the failure of one party to fulfill its obligations or the inability of one party to meet its commitments.
In such cases, one side of the participants seeks relief of its commitments or wishes to terminate the agreement before the term of that agreement has concluded (see, e.g., \cite{salacuse2000renegotiating}).
The consequence of the renegotiation might be the termination of the contract or a new contract with modified terms and clauses.
Our stopping setting, however, is not due to any breach of the contract by any participant or the inability of maintaining the agreement; it is not a consequence of renegotiations of the contract.
Instead, early termination due to the stopping rule is the right of the agent and lies in the commitment of the principal's dynamic mechanism.
Once the agent decides to terminate the contract at a specific period, he cannot break the dynamic contract by refusing the outcomes (i.e., allocations and payments) that have already been realized up to that period.

We consider a finite-horizon Markovian environment, in which the agent can observe the private information, referred to as the state, that arrives dynamically at the beginning of each period.
The dynamic information structure is governed by a stochastic process characterized by the principal's decision rules, the transition kernels, and the agent's strategic behaviors.
After observing his state at each period, the agent chooses a strategy to report his state to the principal and decides whether to stop immediately or to continue. 
Conditioning on the reported information (including the stopping decision), the principal provides an allocation to the agent and induces a payment. 
The principal aims to maximize her ex-ante expected payoff by choosing feasible decision rules including a set of allocation rules and a set of payment rules.
The principal provides three payment rules including an intermediate payment rule that specifies a payment based on the report when the agent decides to continue and two terminal payment rules.
One of the terminal payment rules is state-dependent and the other is posted-price in the sense that this payment rule depends only on the realized stopping time.
The posted-price payment rule enables the principal to influence the agent's stopping decision without taking into account the agent's private information.
This state-independent terminal payment rule could be the early termination fee to disincentivize the agent from early stopping (when the preferences of the principal and the agent are not aligned), or it could be a reward to elicit the agent to stop at certain periods before the final period to fulfill the principal's interests (when the preferences of the principal and the agent are aligned or partially-aligned).
Under some monotone conditions, the optimal stopping rule can be reformulated as a threshold rule with a time-dependent threshold function.
The threshold rule simplifies the principal's design of the posted-price payment rule as well as the agent's reasoning process for decision makings.
%

%
The design problem in this work faces the challenges from the multi-dimensional interdependence of the agent's joint decision of reporting and stopping at the current period and the planned ones for the future.
On one hand, by fixing the current and the planned future reporting strategies, the agent's stopping decision is made by comparing the payoff if he stops immediately and the best-expected payoff he can anticipate from the future.
On the other hand, with a fixed stopping decision, the (current and the planned future) reporting strategies are chosen by comparing the expected payoffs of different reporting strategies, which determines the expected instantaneous payoff at each period up to the effective time horizon pinned down by the stopping decision.
Therefore, the agent's stopping decision enters the principal's characterization of the dynamic incentive compatibility through this dynamic interdependence.
Given the mechanism, the stopping and the reporting decisions together determine the agent's optimal behaviors.
The coupling of these two decisions in the analysis of incentive compatibility distinguishes this work from other dynamic mechanism design problems.

We define the notion of dynamic incentive compatibility in terms of Bellman equations and address the challenge induced by the agent's dynamic multidimensional decision makings via establishing a one-shot deviation principle (see, e.g., \cite{blackwell1965discounted}).
The one-shot deviation principle has uncovered a foundation of optimality in game theory.
It states that if the agent's deviation from truthful reporting is not profitable for one period, then any finite arbitrary deviations from truthfulness are not profitable.
Monotonicity regarding the designer's allocation rules with respect to the agent's private information is an important result for the implementability of mechanism design.
Consider a single-good auction in which the states are bidders' valuations for a single good and the outcomes are the probabilities for the agent to win the good.
Here, the notion of monotonicity is that the probability of winning the good is non-decreasing in the reported state (see, e.g., \cite{myerson1981optimal, berger2010path}).
\cite{myerson1981optimal} has shown that monotonicity is sufficient for implementability in a single-dimensional domain. However, in general, monotonicity acts only as a necessary condition.
\cite{rochet1987necessary} has constructed a necessary and sufficient condition, called cyclic monotonicity, under which one can design a mechanism such that truthful reporting is optimal for the rational agent.
In this work, we describe a set of monotonicity conditions through inequalities characterized by functions of the allocation rules which we call potential functions.
Given the optimality of the stopping rule, we represent each payment rule by the potential functions.
By applying the envelope theorem, we formulate the potential functions in closed form in terms of the allocation rule.
Our main results, Propositions \ref{proposition:sufficient_ICS}-\ref{prop:revenue_equivalence_new}, provide characterizations of the dynamic incentive compatibility when the agent with time-varying private information makes coupled decisions of reporting and stopping in a Markovian dynamic environment.
The characterizations contribute as design principles by explicitly constructing the three payment rules in terms of the allocation rule and formulating the sufficient and the necessary conditions of dynamic incentive compatibility and facilitate a solid analytical view of dynamic mechanism design to elicit the agent's truthful reporting and optimal stopping.
The sufficient and the necessary conditions yield a revenue equivalence property for the dynamic environment.
We also show that given the threshold function and the allocation rule, the state-independent payment rule is unique up to a constant.
We observe that the posted-price payments from the state-independent payment rule are restricted by a class of regular conditions.
Due to the analytical intractability, relaxation approaches are applied to the principal's mechanism design problem.
We also provide a sufficient condition to evaluate the loss of the robustness of the dynamic incentive compatibility due to the relaxations and the approximations used to solve the principal's problem.

\subsection{Related Work}\label{subsec:related_work}

General settings regarding the source of the dynamics in the related literature can be divided into two categories.
On one hand, the literature on dynamic mechanism design considers the dynamic population of participants with static private information.
\cite{parkes2004mdp} have provided an elegant extension of the social-welfare-maximizing (\textit{efficiency}) VCG mechanism to an online mechanism design framework that studies sequential allocation problems in a dynamic-population environment.
In particular, they have considered the setting when each self-interested agent arrives and departs dynamically over time. 
The private information in their model includes the arrival and the departure time as well as the agent's valuation about different outcomes.
However, the agents do not learn new private information or update their private information.
\cite{pai2013optimal} have proposed a dynamic mechanism model of a similar setting but focusing on the profit-maximization (\textit{optimality}) of the designer. 
Other works focusing on this setting include, e.g., \cite{vulcano2002optimal,gallien2006dynamic, gershkov2009dynamic, said2012auctions, pai2013optimal,board2016revenue}.
On the other hand, there is a number of works studying the problems of the static population where the underlying framework is dynamic because of the time-evolution of the private information.
This category of research has been pioneered by the work of
\cite{baron1984regulation} on the regulation of a monopoly and the contributions of \cite{courty2000sequential} on a sequential screening problem.
There is a large amount of work in this category including, for instance, the dynamic pivot mechanisms (e.g., \cite{bergemann2010dynamic,kakade2013optimal}) and dynamic team mechanisms (e.g., \cite{athey2013efficient,bapna2005efficient,nazerzadeh2013dynamic}).
%
%
%
%
\cite{pavan2014dynamic} have provided a general dynamic mechanism model in which the dynamic of the agents' private information is captured by a set of kernels that is applicable for different behaviors of the time evolution including the learning-based and i.i.d. evolution.
They have used a Myersonian approach and designed a profit-maximizing mechanism with monotonic allocation rules.
\cite{kakade2013optimal} have studied a dynamic virtual-pivot mechanism and provided conditions on the dynamics of the agents' private information.
They have shown an optimal mechanism under the environment they call separable.
\cite{bergemann2010dynamic} and  \cite{parkes2007online} have provided surveys of recent advances in dynamic mechanism design.

The challenges of both settings of dynamics described above come from the information asymmetry between the designer and the agents.
Most of the mechanism design problems study the direct revelation mechanism, in which guaranteeing the incentive compatibility becomes essentially important.
In many dynamic-population mechanism problems with static private information, the incentive compatibility constraints are essentially static (\cite{kakade2013optimal}).
The mechanisms with dynamic private information, however, require efforts to guarantee incentive compatibility.
Monotonicity is an important property of incentive-compatible mechanism design that is widely used in the literature on dynamic mechanism design (see, e.g., \cite{esHo2007optimal,kakade2013optimal,pavan2014dynamic}).

Many situations in economics can be modeled as stopping problems.
There is recent literature on the mechanism design with stopping time. 
\cite{kruse2015optimal} have studied optimal stopping as a mechanism design problem with monetary transfers.
In their model, the agent privately observes a stochastic process of a payoff-relevant state.
The incentive-compatible mechanism considered in their work contains two mechanism rules where one rule maps the agent's report into a stopping decision and the other rule maps the agent's report into a payment that is realized only at the stopping time.
Hence, their model considers that the mechanism determines the optimal stopping decision for the agent based on the report and the agent has only one decision (i.e., reporting) to make.
%
%
%
%
\cite{pavan2009dynamic} have described an application of their dynamic mechanism model to the optimal stopping problem, where the allocation rule provided by the principal is the stopping rule.
Basic formats of stopping rules have been summarized in \cite{lovasz1995efficient} and for rigorous mathematical formulations of general stopping problems, see \cite{peskir2006optimal}.

\subsection{Organization}

The rest of the paper is organized as follows.
In Section \ref{sec:dynamic_environment}, we introduce the dynamic environment and specify the key concepts and notations.
Section \ref{sec:dynamic_principal_agent_problem} formulates the dynamic principal-agent problem and determines the optimal stopping rule for the agent.
In Section \ref{sec:implementability_DIC}, we describe the implementability of the dynamic principal-agent model by defining the dynamic incentive compatibility constraint and reformulating the optimal stopping rule as a threshold rule.
Section \ref{sec:IC_characterization} characterizes the dynamic incentive compatibility by establishing the sufficient and the necessary conditions and showing the properties of revenue equivalence.
In Section \ref{sec:example_optimal}, we formally describe the optimization problem of the principal and show examples of its relaxations due to the analytical intractability. 
We also provide an approach to evaluate the loss of robustness of the dynamic incentive compatibility due to the relaxations.
A case study is given in Section \ref{sec:example_relaxed} as a theoretical illustration.
Section \ref{sec:conclusion} concludes the work.

%

\section{Dynamic Environment}\label{sec:dynamic_environment}

In this section, we describe the dynamic environment of the model.
As conventions, let $\tilde{x}_{t}$ represent the random variable such that $x_{t}\in X_{t}$ is a realized sample of $\tilde{x}$.
By $h^{x}_{s,t}\equiv \{x_{s},...x_{t}\}\in\prod_{k=s}^{t}X_{k}$, we denote the \textit{history} of $x$ from period $s$ up to $t$ (including $t$), with $h^{x}_{t}\equiv h^{x}_{1,t}$, $h^{x}_{t,t}\equiv h^{x}_{t,t-1} \equiv x_{t}$, for all $t\in\mathbb{T}$.
For any measurable set $X$, $\Delta(X)$ is the set of probability measures over $X$.  
Table \ref{table:notations} summarizes the main notations of this paper.

There are two rational (risk neutral, under the expected utility hypothesis) participants in the mechanism: a principal (indexed by $0$, she) and an agent (indexed by $1$, he). 
We consider a finite-time horizon with discrete time $t \in \mathbb{T} \equiv \{1,2,\dots, T\}$.
Let $\mathbb{T}_{t}\equiv \{t,t+1,\dots, T\}$, for all $t\in\mathbb{T}$.
At each period $t$, the agent observes his \textit{state} $\theta_{t}\in \Theta_{t}\equiv{[\underline{\theta}_{t}, \Bar{\theta}_{t}]} \subset \mathbb{R}$. 
Based on $\theta_{t}$, the agent sends a \textit{report} $m_{t}\in M_{t}$ to the principal.
In this paper, we restrict our attention on the \textit{direct mechanism} (\cite{GametheoryDJ}), in which the message space coincides with the state space, i.e., $M_{t}=\Theta_{t}$, for all $t\in\mathbb{T}$.
%
%
%
%
%
Upon receiving the report $\hat{\theta}_{t}$, the principal specifies an \textit{allocation} $a_{t}\in A_{t}$ and a payment $p_{t}\in \mathbb{R}$.
Each $A_{t}$ is a measurable space of all possible allocations.
Let $\sigma_{t}: \prod_{s=1}^{t}\Theta_{s} \times \prod_{s=1}^{t-1}\Theta_{s} \times \prod_{s=1}^{t-1} A_{s} \mapsto \Theta_{t}$ be the \textit{reporting strategy} at $t$, such that $\hat{\theta}_{t}=\sigma_{t} (\theta_{t}|h^{\theta}_{t-1}, h^{\hat{\theta}}_{t-1}, h^{a}_{t-1})$ is the report sent to the principal at period $t$, given his period-$t$ true state $\theta_{t}$ and the histories $h^{\theta}_{t-1}, h^{\hat{\theta}}_{t-1}, h^{a}_{t-1}$.
%
%
The allocation and the payment are chosen, respectively, by the decision rules $\{\alpha, \phi\}$, where $\alpha\equiv \{\alpha_{t} \}_{t\in \mathbb{T} }$ is a collection of (instantaneous) \textit{allocation rules} $\alpha_{t}: \prod_{s=1}^{t}\Theta_{s} \mapsto A_{t}$ and  $\phi\equiv \{\phi_{t} \}_{t\in\mathbb{T}}$ is a collection of (instantaneous) \textit{payment} rules $\phi_{t}:\prod_{s=1}^{t}\Theta_{s}\mapsto \mathbb{R}$, such that the principal specifies an allocation $a_{t}=\alpha_{t}(\hat{\theta}_{t}|h^{\hat{\theta}}_{t-1})$ and a payment $p_{t} = \phi_{t}(\hat{\theta}_{t}|h^{\hat{\theta}}_{t-1})$ when $\hat{\theta}_{t}$ is reported at the current period and $h^{\hat{\theta}}_{t-1}$ has been reported up to $t-1$.

The mechanism allows the agent to leave the mechanism at any period $t\in\mathbb{T}$ by deciding whether to stop or continue at each period according to his optimal stopping rule.
Hence, the agent's decision of $\sigma_{t}$ is coupled with his stopping decision at each period $t$.
Therefore, the principal's characterization of incentive compatibility has to take into account the agent's coupled decision makings.
%
%
To influence the agent's stopping decision, the principal uses a \textit{terminal payment rule} $\rho: \mathbb{T} \mapsto \mathbb{R}$, with $\rho(T)=0$, which is independent of agent's reports and specifies an additional payment $\rho(t)$ at period $t$ if the agent decides to stop at $t$.
To distinguish the intermediate periods and the terminal period, let $\xi_{t}\equiv\phi_{t}$ when the agent realizes his stopping time at period $t$, such that the agent receives a payment $p_{t} = \xi_{t}(\hat{\theta}_{t}|h^{\hat{\theta}}_{t-1}) + \rho(t)$.
%

The mechanism is information-asymmetric because $\theta_{t}$ is privately possessed by the agent for every $t\in \mathbb{T}$ and the principal can learn the true state only through the report $\hat{\theta}_{t}$.
The mechanism is dynamic because the agent's state $\theta_{t}$ evolves endogenously over time and the decisions of both the agent and the principal are made over multiple periods.
The state dynamics lead to the time evolution of probability measures of the states. 
As a result, the expectations of future behaviors at different periods are in general different from each other.
%

\renewcommand\arraystretch{1}  %


\begin{table}[htbp]
	\centering
	\caption{Summary of Notations}\label{table:notations}
	\begin{tabular}{c|c}
		\toprule  
		Notations & Meaning\\
		\toprule  
		$\tilde{\theta}_{t}, \theta_{t}$ & random variable of (period-$t$) state, realization of state \\ 
		\midrule 
		$\sigma=\{\sigma_{t}\}_{t\in\mathbb{T}}$ & reporting strategy \\
		\midrule 
		$\hat{\sigma}[t]=\{\hat{\sigma}[t]_{s}\}_{s\in \mathbb{T}}$ & one-shot deviation (reporting) strategy at $t$\\
		\midrule 
		$\alpha=\{\alpha_{t}\}_{t\in\mathbb{T}}, a_{t}$ & allocation rule, allocation\\
		\midrule 
		$\phi = \{\phi_{t}\}_{t\in \mathbb{T}\backslash\{T\}}$& intermediate payment rule\\
			\midrule 
		$\xi=\{\xi_{t}\}_{t\in \mathbb{T}}$	& state-dependent terminal payment rule \\
		\midrule 
		$\rho$ & state-independent terminal payment rule\\
		\midrule 
		$h^{x}_{s,t}$& history of $x$ from $s$ to $t$, $x\in\{\theta, \tilde{\theta}, \hat{\theta}, a, \tilde{a},\alpha\}$\\
		\midrule 
		 $K_{t}, F_{t}, f_{t}$& transition kernel, c.d.f., and p.d.f. \\
		 \midrule 
		 $\Xi_{\alpha;\sigma}$  & stochastic process of the state \\
		 \midrule 
		 $\Xi_{\alpha;\sigma}\big[h^{\theta}_{t}\big]$ & period-$t$ interim process of the state\\
		 \midrule 
		 $\alpha|\theta_{t}, \hat{\theta}_{t}$ & simplified notation of $\Xi_{\alpha;\hat{\sigma}[t]}\big[h^{\theta}_{t}\big]$ when $\hat{\sigma}[t]_{t}$ reports
		 $\hat{\theta}_{t}$\\
		 \midrule 
		$J^{\alpha, \phi, \xi, \rho}_{0}(\tau;\sigma)$ & the principal's ex-ante expected payoff for (horizon) $\tau\in \mathbb{T}$\\
		 \midrule 
		 $J^{\alpha, \phi, \xi, \rho}_{1}(\tau;\sigma)$ &the agent's ex-ante expected payoff for $\tau\in \mathbb{T}$\\
		 \midrule 
		 $J^{\alpha, \phi, \xi, \rho}_{1,t}(\tau, \theta_{t},\sigma_{t}(\theta_{t});\sigma|h^{\theta}_{t-1})$ & the agent's period-$t$ interim expected payoff for $\tau\in\mathbb{T}_{t}$\\
		 \midrule
		 $\Omega[\sigma]$ & stopping rule when agent adopts the reporting strategy $\sigma$\\
		\bottomrule  
	\end{tabular}
\end{table}


\subsection{Markovian Dynamics}

We consider that the agent's state endogenously evolves over time in a Markovian environment and describe the details of the dynamics and the underlying stochastic process that governs the state dynamics.


\begin{definition}[Markovian Dynamics]\label{def:markov_dynamic}
The Markovian (endogenous) dynamics are characterized by a set of transition kernels $K=\{K_{t}\}_{t\in\mathbb{T}}$, where $K_{t}: \Theta_{t-1} \times \prod_{s=1}^{t-1} A_{s} \mapsto \Delta(\Theta_{t})$ is the period-$t$ transition kernel of the state, i.e., $\tilde{\theta}_{t}\sim K_{t}(\theta_{t-1}, h^{a}_{t-1})$.
Let $F_{t}(\cdot|$ $\theta_{t-1}, h^{a}_{t-1} )$ be the cumulative distribution function (c.d.f.) of $\tilde{\theta}_{t}$, with $f_{t}(\cdot| \theta_{t-1}, h^{a}_{t-1})$ as the probability density function (p.d.f.).
$F_{1}$ and $f_{1}$ are given at the initial period.
\end{definition}

In Markovian dynamics, the generation of the next-period state $\theta_{t+1}$ depends on the current-period state $\theta_{t}$ and the history of allocations $h^{a}_{t}$.
Hence, $\theta_{t+1}$ is independent of the history of true states $h^{\theta}_{t-1}$ given $\theta_{t}$.
From Ionescu-Tulcea theorem (see, e.g., \cite{ash2014real}), the transition kernel $K$, the allocation rule $\alpha$, and the agent's reporting strategy $\sigma$ define a unique stochastic process (i.e., a probability measure) that governs the state dynamics. 
Let $\Xi_{\alpha; \sigma}$ denote the stochastic process.
Given any realizations of period-$t$ state and history of states, we define the period-$t$ \textit{interim process} in the following definition.

\begin{definition}[Interim Process]\label{def:interim_process}
The \textit{interim process} $\Xi_{\alpha;\sigma}\big[h^{\theta}_{t}\big]$ at period $t\in \mathbb{T}$ consists of
\begin{itemize}
    \item[(i)] a deterministic process of the realized $h^{\theta}_{t}\in \prod_{s=1}^{t} \Theta_{s}$ up to time $t$, and
    \item[(ii)] a stochastic process starting from $t+1$ that is uniquely characterized by period-$t$ state $\theta_{t}$, history of allocations $h^{a}_{t}\in \prod_{s=1}^{t} A_{s}$, the allocation rule  $\alpha^{T}_{t+1}\equiv\{\alpha_{s}\}_{s\in \mathbb{T}_{t+1}}$, and the planned reporting strategies $\sigma^{T}_{t+1}=\{\sigma_{s}\}_{s\in \mathbb{T}_{t+1}}$, given the kernels $K^{T}_{t+1}\equiv \{K_{s}\}_{s\in \mathbb{T}_{t+1}}$.
\end{itemize}
\end{definition}

Given the processes $\Xi_{\alpha;\sigma}$ and $\Xi_{\alpha;\sigma}[h^{\theta}_{t}]$, we specify the timing of the mechanism as follows.
\begin{itemize}
    \item[\rom{1}.] \textit{Ex-ante} stage: there is no realization of state, i.e., before the mechanism starts. At this stage, the randomness of the future is characterized by $\Xi_{\alpha;\sigma}$. 
    %
    %
    \item[\rom{2}.] \textit{Interim} stage (at period $t$): $h^{\theta}_{t}$ are realized according to the Markovian dynamics. At each (period-$t$) interim stage:
    \begin{itemize}
        \item[1.] Given the current state $\theta_{t}$, the agent chooses a report $\hat{\theta}_{t}=\sigma_{t}(\theta_{t}|h^{\theta}_{t-1}, h^{\hat{\theta}}_{t-1}, h^{a}_{t-1})$, and decides whether to stop immediately or continue to the next period.
        \item[2.] Upon receiving $\hat{\theta}_{t}$, the principal specifies an allocation $a_{t}= \alpha_{t}(\hat{\theta}_{t}|h^{\hat{\theta}}_{t-1})$, a payment $p_{t} = \phi_{t}(\hat{\theta}_{t}|h^{\hat{\theta}}_{t-1})$ or $p_{t}=\xi_{t}(\hat{\theta}_{t}|h^{\hat{\theta}}_{t-1})+\rho(t)$ if the agent decides to continue or to stop, respectively.
    \end{itemize}
    At each interim stage, the randomness of the future is characterized by $\Xi_{\alpha;\sigma}[h^{\theta}_{t}]$. We denote the corresponding expectations of $\Xi_{\alpha;\sigma}$ and $\Xi_{\alpha;\sigma}[h^{\theta}_{t}]$ by $\mathbb{E}^{\Xi_{\alpha;\sigma}}\big[\cdot \big]$ and $\mathbb{E}^{\Xi_{\alpha;\sigma}[h^{\theta}_{t}]}\big[\cdot \big]$, respectively. We suppress the notation $\sigma$ when it is a truthful reporting strategy.
    %
\end{itemize}

\section{Dynamic Principal-Agent Problem}\label{sec:dynamic_principal_agent_problem}

In this section, we describe the principal-agent problem by identifying their respective objectives.
Let $u_{i,t}: \Theta_{t} \times A_{t} \mapsto \mathbb{R}$ denote the (instantaneous) \textit{utility} of the participant $i$ for $i\in \{0,1\}$ such that $u_{i,t}(\theta_{t}, a_{t})$ is the utility that the participant $i$ receives when the agent's true state is $\theta_{t}$ and the allocation is $a_{t}$ for all $t\in \mathbb{T}$. 
Given any allocation rule $\alpha$, we assume that $u_{i,t}$ is Lipschitz continuous in $\theta_{t}$, for all $\theta_{t}\in \Theta_{t}$, all $t\in \mathbb{T}$.

Given any reporting strategy $\sigma$, define the \textit{ex-ante expected values} of the principal and the agent, respectively, for any time horizon $\tau\in \mathbb{T}$ as follows:
$$
\begin{aligned}
    Z^{\alpha,\phi,\xi}_{0}(\tau;\sigma)
    \equiv&\mathbb{E}^{\Xi_{\alpha;\sigma } }\Big[ \delta^{\tau}\big[ u_{0,\tau}(\tilde{\theta}_{\tau}, \alpha_{\tau}(\sigma_{\tau}(\tilde{\theta}_{\tau}|h^{\tilde{\theta}}_{\tau-1}, h^{\hat{\theta}}_{\tau-1}, h^{\tilde{a}}_{\tau-1})|h^{\hat{\theta}}_{\tau-1}))\\
    -&\xi_{\tau}(\sigma_{\tau}(\tilde{\theta}_{\tau}|h^{\tilde{\theta}}_{\tau-1}, h^{\hat{\theta}}_{\tau-1}, h^{\tilde{a}}_{\tau-1} )|h^{\hat{\theta}}_{\tau-1})\big]\\
    +& \sum_{t=1}^{\tau-1}\delta^{t} \big[u_{0,t}( \tilde{\theta}_{t},   \alpha_{t} (\sigma_{t}(\tilde{\theta}_{t}|h^{\tilde{\theta}}_{t-1}, h^{\hat{\theta}}_{t-1}, h^{\tilde{a}}_{t-1} )| h^{\hat{\theta}}_{t-1})) \\
    -& \phi_{t}(\sigma_{t}(\tilde{\theta}_{t} |h^{\tilde{\theta}}_{t-1}, h^{\hat{\theta}}_{t-1}, h^{\tilde{a}}_{t-1} ) |h^{\hat{\theta}}_{t-1}) \big]  \Big],
\end{aligned}
$$
and
$$
\begin{aligned}
    Z^{\alpha, \phi, \xi}_{1}(\tau;\sigma)
    \equiv& \mathbb{E}^{\Xi_{\alpha;\sigma} } \Big[ \delta^{\tau}\big[ u_{1,\tau}(\tilde{\theta}_{t}, \alpha_{\tau}( \sigma_{\tau}(\tilde{\theta}_{\tau}|h^{\tilde{\theta}}_{\tau-1}, h^{\hat{\theta}}_{\tau-1}, h^{\tilde{a}}_{\tau-1} )|h^{\hat{\theta}}_{\tau-1} )\\
    +&\xi_{\tau}(\sigma_{\tau}(\tilde{\theta}_{\tau}|h^{\tilde{\theta}}_{\tau-1}, h^{\hat{\theta}}_{\tau-1}, h^{\tilde{a}}_{\tau-1} )| h^{\hat{\theta}}_{\tau-1})\big] \\
    +& \sum_{t=1}^{\tau-1} \delta^{t} \big[u_{1,t}( \tilde{\theta}_{t},  \alpha_{t} (\sigma_{t}(\tilde{\theta}_{t}|h^{\tilde{\theta}}_{t-1}, h^{\hat{\theta}}_{t-1}, h^{\tilde{a}}_{t-1} )|h^{\hat{\theta}}_{t-1})) \\
    +& \phi_{t}(\sigma_{t}(\tilde{\theta}_{t}|h^{\tilde{\theta}}_{t-1}, h^{\hat{\theta}}_{t-1}, h^{\tilde{a}}_{t-1} )|h^{\hat{\theta}}_{t-1}) \big]   \Big],
\end{aligned}
$$
where $\delta\in (0,1]$ is the discount factor.
Let $J^{\alpha,\phi,\xi,\rho}_{0}(\cdot;\sigma):\mathbb{T} \mapsto \mathbb{R}$ and $J^{\alpha,\phi,\xi,\rho}_{1}(\cdot;\sigma): \mathbb{T} \mapsto \mathbb{R}$ denote the \textit{ex-ante expected payoffs} of the principal and the agent, respectively, given as follows:
\begin{itemize}
    \item[] \textbf{Principal:}
\begin{equation}\label{eq:ex_ante_payoff_principal}
    \begin{aligned}
        J^{\alpha, \phi, \xi, \rho}_{0}(\tau;\sigma)\; \equiv& Z^{\alpha, \phi, \xi}_{0}(\tau;\sigma) -  \rho(\tau) ; 
    \end{aligned}
\end{equation}
\item[] \textbf{Agent:}
\begin{equation}\label{eq:ex_ante_payoff_agent}
    \begin{aligned}
        J^{\alpha, \phi, \xi, \rho}_{1}(\tau;\sigma) \;\equiv&  Z^{\alpha, \phi, \xi}_{1}(\tau;\sigma) +  \rho(\tau).
    \end{aligned}
\end{equation}
\end{itemize}

Similarly, we define the \textit{interim expected value} of the agent evaluated at period $t$, when the histories are $(h^{\theta}_{t-1}, h^{\hat{\theta}}_{t-1}, h^{a}_{t-1})$, $\theta_{t}$ is observed, and $\sigma_{t}(\theta_{t})$ is reported, as follows, for any $\tau\in\mathbb{T}_{t}$:
\begin{equation}\label{eq:interim_payoff_agent}
    \begin{aligned}
        Z^{\alpha, \phi, \xi}_{1,t}(\tau, \theta_{t},&\sigma_{t}(\theta_{t}|h^{\theta}_{t-1}, h^{\hat{\theta}}_{t-1}, h^{a}_{t-1});\sigma| h^{\theta}_{t-1}) \\
        \equiv& \mathbb{E}^{\Xi_{\alpha;\sigma}[h^{\theta}_{t}]}\Big[ \delta^{\tau}\big[ u_{1,\tau}(\tilde{\theta}_{\tau},\alpha_{\tau}( \sigma_{\tau}(\tilde{\theta}_{\tau}|h^{\theta}_{\tau-1}, h^{\hat{\theta}}_{\tau-1}, h^{a}_{\tau-1})| h^{\hat{\theta}}_{\tau-1}))\\
        +&\xi_{\tau}( \sigma_{\tau}(\tilde{\theta}_{\tau}|h^{\tilde{\theta}}_{\tau-1}, h^{\hat{\theta}}_{\tau-1}, h^{\tilde{a}}_{\tau-1}) |h^{\hat{\theta}}_{\tau-1})) \big]\\
        +& \sum_{s=1}^{\tau-1} \delta^{s} \big[u_{1,s}( \tilde{\theta}_{s},  \alpha_{s} ( \sigma_{s}(\tilde{\theta}_{s}|h^{\theta}_{s-1}, h^{\hat{\theta}}_{s-1}, h^{a}_{s-1})|h^{\hat{\theta}}_{s-1})\\
        +& \phi_{s}(\sigma_{s}(\tilde{\theta}_{s}|h^{\tilde{\theta}}_{s-1}, h^{\hat{\theta}}_{s-1}, h^{\tilde{a}}_{s-1})|h^{\hat{\theta}}_{s-1}) \big]   \Big],
    \end{aligned}
\end{equation}
with $Z^{\alpha, \phi, \xi }_{1,t}(\tau, \theta_{t}|h^{\theta}_{t-1})\equiv Z^{\alpha, \phi, \xi }_{1,t}(\tau, \theta_{t},\theta_{t};\sigma|h^{\theta}_{t-1})$ when $\sigma$ is truthful at every period. 
Then, the corresponding period-$t$ \textit{interim expected payoff} of the agent can be defined as follows: 
\begin{equation}\label{eq:interim_expected_payoff_sigma}
    \begin{aligned}
       J^{\alpha, \phi, \xi, \rho}_{1,t}(\tau, \theta_{t},&\sigma_{t}(\theta_{t}|h^{\theta}_{t-1}, h^{\hat{\theta}}_{t-1}, h^{a}_{t-1});\sigma|h^{\theta}_{t-1}) \\
       =& Z^{ \alpha, \phi, \xi }_{1,t}(\tau, \theta_{t},\sigma_{t}(\theta_{t}|h^{\theta}_{t-1}, h^{\hat{\theta}}_{t-1}, h^{a}_{t-1});\sigma|h^{\theta}_{t-1}) + \rho(\tau),
    \end{aligned}
\end{equation}
with $J^{\alpha, \phi, \xi, \rho}_{1,t}(\tau, \theta_{t}|h^{\theta}_{t-1}) \equiv J^{\alpha, \phi, \xi, \rho}_{1,t}(\tau, \theta_{t},\theta_{t};\sigma|h^{\theta}_{t-1})$ when $\sigma$ is truthful at every period.

%

At the ex-ante stage, the principal provides a \textit{take-it-or-leave-it} offer to the agent by taking into account the agent's stopping rule.
Given the stopping rule, the time horizon $\tau$ of the ex-ante expected payoffs can be calculated.
The principal aims to maximize her ex-ante expected payoff (\ref{eq:ex_ante_payoff_principal}) by anticipating the agent's planned reporting strategy and the stopping time, by choosing the decision rules $\{\alpha,\phi,\xi, \rho\}$.
The agent with his planned reporting strategy $\sigma$ decides whether to accept the offer by checking the following rational participation (RP) constraint:
\begin{equation}\label{eq:rational_part_0}
    \text{RP:}\;\;\;J^{\alpha,\phi,\xi,\rho}_{1}(\tau;\sigma)\geq 0,
\end{equation}
under which, the agent expects a non-negative payoff by participating.
Besides the RP constraint, the principal also wants to incentivize the agent to truthfully reveal his true state in the direct mechanism.
%
%
To achieve this, the principal needs to impose the incentive compatibility (IC) constraint in addition to the RP constraint to guarantee that truthful reporting at each period is for the agent's best interests.


With the knowledge of the Markovian dynamics, the agent can determine his current reporting strategy and plan his future behaviors.
Suppose at period $t$, the agent observes $\theta_{t}$ and reports $\hat{\theta}_{t}$.
The report $\hat{\theta}_{t}$ leads to a current-period allocation $a_{t}=\alpha_{t}(\hat{\theta}_{t}|h^{\hat{\theta}}_{t-1})$.
As a rational participant, the current-period decision making of the agent aims to maximize his interim expected payoff (\ref{eq:interim_expected_payoff_sigma}).
%
Due to the Markovian dynamics, the probability measures of each of the future states depends on the current state $\theta_{t}$ and the history of allocations $h^{a}_{t}$.
%
Hence, the agent's decision of how to report his current-period state is independent of the past true states but depends on $h^{a}_{t-1}$ and $h^{\hat{\theta}}_{t-1}$ (through $h^{a}_{t-1}$), i.e.,
$\hat{\theta}_{t}=\sigma_{t}(\theta_{t}| h^{\hat{\theta}}_{t-1}, h^{a}_{t-1})$. (Hereafter, we will drop $h^{\theta}$ in the notation of the agent's reporting strategy.) 
Therefore, guaranteeing IC when the agent has reported truthfully at all past periods is sufficient to ensure IC when the agent has a history of arbitrary reports.
Given the dependence of the future states on the histories, the agent can plan his future reporting strategies from period $t+1$ onward to obtain an optimal interim expected payoff.
As a result, the principal's craft of IC constraint at each period needs to be coupled with the agent's current and his planned future behaviors.
As described in Section \ref{sec:optimal_stopping}, the agent's adoption of stopping rule further complicates the principal's guarantee of IC.


\subsection{Main Assumptions }

We introduce the following assumptions to support our theoretical analysis of the dynamic model with optimal stopping time.

\begin{assumption}\label{assp:stopping_condition}
Given any $\{\alpha, \phi, \xi, \rho\}$ and $\sigma$, the following holds, for all $t\in\mathbb{T}$,
%
%
\begin{equation}\label{eq:assp_stopping_condition}
    \mathbb{E}^{\Xi_{\alpha;\sigma}}\Big[\sup_{\tau\in\mathbb{T}_{t}}|J^{\alpha,\phi,\xi,\rho}_{1,t}(\tau, \tilde{\theta}_{t}, \sigma_{t}(\tilde{\theta}_{t}| h^{\hat{\theta}}_{t-1}, h^{a}_{t-1}); \sigma|h^{\tilde{\theta}}_{t-1}) | \Big] <\infty.
\end{equation}
\end{assumption}
Assumption \ref{assp:stopping_condition} guarantees that the expected absolute value of any period-$t$ interim expected payoff is bounded for any time horizon $\tau\in \mathbb{T}_{t}$.
This assumption can be satisfied, for example, when the instantaneous payoff $u_{i,t}(\theta_{t}, a_{t}) + p_{t}$ is uniformly bounded for any $\theta_{t}\in \Theta_{t}$, $t\in \mathbb{T}$, given $\{\alpha, \phi, \xi, \rho\}$ and $\sigma$.
This assumption is essential to the existence of optimal stopping rule described in Section \ref{sec:optimal_stopping}.

Define $\chi^{\alpha,\phi,\xi}_{1,t}$ as the difference between the interim expected values when $\tau=t$ and $\tau=t+1$, respectively, evaluated at period $t$ as follows, 
%
\begin{equation}\label{eq:difference_Z}
    \chi^{\alpha, \phi, \xi }_{1,t}(\theta_{t}) \equiv 
    Z^{\alpha, \phi, \xi}_{1,t}(t+1, \theta_{t}|h^{\theta }_{t-1}) - Z^{\alpha, \phi, \xi}_{1,t}(t, \theta_{t}|h^{\theta }_{t-1}).
\end{equation}
Here, $\chi^{\alpha, \phi, \xi }_{1,t}(\theta_{t})$ can be interpreted as the expected marginal change of deferring the stopping time from current period $t$ to the next period $t+1$.
Hence, $\chi^{\alpha, \phi, \xi }_{1,t}(\theta_{t})$ characterizes the marginal incentive of agent's stopping decision at period $t$.
General optimal stopping problems do not require monotonicity assumptions on $\chi^{\alpha, \phi, \xi }_{1,t}(\theta_{t})$.
However, as imposed in many economic models, there is a monotonicity condition referred to as the single crossing property (see, e.g., \cite{quah2012aggregating}): $\chi^{\alpha, \phi, \xi }_{1,t}(\theta_{t})$ crosses the horizontal axis just once, from negative to positive (resp. from positive to negative), as the state $\theta_{t}$ increases (resp. decreases).
A single crossing property applicable in our dynamic model is shown in the following assumption.
\begin{assumption}\label{assp:single_crossing} 
%
$\chi^{\alpha, \phi, \xi }_{1,t}(\theta_{t})$ is non-decreasing in $\theta_{t}$ for all $t\in\mathbb{T}$.
\end{assumption}
This assumption is for the desideratum of establishing a threshold-based optimal stopping rule in Section \ref{sec:threshold_rule}.
To ensure that Assumption \ref{assp:single_crossing} holds, we can, for example, make the derivative of $\chi^{\alpha,\phi,\xi}_{1,t}$ with respective to $\theta_{t}$ non-negative or impose the monotonicity of the utility function, payment rules, and their combinations, e.g., each of $\mathbb{E}^{\Xi_{\alpha;\sigma}[h^{\theta}_{t}] }\Big[\xi_{t+1}(\tilde{\theta}_{t+1}|h^{\theta}_{t})-\xi_{t}(\theta_{t}|h^{\theta}_{t-1})\Big]$, $\phi_{t}(\theta_{t})$, and $\mathbb{E}^{\Xi_{\alpha;\sigma}[h^{\theta}_{t}] }$ $\Big[u_{1,t+1}(\tilde{\theta}_{t+1},$ $ \alpha_{t+1}(\tilde{\theta}_{t+1}|h^{\theta }_{t}))\Big]$ is non-decreasing in $\theta_{t}$.
This assumption can be naturally interpreted in many economic scenarios. For example, suppose that the principal aims to establish a multi-period cooperation with an agent, who has time-evolving productivity (i.e., state) to finish the tasks assigned by the principal over time. Assumption \ref{assp:single_crossing} assumes that the agent with higher productivity is more incentivized by the mechanism to continue instead of stopping immediately than the agent with lower productivity.

\begin{assumption}\label{assp:full_support}
The probability density $f_{t}(\theta_{t}|\theta_{t-1}, h^{a}_{t-1})>0$ for all $\theta_{t}\in \Theta_{t}$, $\theta_{t-1}\in \Theta_{t-1}$, $h^{a}_{t-1}\in \prod_{s=1}^{t-1} A_{s}$, $t\in\mathbb{T}\backslash\{1\}$.
\end{assumption}

Assumption \ref{assp:full_support} restricts our attention to a full support environment where each state has a strictly positive probability to be realized at each period. This assumption is imposed to support the uniqueness of the threshold-based optimal stopping rule and to ensure that the marginal effect of the change of the current state on the future states is bounded (as used in Lemma \ref{lemma:non-decreasing_J}).
Finally, the following assumption imposes a first-order stochastic dominance on the dynamics of the agent's state.

\begin{assumption}\label{assp_monotone_transistions}
For all $\theta'_t \geq \theta_t\in\Theta_{t}$, $\bar{\theta}_{t+1}\in\Theta_{t+1}$, $h^{a}_{t}\in \prod_{s=1}^{t} A_{s}$, $t\in\mathbb{T}\backslash\{T\}$, %
\begin{equation}\label{eq:assp_stochastic_dominance}
F_{t+1}(\bar{\theta}_{t+1}| \theta'_{t}, h^{a}_{t}) \leq F_{t+1}(\bar{\theta}_{t+1}|\theta_{t}, h^{a}_{t}).
\end{equation}
\end{assumption}
Assumption \ref{assp_monotone_transistions} assumes that a larger state at current period leads to a larger state at the next period in the sense of a first-order stochastic dominance, i.e., the partial derivative of $F_{t+1}$ with respect to $\theta_{t}$ is non-positive.
It is straightforward to see that Assumptions \ref{assp:stopping_condition}, \ref{assp:full_support}, and \ref{assp_monotone_transistions} are spontaneously compatible to each other.
%
%
However, since the monotonicity of the state dynamics influences the probability measures of the expectations taken in the term $Z^{\alpha,\phi,\xi}_{1,t}$, the monotonicity of $\chi^{\alpha,\phi,\xi}_{1,t}$ specified by Assumption \ref{assp:single_crossing} is correlated to Assumption \ref{assp_monotone_transistions}.
Given its definition in (\ref{eq:difference_Z}), the monotonicity of $\chi^{\alpha,\phi,\xi}_{1,t}$ also relates to the monotonicity of the utility function, the allocation rule, and the payment rules.
In Section \ref{sec:IC_characterization}, the characterizations of the dynamic incentive compatibility consider the case when all the above four assumptions are satisfied as well as a more general case without Assumption \ref{assp:single_crossing}.

\subsection{Optimal Stopping Rule}\label{sec:optimal_stopping}

In this section, we construct the optimal stopping rule for the agent and identify the dynamic incentive compatibility constraints when the agent makes coupled decisions of reporting and stopping at each period.
Let $\Omega[\sigma]$ denote the agent's stopping rule when he uses reporting strategy $\sigma$.
The optimality of $\Omega[\sigma]$ is defined as follows.

\begin{definition}[Optimal Stopping]\label{def:optimal_stopping}
Given $\{\alpha,\phi,\xi, \rho\}$ and any reporting strategy $\sigma$, the agent's stopping rule $\Omega[\sigma]$ is optimal if there exists a $\tau^{*}$ such that
\begin{equation}\label{eq:optimal_stopping_def}
    \sup_{\tau\in \mathbb{T}} J^{\alpha, \phi, \xi,\rho}_{1}(\tau;\sigma) = J^{\alpha, \phi, \xi, \rho }_{1}(\tau^{*};\sigma).
\end{equation}
Let $\Omega^{*}[\sigma]$ denote the optimal stopping rule when the agent adopts $\sigma$.
%
\end{definition}
Fix $\{\alpha,\phi,\xi,\rho\}$.
To study the optimal stopping problem (\ref{eq:optimal_stopping_def}), we introduce the agent's \textit{valuation function} at period $t$ as follows, for all $\theta_{t}\in \Theta_{t}$, $t\in \mathbb{T}$, any stopping rule $\sigma$:
\begin{equation}\label{eq:stopping_problem_1}
    V^{\alpha,\phi,\xi,\rho}_{t}(\theta_{t};\sigma)\equiv \sup_{\tau\in\mathbb{T}_{t}} J^{\alpha,\phi,\xi,\rho}_{1,t}(\tau, \theta_{t},\sigma_{t}(\theta_{t}| h^{\hat{\theta}}_{t-1}, h^{a}_{t-1});\sigma|h^{\theta}_{t-1}),
\end{equation}
$V^{\alpha,\phi,\xi,\rho}_{t}(\theta_{t})=V^{\alpha,\phi,\xi,\rho}_{t}(\theta_{t};\sigma)$ when $\sigma$ is truthful, 
where the supremum is taken over all time horizon $\tau$ of the process $\Xi_{\alpha}[h^{\theta}_{t}]$ starting from $t$. 
The valuation function $V^{\alpha, \phi, \xi, \rho}_{t}(\theta_{t}; \sigma)$ is the maximum interim expected payoff the agent can expect at period $t$ by varying the time horizon when he has reported $h^{\hat{\theta}}_{t-1}$, observes the current state $\theta_{t}$ and reports $\hat{\theta}_{t}=\sigma_{t}(\theta_{t}| h^{\hat{\theta}}_{t-1}, h^{a}_{t-1})$, and plans to report future states by $\{\sigma_{s}\}_{s\in \mathbb{T}_{t+1}}$.
Note that the agent can modify his planned future reporting strategies at different periods to evaluate the valuation (\ref{eq:stopping_problem_1}).

%

Suppose that Assumption \ref{assp:stopping_condition} holds.
Backward induction leads to the following Bellman equation, for any reporting strategy $\sigma$: 
%

\begin{itemize}
    \item[(i)] for all $t\in \mathbb{T}\backslash \{T\}$, 
    \begin{equation}\label{eq:recurive_problem}
   \begin{aligned}
       &V^{\alpha,\phi,\xi,\rho}_{t}(\theta_{t}; \sigma) \\
       =& \max\Big(  J^{\alpha,\phi,\xi,\rho}_{1,t}(t, \theta_{t},\sigma_{t}(\theta_{t}| h^{\hat{\theta}}_{t-1}, h^{a}_{t-1});\sigma|h^{\theta}_{t-1}), \mathbb{E}^{\Xi_{\alpha}[h^{\theta}_{t}]}\big[ V^{\alpha,\phi,\xi,\rho}_{t+1}(\tilde{\theta}_{t+1};\sigma)\big]  \Big);
   \end{aligned}
\end{equation}
   \item[(ii)] for $t=T$, 
   \begin{equation}\label{eq:final_value_function}
       V^{\alpha,\phi,\xi,\rho}_{T}(\theta_{T};\sigma) = J^{\alpha,\phi,\xi,\rho}_{1,T}(T, \theta_{T},\sigma_{T}(\theta_{T}| h^{\hat{\theta}}_{T-1}, h^{a}_{T-1}); \sigma|h^{\theta}_{T-1}).
   \end{equation}
\end{itemize}
The backward induction yields an equivalent representation of $V^{\alpha,\phi,\xi,\rho}_{t}$ by comparing interim expected payoff if the agent stops at $t$ and the expected valuation of the next period if the agent continues.
The formulations (\ref{eq:recurive_problem}) and (\ref{eq:final_value_function}) naturally lead to the following definition of a stopping region, for any reporting strategy $\sigma$, $t\in\mathbb{T}$:
\begin{equation}\label{eq:original_stopping_zone} 
\begin{split}
    &\Lambda^{\alpha,\phi,\xi,\rho}_{1,t}(t;\sigma)\\
    \equiv& \{\theta_{t}\in\Theta_{t}:V^{\alpha,\phi,\xi,\rho}_{1,t}(\theta_{t};\sigma) =  J^{\alpha,\phi,\xi,\rho}_{1,t}(t, \theta_{t}, \sigma_{t}(\theta_{t}| h^{\hat{\theta}}_{t-1}, h^{a}_{t-1});\sigma|h^{\theta}_{t-1})\}.
\end{split}
\end{equation}
Hence, the period-$t$ stopping region is a set of states that are realized at period $t$ such that the period-$t$ valuation equals the period-$t$ interim expected payoff if the agent stops at $t$.
Based on the stopping region (\ref{eq:original_stopping_zone}), we define the stopping rule, for any reporting strategy $\sigma$:
\begin{equation}\label{eq:stopping_rule_formulation}
\begin{split}
    \Omega^{*}&[\sigma]: \exists \tau\in \mathbb{T}, \text{ s.t., } 
    \tau =\inf\{t\in \mathbb{T}:\theta_{t}\in \Lambda^{\alpha,\phi,\xi,\rho}_{1,t}(t;\sigma)  \}.
\end{split}
\end{equation}
%
%
The stopping rule $\Omega^{*}[\sigma]$ shown in (\ref{eq:stopping_rule_formulation}) calls for stopping at period $t$ if the realized state $\theta_{t}$ is in the stopping region.
Theorem 1.9 of \cite{peskir2006optimal} has shown that $\Omega^{*}[\sigma]$ given in (\ref{eq:stopping_rule_formulation}) solves (\ref{eq:stopping_problem_1}). 
Since $J^{\alpha,\phi, \xi, \rho}(\tau^{*};\sigma) =$  $ \mathbb{E}^{\tilde{\theta}_{1}\sim K_{1}}\Big[$  $V_{1,1}^{\alpha,\phi, \xi, \rho}(\tilde{\theta}_{1};\sigma)\Big]$, $\Omega^{*}[\sigma]$ given in (\ref{eq:stopping_rule_formulation}) is an optimal stopping rule.
Interested readers may refer to Chapter 1 of \cite{peskir2006optimal} for a rigorous characterization of general optimal stopping problems for Markovian processes. 
The coupling of the reporting strategy and the optimal stopping decision is captured in (\ref{eq:original_stopping_zone}) and (\ref{eq:stopping_rule_formulation}).
Given $\{\alpha,\phi,\xi,\rho\}$, the agent chooses the current $\sigma_{t}$ and plans the future $\{\sigma_{s}\}_{s\in\mathbb{T}_{t+1}}$ to maximize the current period valuation
$V^{\alpha,\phi,\xi, \rho}_{t}(\theta_{t};\sigma)$.
For any reporting strategy $\sigma$, $\theta_{t}\in \Theta_{t}$, $t\in\mathbb{T}$, we introduce 
%
%
\begin{equation}\label{eq:sup_time_horizon}
    \begin{split}
        \tau^{\sup}_{t}[\theta_{t};\sigma] = \inf\{\argsup_{\tau\in \mathbb{T}_{t}} J^{\alpha,\phi,\xi,\rho}_{1,t}(\tau, \theta_{t},\sigma_{t}(\theta_{t}| h^{\hat{\theta}}_{t-1}, h^{a}_{t-1});\sigma|h^{\theta}_{t-1})  \}.
    \end{split}
\end{equation}
Here, $\tau^{\sup}_{t}[\theta_{t};\sigma]$ is the smallest time horizon that leads to the maximum period-$t$ interim expected payoff, when the agent uses $\sigma$ and his true state is $\theta_{t}$.
On one hand, the agent's current $\sigma_{t}$ and the planned future $\{\sigma_{s}\}_{s\in\mathbb{T}_{t+1}}$ determine $\tau^{\sup}_{t}[\theta_{t};\sigma]$.
On the other hand, $\tau^{\sup}_{t}[\theta_{t};\sigma]$ determines how far the agent should look into the future to maximize $V^{\alpha,\phi,\xi, \rho}_{t}(\theta_{t};\sigma)$.
This coupling complicates the incentive compatibility and we address the challenge in the next section.




\section{Implementability}\label{sec:implementability_DIC}

In this section, we formally define the incentive compatibility in our dynamic environment and formulate the optimal stopping rule as a threshold-based rule.

\begin{definition}[Dynamic Incentive Compatibility]\label{def:IC_stopping}
The dynamic mechanism $\{\alpha,\phi,$ $\xi, \rho\}$ is dynamic incentive-compatible (DIC) if, for all reporting strategy $\sigma$, 
\begin{itemize}
    \item[(1)] for $t\in\mathbb{T}\backslash\{T\}$,
    \begin{equation}\label{eq:ICS_bellman}
    \begin{aligned}
        &\max\Big(J^{\alpha,\phi,\xi,\rho}_{1,t}(t, \theta_{t}|h^{\theta}_{t-1}),\;\;\; \mathbb{E}^{\Xi_{\alpha}[h^{\theta}_{t}]}\big[ V^{\alpha,\phi,\xi,\rho}_{t+1}(\tilde{\theta}_{t+1})\big]  \Big)\\
        \geq& \max\Big(J^{\alpha,\phi,\xi,\rho}_{1,t}(t, \theta_{t}, \sigma_{t}(\theta_{t}| h^{\hat{\theta}}_{t-1}, h^{a}_{t-1}); \sigma|h^{\theta}_{t-1}), \mathbb{E}^{\Xi_{\alpha;\sigma}[h^{\theta}_{t}]}\big[ V^{\alpha,\phi,\xi,\rho}_{t+1}(\tilde{\theta}_{t+1};\sigma)\big]  \Big);
    \end{aligned}
\end{equation}
\item[(2)] for $t=T$,
\begin{equation}\label{eq:ICS_T}
    J^{\alpha,\phi,\xi,\rho}_{1,T}(T, \theta_{T}|h^{\theta}_{T-1}) \geq J^{\alpha,\phi,\xi,\rho}_{1,T}(T, \theta_{T}, \sigma_{T}(\theta_{T}| h^{\hat{\theta}}_{T-1}, h^{a}_{T-1}); \sigma|h^{\theta}_{T-1});
\end{equation}
\end{itemize}
i.e., the agent of state $\theta_{t}$ maximizes his value at period $t$ by reporting truthfully at all $t\in\mathbb{T}$.
The rules $\{\alpha,\phi,\xi,\rho\}$ are called implementable if the corresponding mechanism is dynamic incentive-compatible.
\end{definition}

At the final period $t=T$, the agent stops with probability $1$.
His incentive to misreport the true state is captured by the immediate instantaneous payoff.
The condition (\ref{eq:ICS_T}) guarantees the non-profitability of misreporting and thus disincentivizes the agent from untruthful reporting.
At each non-final period $t\in \mathbb{T}\backslash\{T\}$, agent's exploration of profitable deviations from truthful reporting strategy takes into account the misreporting of the current state as well as any possible planned future misreporting. 
In particular, when the agent's optimal stopping rule calls for stopping if he reports truthfully, $V_{t}^{\alpha,\phi,\xi, \rho}(\theta_{t}) = J^{\alpha,\phi,\xi, \rho}_{1,t}(t,\theta_{t}|h^{\theta}_{t-1})$; when the agent's optimal stopping rule calls for continuing if he reports truthfully, $V_{t}^{\alpha,\phi,\xi, \rho}(\theta_{t}) = \mathbb{E}^{\Xi_{\alpha}[h^{\theta}_{t}]}\big[ V^{\alpha,\phi,\xi,\rho}_{t+1}(\tilde{\theta}_{t+1}) \big]$.
There are three situations of deviations: \textit{(i)} misreporting at $t$ and stopping at $t$, \textit{(ii)} misreporting at $t$, planned misreporting in the future and continuing, \textit{(iii)} truthful reporting at $t$, planned misreporting in the future and continuing.
The condition (\ref{eq:ICS_bellman}) ensures that no such deviations from truthfully reporting are profitable.

Let $\hat{\sigma}[t]=\{\hat{\sigma}[t]_{s}\}_{s\in\mathbb{T}}$ be any reporting strategy that differs from the truthful reporting strategy $\sigma^{*}$ at only one period $t\in\mathbb{T}$, i.e., $\hat{\sigma}[t]$ is truthful at all periods before $t$ and is planned to be truthful at all periods after $t$.
%
%
We call $\hat{\sigma}[t]$ as a \textit{one-shot deviation} strategy at $t$, for any $t\in\mathbb{T}$. 
To simplify the notation, we denote the process $\Xi_{\alpha;\hat{\sigma}[t]}[h^{\theta}_{t}]$ with $\hat{\theta}_{t}=\hat{\sigma}[t]_{t}(\theta_{t}| h^{\theta}_{t-1}, h^{a}_{t-1})$ by $\alpha|\theta_{t}, \hat{\theta}_{t}$ and let $\alpha|\theta_{t} = \alpha|\theta_{t}, \theta_{t}$ when the agent reports truthfully.
Additionally, we omit $\hat{\sigma}[t]$ in the payoff and the valuation functions and only show $\hat{\theta}_{t}$ as a typical reported state using $\hat{\sigma}[t]_{t}$ unless otherwise stated. 
We denote the corresponding expectation by $\mathbb{E}^{\alpha|\theta_{t}, \hat{\theta}_{t}}[\cdot]$.

\begin{proposition}\label{prop:one_shot_deviation_ICS}
Suppose that Assumption \ref{assp:stopping_condition} holds. In Markovian dynamic environment, if the agent obtains no gain from untruthful reporting by using any one-shot deviation strategy $\hat{\sigma}[t]$, for any $t\in\mathbb{T}$, $\theta_{t}\in \Theta_{t}$, and any history of truthful reports $h^{\theta}_{t-1}\in \prod_{s=1}^{t-1}\Theta_{s}$, then he obtains no gain by using any untruthful reporting strategy for any arbitrary report history $h^{\hat{\theta}}_{t-1}\in \prod_{s=1}^{t-1}\Theta_{s}$, i.e., 
%
%
%
%
\begin{itemize}
    \item[(i)] for $t\in\mathbb{T}\backslash\{T\}$,
    \begin{equation}\label{eq:ICS_bellman_1SD}
    \begin{aligned}
        \max\Big(J^{\alpha,\phi,\xi,\rho}_{1,t}(t, \theta_{t}|h^{\theta}_{t-1}&),\;\;\; \mathbb{E}^{\alpha|\theta_{t} }\big[ V^{\alpha,\phi,\xi,\rho}_{t+1}(\tilde{\theta}_{t+1})\big]  \Big)\\
        \geq \max\Big(&J^{\alpha,\phi,\xi,\rho}_{1,t}(t, \theta_{t},\hat{\theta}_{t}|h^{\theta}_{t-1}),\;\;\; \mathbb{E}^{\alpha|\theta_{t},\hat{\theta}_{t}}\big[ V^{\alpha,\phi,\xi,\rho}_{t+1}(\tilde{\theta}_{t+1})\big]  \Big),
    \end{aligned}
\end{equation}
\item[(ii)] for $t=T$,
\begin{equation}\label{eq:ICS_T_1SD}
    J^{\alpha,\phi,\xi,\rho}_{1,T}(T, \theta_{T}|h^{\theta}_{T-1}) \geq J^{\alpha,\phi,\xi,\rho}_{1,T}(T, \theta_{T}, \hat{\theta}_{T}|h^{\theta}_{T-1}).
\end{equation}
\end{itemize}
\end{proposition}
\proof
See Appendix \ref{app:proof_lemma_one_shot_deviation_ICS}.
\endproof

Proposition \ref{prop:one_shot_deviation_ICS} establishes a \textit{one-shot deviation principle} for our dynamic mechanism.
The one-shot deviation principle enables us to reduce the complexity of the characterization of the DIC and to focus on the analysis of the agent's incentive compatibility at each period when he has reported truthfully at all past periods and plans to report truthfully at all future periods.
As a result, we can restrict attention to the conditions (\ref{eq:ICS_bellman_1SD}) and (\ref{eq:ICS_T_1SD}) as the DIC constraints when the agent uses any one-shot deviation strategy at each period $t\in\mathbb{T}$.
%
%
In the rest of this paper, when it comes to the agent's deviation from truthfulness, we focus on his one-shot deviation strategy $\hat{\sigma}[t]$ that reports $\hat{\theta}_{t} = \hat{\sigma}[t]_{t}(\theta_{t}| h^{\theta}_{t-1}, h^{a}_{t-1})$ at any single period $t\in\mathbb{T}$.
With a slight abuse of notation, we use $h^{\hat{\theta}_{t}}_{s}$, $s\geq t$, to denote the history of reports (including the planned reports) when the agent uses one-shot deviation strategy at $t$ and reports $\hat{\theta}_{t}$.

Define the \textit{continuing value} as, for any $\theta_{t}, \hat{\theta}_{t}\in \Theta_{t}$, $t\in\mathbb{T}\backslash\{T\}$,
\begin{equation}\label{eq:continuing_value}
    \mu^{\alpha,\phi,\xi,\rho}_{t}(\theta_{t},\hat{\theta}_{t}) \equiv \sup_{\tau\in\mathbb{T}_{t+1}}\mathbb{E}^{\alpha|\theta_{t},\hat{\theta}_{t}}\Big[ J^{\alpha,\phi,\xi,\rho}_{1,t+1}(\tau,\tilde{\theta}_{t+1}|h^{\theta}_{t}) \Big] - J^{\alpha,\phi,\xi,\rho}_{1,t}(t, \theta_{t},\hat{\theta}_{t}|h^{\theta}_{t-1}),
\end{equation}
with $\mu^{\alpha,\phi,\xi,\rho}_{t}(\theta_{t},\theta_{t}) \equiv \mu^{\alpha,\phi,\xi,\rho}_{t}(\theta_{t})$ when the agent uses the truthful reporting strategy.
The continuing value captures the agent's maximum expected gain when he decides to continue to the next period instead of stopping immediately.
Then, the stopping rule in (\ref{eq:stopping_rule_formulation}) can be characterized by the continuing value as follows, for any one-shot deviation strategy $\hat{\sigma}[t]$ at any $t\in\mathbb{T}$:
\begin{equation}\label{eq:optimal_stopping_rewrite_1}
    \Omega^{*}[\hat{\sigma}[t]]: \exists \tau\in\mathbb{T}, \text{ s.t., } \tau = \inf\{t'\in\mathbb{T}: \mu^{\alpha,\phi,\xi,\rho}_{t'}(\theta_{t'},\hat{\theta}_{t'}) \leq  0 \}.
\end{equation}
Define the \textit{marginal value} $L^{\alpha, \phi, \xi,\rho }_{t}$ as follows, for any $\theta_{t}, \hat{\theta}_{t}\in\Theta_{t}$, $t\in \mathbb{T}\backslash \{T\}$:
\begin{equation}\label{eq:change_interim_payoff_agent}
    \begin{aligned}
        L^{\alpha, \phi, \xi,\rho }_{t}(\theta_{t}, \Hat{\theta}_{t})\equiv   \mathbb{E}^{\alpha|\theta_{t}, \Hat{\theta}_{t}}\Big[\delta^{t+1}\big[ u_{1,t+1}(\tilde{\theta}_{t+1}, \alpha_{t+1}(\tilde{\theta}_{t+1}|h^{\hat{\theta}_{t}}_{t} )) + \xi_{t+1}(\tilde{\theta}_{t+1}|h^{\hat{\theta}_{t}}_{t})\big]&\\
        +\rho(t+1)\Big] + \delta^{t}\big[\phi_{t}(\hat{\theta}_{t}|h^{\theta}_{t-1})-\xi_{t}(\Hat{\theta}_{t}|h^{\theta}_{t-1})\big]&,
    \end{aligned}
\end{equation}
with $L^{\alpha, \phi, \xi,\rho}_{t}(\theta_{t})=$  $L^{\alpha, \phi, \xi,\rho }_{t}(\theta_{t}, \theta_{t})$ when the agent uses the truthful reporting strategy.
Hence, $L^{\alpha, \phi, \xi,\rho }_{t}$ $(\theta_{t},\Hat{\theta}_{t})$ $- \rho(t)$ captures the marginal change in the interim expected payoffs evaluated at $t$ if the agent plans to stop at period $t+1$ instead of stopping immediately at $t$ when his true state is $\theta_{t}$ and he reports $\hat{\theta}_{t}$.

\begin{lemma}\label{lemma:marginal_value_represt} 
Given $\{\alpha, \phi, \xi,\rho\}$, we have, for any $\tau\in \mathbb{T}$,
\begin{equation}\label{eq:ex_ante_payoff_marginal}
    J^{\alpha,\phi, \xi,\rho}_{1}(\tau) = \mathbb{E}^{\Xi_{\alpha}}\Big[ \sum_{s=1}^{\tau-1} L^{\alpha, \phi, \xi, \rho }_{s}(\tilde{\theta}_{s}) -\rho(s)\Big] + J^{\alpha, \phi, \xi, \rho}_{1}(1),
\end{equation} 
and, for any $\theta_{t}\in \Theta_{t}$, $h^{\theta}_{t}\in \prod_{s=1}^{t-1} \Theta_{s}$, $t\in\mathbb{T}$, $\tau\in\mathbb{T}_{t}$,
\begin{equation}\label{eq:interim_payoff_marginal}
    \begin{aligned}
        J^{\alpha,\phi, \xi, \rho }_{1,t}(\tau, \theta_{t}|h^{\theta}_{t-1}) =& \mathbb{E}^{\alpha|\theta_{t}}\Big[\sum_{s=t}^{\tau-1}L_{s}^{\alpha, \phi, \xi, \rho }(\tilde{\theta}_{s}) -\rho(s) \Big] + J^{\alpha, \phi, \xi,\rho }_{1,t}(t, \theta_{t}|h^{\theta}_{t-1}).
    \end{aligned}
\end{equation}
\end{lemma}
\proof
See Appendix \ref{app:proof_lemma:marginal_value_represt}.
\endproof

\begin{lemma}\label{lemma:monotonicity_of_L} Suppose that Assumption \ref{assp:single_crossing} holds. Then, $L^{\alpha, \phi, \xi, \rho }_{t}(\theta_{t})$ is non-decreasing in $\theta_{t}$, for all $t\in \mathbb{T}\backslash \{T\}$.
\end{lemma}

The proof of Lemma \ref{lemma:monotonicity_of_L} directly follows the formulation of $J^{\alpha,\phi,\xi, \rho}_{1}$ in Lemma \ref{lemma:marginal_value_represt}.
Lemma \ref{lemma:marginal_value_represt} shows that the agent's ex-ante expected payoff and his period-$t$ interim expected payoff can be represented in terms of $L^{\alpha, \phi, \xi, \rho }$, $\rho$, and the expected payoffs when the agent stops at the starting periods (i.e., period $1$ or period $t$, respectively).
%
Lemma \ref{lemma:monotonicity_of_L} establishes a single crossing condition that is necessary for the existence of the threshold-based stopping rule (see Section \ref{sec:threshold_rule}).

From Lemma \ref{lemma:marginal_value_represt}, we can represent the continuing value in terms of $L^{\alpha, \phi, \xi, \rho}_{t}$ as follows, for any $\theta_{t}, \hat{\theta}_{t}\in \Theta_{t}$, $t\in\mathbb{T}\backslash \{T\}$:
\begin{equation}
    \begin{aligned}
       \mu^{\alpha, \phi, \xi, \rho}_{t}(\theta_{t}, \hat{\theta}_{t})\equiv \sup_{\tau\in\mathbb{T}_{t+1}}\Big[\mathbb{E}^{\alpha|\theta_{t},\hat{\theta}_{t}}\big[\sum_{s=t+1}^{\tau-1}L^{\alpha,\phi,\xi,\rho}_{s}(\tilde{\theta}_{s})&- \rho(s)  \big]  \Big] \\
       &+ L^{\alpha,\phi,\xi,\rho}_{t}(\theta_{t},\hat{\theta}_{t}) -\rho(t).
    \end{aligned}
\end{equation}
Let, for any $\theta_{t}, \hat{\theta}_{t}\in \Theta_{t}$, $t\in\mathbb{T}\backslash \{T\}$, %
\begin{equation}\label{eq:short_continuing_value}
    \bar{\mu}^{\alpha,\phi,\xi,\rho}_{t}(\theta_{t},\hat{\theta}_{t}) \equiv \mu^{\alpha,\phi,\xi,\rho}_{t}(\theta_{t}, \hat{\theta}_{t}) + \rho(t),
\end{equation}
with $\bar{\mu}^{\alpha,\phi,\xi,\rho}_{t}(\theta_{t},\theta_{t}) \equiv \bar{\mu}^{\alpha,\phi,\xi,\rho}_{t}(\theta_{t})$.
%
%
Then, we define the following auxiliary functions, for any $\theta_{t}, \hat{\theta}_{t}\in \Theta_{t}$, $h^{\theta}_{t-1}\in \prod_{s=1}^{t-1}\Theta_{s}$, $t\in\mathbb{T}$, 
\begin{equation}\label{eq:lemma_rationable_1}
    U_{S,t}^{\alpha,\phi, \xi, \rho}(\theta_{t},\hat{\theta}_{t}|h^{\theta}_{t-1})\equiv \delta^{t}\big[ u_{1,t}(\theta_{t}, \alpha_{t}(\hat{\theta}_{t}|h^{\theta}_{t-1})) + \xi_{t}(\hat{\theta}_{t}|h^{\theta }_{t-1})\big]+\rho(t),
\end{equation}
and
\begin{equation}\label{eq:lemma_rationable_2}
    \begin{aligned}
        U_{\Bar{S},t}^{\alpha,\phi, \xi, \rho}(\theta_{t},\hat{\theta}_{t}|h^{\theta}_{t-1})\equiv  \delta^{t}\big[u_{1,t}(\theta_{t}, \alpha_{t}(\hat{\theta}_{t}|h^{\theta}_{t-1})) &+\xi_{t}(\hat{\theta}_{t}|h^{\theta }_{t-1})\big]+ \bar{\mu}^{\alpha,\phi, \xi,\rho}_{t}(\theta_{t}, \hat{\theta}_{t}),
    \end{aligned}
\end{equation}
with $U_{j,t}^{\alpha,\phi, \xi, \rho}(\theta_{t},\theta_{t}|h^{\theta}_{t-1}) \equiv U_{j,t}^{\alpha,\phi, \xi, \rho}(\theta_{t}|h^{\theta}_{t-1})$, for $j\in\{S,\bar{S}\}$, where
the subscripts $S$ and $\Bar{S}$ represent ``stop" and ``non-stop", respectively.
Basically, $U^{\alpha,\phi, \xi, \rho}_{S,t}$ and $U^{\alpha,\phi, \xi, \rho}_{\bar{S},t}$ are the agent's expected payoffs that are directly determined by his current reporting strategy $\hat{\sigma}[t]_{t}$ and planned truthful reporting strategies $\{\hat{\sigma}[t]_{s}\}_{s\in \mathbb{T}_{t+1}}$, when he decides to stop at $t$ and continue to $t+1$, respectively.
We can rewrite the DIC in Proposition \ref{prop:one_shot_deviation_ICS} in terms of these auxiliary functions in the following lemma.
\begin{lemma}\label{lemma:ICS_equivalent}
The IC constraints (\ref{eq:ICS_bellman_1SD}) and (\ref{eq:ICS_T_1SD}) are equivalent to the following
\begin{equation}\label{eq:O-LA-IC_1}
    U_{S,t}^{\alpha,\phi, \xi, \rho}(\theta_{t}|h^{\theta}_{t-1})\geq U_{S,t}^{\alpha,\phi, \xi, \rho}(\theta_{t}, \hat{\theta}_{t}|h^{\theta}_{t-1}),
\end{equation}
and
\begin{equation}\label{eq:O-LA-IC_2}
    U_{\bar{S},t}^{\alpha,\phi, \xi, \rho}(\theta_{t}|h^{\theta}_{t-1})\geq U_{\bar{S},t}^{\alpha,\phi, \xi, \rho}(\theta_{t}, \hat{\theta}_{t}|h^{\theta}_{t-1}),
\end{equation}
for all $\theta_{t}$, $ \hat{\theta}_{t}\in \Theta_{t}$, $h^{\theta}_{t-1}\in \prod_{s=1}^{t-1} \Theta_{s}$, $t\in\mathbb{T}$.
\end{lemma}

Conditions (\ref{eq:O-LA-IC_1}) and (\ref{eq:O-LA-IC_2}) reformulate the incentive constraints in (\ref{eq:ICS_bellman_1SD}) and (\ref{eq:ICS_T_1SD}). They ensure that misreporting is not profitable in the instantaneous payoff when the agent stops and continues, respectively.
If (\ref{eq:O-LA-IC_1}) (resp. (\ref{eq:O-LA-IC_2})) is satisfied, the optimality of the stopping rule guarantees that continuing (resp. stopping) is not profitable when the stopping rule calls for stopping (resp. continuing).

\subsection{Threshold Rule}\label{sec:threshold_rule}

In this section, we revisit the optimal stopping rule and introduce a class of threshold-based stopping rule (see, e.g., \cite{kruse2018inverse,jacka1992finite, villeneuve2007threshold}).
Given the definition of $\bar{\mu}^{\alpha, \phi,\xi, \rho}_{t}(\theta_{t}, \hat{\theta}_{t})$ in (\ref{eq:short_continuing_value}), we can rewrite the optimal stopping rule in (\ref{eq:optimal_stopping_rewrite_1}) as, for any $\hat{\sigma}[t]$ that reports $\hat{\theta}_{t}\in\Theta_{t}$, $t\in\mathbb{T}$,
\begin{equation}\label{eq:optimal_stoppinng_represent_2}
\begin{split}
    \Omega^{*}[\hat{\sigma}[t]]: \exists \tau, \text{ s.t., } \tau = \inf\{t\in \mathbb{T}:\bar{\mu}^{\alpha,\phi,\xi,\rho}_{t}(\theta_{t}, \hat{\theta}_{t}) \leq \rho(t) \},
\end{split}
\end{equation}
with the corresponding stopping region
\begin{equation}\label{eq:rewrite_stopping_region}
    \Lambda^{\alpha,\phi,\xi,\rho}_{1,t}(t;\hat{\theta}_{t}) = \{\theta_{t}\in\Theta_{t}: \bar{\mu}^{\alpha,\phi,\xi,\rho}_{t}(\theta_{t}, \hat{\theta}_{t}) \leq \rho(t) \}.
\end{equation}
Hence, the principal can adjust $\rho(t)$ to influence the agent's stopping time by changing the stopping region $\Lambda^{\alpha,\phi,\xi,\rho}_{1,t}(t;\hat{\theta}_{t})$.
In general, the stopping rule partitions the state space $\Theta_{t}$ into multiple zones, which leads to multiple stopping sub-regions. Suppose that, given $<\alpha,\phi,\xi,\rho>$ and $\hat{\sigma}[t]$, there are $n_{t}$ stopping sub-regions at period $t$. Let $\{\theta^{\ell;k}_{t}, \theta^{r;k}_{t}\}$, $\underline{\theta}_{t}\leq\theta^{\ell;k}_{t}\leq \theta^{r;k}_{t} \leq \bar{\theta}_{t}$, denote the boundaries of the $k$-th stopping sub-region, such that the stopping region (\ref{eq:rewrite_stopping_region}) is equivalent to 
\begin{equation}\label{eq:rewrite_stopping_region_v2}
\vec{\Lambda}^{\alpha,\phi,\xi,\rho}_{1,t}(t;\hat{\theta}_{t}|n_{t}) =\{\theta_{t}\in [\theta^{\ell;k}_{t}, \theta^{r;k}_{t}], \forall k = 1,\dots, n_{t}\},
\end{equation}
with $\vec{\Lambda}^{\alpha,\phi,\xi,\rho}_{1,t}(t|n_{t})=\vec{\Lambda}^{\alpha,\phi,\xi,\rho}_{1,t}(t;\theta_{t}|n_{t})$ when the agent reports truthfully.

If there is only one stopping sub-region, i.e., the stopping rule partitions the state space $\Theta_{t}$ into two regions, then the stopping rule is a threshold rule.
The existence of a threshold rule depends on the monotonicity of $\bar{\mu}^{\alpha,\phi,\xi,\rho}_{t}$ with respect to $\theta_{t}$. The following lemma directly follows Lemma \ref{lemma:monotonicity_of_L}. 
\begin{lemma}\label{lemma:nonotone_continuing_value}
Suppose that Assumption \ref{assp:single_crossing} holds.
Then, $\bar{\mu}^{\alpha, \phi, \xi, \rho }_{t}(\theta_{t}, \hat{\theta}_{t})$ is non-decreasing in $\theta_{t}$, for all $t\in \mathbb{T}\backslash \{T\}$, any $\hat{\theta}_{t}\in \Theta_{t}$.
\end{lemma}

Since $\rho(t)$ is independent of states or reports, the monotonicity in Lemma \ref{lemma:nonotone_continuing_value} suggests a threshold-based stopping rule for the agent. 
%
%
Let $\eta: \mathbb{T} \rightarrow \Theta_{t}$ be the \textit{threshold function}, for all $t\in\mathbb{T}$, such that the agent chooses to stop the first time the state $\theta_{t}\leq \eta(t)$.
Since the agent has to stop at the final period, we require $\eta(T)=\bar{\theta}_{T}$. 
Let $\Omega[\hat{\sigma}[t]]|\eta$ denote the stopping rule with the threshold function $\eta$.
\begin{definition}[Threshold Rule]\label{def:threshold_rule}
Fix any $\hat{\sigma}[t]$, $t\in\mathbb{T}$.
We say that the stopping rule $\Omega[\hat{\sigma}[t]]|\eta$ is a threshold rule if there exists a threshold function $\eta$ such that
\begin{equation}\label{eq:threshold_based_stopping}
\Omega[\hat{\sigma}[t]]|\eta: \exists \tau, \text{ s.t, } \tau= \inf\{t\in \mathbb{T}: \theta_t \leq \eta(t)\}.
\end{equation}
\end{definition}

We establish that the optimal stopping rule shown in (\ref{eq:optimal_stoppinng_represent_2}) is a threshold rule in Lemma \ref{prop:threshold_rule}.

\begin{lemma}\label{prop:threshold_rule} Suppose that Assumption \ref{assp:single_crossing} holds.
If the stopping rule $\Omega[\hat{\sigma}[t]]$ is optimal in the mechanism $\{\alpha,\phi,\xi,\rho\}$, then it is a threshold rule with a threshold function $\eta$, denoted as $\Omega[\hat{\sigma}[t]]|\eta$.
\end{lemma}
\proof
See Appendix \ref{app:prop_threshold_rule}. 
\endproof

Lemma \ref{lemma:threshold_rule} shows the uniqueness of the threshold function.
\begin{lemma}\label{lemma:threshold_rule}
Suppose that Assumptions \ref{assp:single_crossing} and \ref{assp:full_support} hold. Then, each threshold rule has a unique threshold function.
\end{lemma}
\proof
See Appendix \ref{app:lemma_threshold_rule}.
\endproof

From (\ref{eq:threshold_based_stopping}), the payment rule $\rho$ can be fully characterized by the threshold function $\eta$ such that the principal can influence the agent's stopping time by manipulating $\eta$. For example, setting $\eta(t)=\underline{\theta}_{t}$ (resp. $\eta_{t} = \bar{\theta}_{t}$) forces the agent to continue (resp. stop) with probability $1$.

%

%


\section{Characterization of Incentive Compatibility}\label{sec:IC_characterization}

%

In this section, we characterize the dynamic incentive compatibility of the mechanism. We first introduce the length functions and the potential functions and show a sufficient condition by constructing the payment rules in terms of the potential functions based on the relationship between the length and the potential functions.
Second, we pin down the payment rules in terms of the allocation rule by applying the envelope theorem.
%

Let, for any $\theta_{t}$, $\hat{\theta}_{t}\in\Theta_{t}$, $\tau\in\mathbb{T}_{t}$, $t\in\mathbb{T}\backslash{\{T\}}$,
\begin{equation}\label{eq:sub_path_length}
    \begin{split}
        \pi^{\alpha}_{t}(\theta_{t},\hat{\theta}_{t};\tau) = \mathbb{E}^{\alpha|\theta_{t},\hat{\theta}_{t}}\Big[\delta^{t}u_{1,t}(&\theta_{t}, \alpha_{t}(\hat{\theta}_{t}|h^{\theta}_{t-1}))+\sum_{s=t+1}^{\tau}\delta^{s}u_{1,s}(\tilde{\theta}_{s}, \alpha_{s}(\tilde{\theta}_{s}|h^{\hat{\theta}_{t}}_{s-1})) \\
        &+\sum_{s=t+1}^{\tau-1} \delta^{s}\phi_{s}(\tilde{\theta}_{s}|h^{\hat{\theta}_{t}}_{s-1}) +\delta^{\tau}\xi_{\tau}(\tilde{\theta}_{\tau}|h^{\hat{\theta}_{t}}_{\tau-1})  +\rho(\tau)\Big],
    \end{split}
\end{equation}
with $\pi^{\alpha}_{t}(\theta_{t};\tau)=\pi^{\alpha}_{t}(\theta_{t},\theta_{t};\tau)$.
The term $\pi^{\alpha}_{t}(\cdot;\tau)$ denotes the agent's period-$t$ expected payoff to go for a time horizon $\tau>t$ without the current-period payment specified by $\phi_{t}$. When $\tau=t$, the term $\pi^{\alpha}_{t}(\cdot;\tau)$ is the agent's period-$t$ instantaneous payoff if he stops at $t$.
Define the \textit{length functions} as, for any $\theta_{t}$, $\hat{\theta}_{t}\in\Theta_{t}$, $t\in\mathbb{T}$,
\begin{equation}\label{eq:path_length_eta_stop}
    \ell^{\alpha}_{S,t}(\hat{\theta}_{t}, \theta_{t}) \equiv \delta^{t} u_{1,t}(\hat{\theta}_{t}, \alpha_{t}(\hat{\theta}_{t}|h^{\theta}_{t-1})) - \delta^{t} u_{1,t}(\theta_{t}, \alpha_{t}(\hat{\theta}_{t}|h^{\theta}_{t-1})),
\end{equation}
and, for any $\tau\in\mathbb{T}_{t}$,
\begin{equation}\label{eq:path_length_eta_continue}
    \begin{aligned}
        \ell^{\alpha}_{\bar{S},t}(\hat{\theta}_{t}, \theta_{t};\tau) =& \pi^{\alpha}_{t}(\hat{\theta}_{t};\tau) - \pi^{\alpha}_{t}(\theta_{t},\hat{\theta}_{t};\tau).
    \end{aligned}
\end{equation}
%
%
The length function $\ell_{\bar{S},t}(\hat{\theta}_{t},\theta_{t})$ (resp. $\ell_{S,t}(\hat{\theta}_{t},\theta_{t})$) describes the change in the value of $\pi^{\alpha}_{t}$ (resp. $\delta u_{1,t}$) while keeping reporting $\hat{\theta}_{t}$ when the agent observes $\hat{\theta}_{t}$ instead of observing $\theta_{t}$.
From the definition of $\pi^{\alpha}_{t}$, it is straightforward to see that $\ell^{\alpha}_{S,t}(\hat{\theta}_{t}, \theta_{t}) = \ell^{\alpha}_{\bar{S},t}(\hat{\theta}_{t}, \theta_{t};t)$.
%

Let $\beta^{\alpha}_{S,t}(\cdot):\Theta_{t}\rightarrow \mathbb{R}$ and $\beta^{\alpha}_{\bar{S},t}(\cdot):\Theta_{t}\rightarrow \mathbb{R}$ be the \textit{potential functions} that depend only on $\alpha$.
Proposition \ref{proposition:sufficient_ICS} shows a sufficient condition for DIC by showing the constructions of the payment rules in terms of the potential functions and the threshold function (for $\rho$) and specifying the relationships between the potential functions and the length functions.



\begin{proposition}\label{proposition:sufficient_ICS}
Fix an allocation rule $\alpha$ and a threshold function $\eta$. 
Suppose that Assumptions \ref{assp:stopping_condition}, \ref{assp:full_support}, and \ref{assp_monotone_transistions} hold.
The dynamic mechanism is dynamic incentive-compatible if the following statements are satisfied:

%
\begin{itemize}
    \item[(i)] The payment rules $\phi$, $\xi$, and $\rho$, respectively, are constructed as, for all $t\in\mathbb{T}$,
\begin{equation}\label{eq:payment_phi_eta_1}
    \begin{aligned}
        \phi_{t}(\theta_{t}) =&  \delta^{-t}\beta^{\alpha}_{\Bar{S},t}(\theta_{t}) - \delta^{-t}\mathbb{E}^{\alpha|\theta_{t}}\Big[\beta^{\alpha}_{\bar{S},t+1}(\tilde{\theta}_{t+1}) \Big] -u_{1,t}(\theta_{t}, \alpha_{t}(\theta_{t}|h^{\theta}_{t-1} )),
    \end{aligned}
\end{equation}
\begin{equation}\label{eq:payment_xi_eta}
    \xi_{t}(\theta_{t}) = \delta^{-t}\beta^{\alpha}_{S,t}(\theta_{t})-u_{1,t}(\theta_{t},\alpha_{t}(\theta_{t}|h^{\theta}_{t-1})),
\end{equation}
\begin{equation}\label{eq:construct_rho_primal}
    \begin{split}
        \rho(t) = \delta^{-t}\mathbb{E}^{\alpha|\eta(t)}\Big[\sum_{s=t}^{T-1} \big(&\beta^{\alpha}_{S,s+1}(\tilde{\theta}_{s+1}\vee \eta(s+1))-\beta^{\alpha}_{S,s}(\tilde{\theta}_{s}\vee \eta(s)) \big) \\
        &- \big(\beta^{\alpha}_{\bar{S},s+1}(\tilde{\theta}_{s+1}\vee \eta(s+1))-\beta^{\alpha}_{\bar{S},s}(\tilde{\theta}_{s}\vee \eta(s))   \big) \Big].
    \end{split}
\end{equation}

\item[(ii)] $\beta^{\alpha}_{\bar{S},t}$ and $\beta^{\alpha}_{S,t}$ satisfy,
%
for all $\theta_{t}$, $\hat{\theta}_{t}\in \Theta_{t}$, $t\in\mathbb{T}$,
\begin{equation}\label{eq:potential_function_relation_S}
    \begin{aligned}
    &\beta^{\alpha}_{S,t}(\hat{\theta}_{t})- \beta_{S,t}^{\alpha}(\theta_{t})\leq \ell^{\alpha}_{S,t}(\hat{\theta}_{t},\theta_{t}),
    \end{aligned}
\end{equation}
\begin{equation}\label{eq:potential_function_relation_nS}
    \beta^{\alpha}_{\bar{S},t}(\hat{\theta}_{t})-\beta^{\alpha}_{\bar{S},t}(\theta_{t})\leq \inf_{\tau\in\mathbb{T}_{t}}\Big\{\ell^{\alpha}_{\bar{S},t}(\hat{\theta}_{t},\theta_{t};\tau)\Big\}-\sup_{\tau'\in\mathbb{T}_{t}}\rho(\tau'),
\end{equation}
and 
\begin{equation}\label{eq:potential_function_relation_ine}
    \beta^{\alpha}_{\bar{S},t}(\theta_{t}) \geq \beta^{\alpha}_{S,t}(\theta_{t}) \text{ and } \beta^{\alpha}_{\bar{S},T}(\theta_{T}) = \beta^{\alpha}_{S,T}(\theta_{T}).
\end{equation}

\item[(iii)] $\beta^{\alpha}_{\bar{S},t}$ and $\beta^{\alpha}_{S,t}$ are constructed in terms of $\alpha$ such that the following term is non-decreasing in $\theta_{t}$, for all $t\in\mathbb{T}$:
\begin{equation}\label{eq:monotonicity_potential_assumption_2}
    \begin{aligned}
   \chi^{\alpha}_{1,t}(\theta_{t}) = \mathbb{E}^{\alpha|\theta_{t}}\Big[\beta^{\alpha}_{S,t+1}(\tilde{\theta}_{t+1}) - \beta^{\alpha}_{\bar{S},t+1}(\tilde{\theta}_{t+1}) \Big] -\big(\beta^{\alpha}_{S,t}(\theta_{t})-\beta^{\alpha}_{\bar{S},t}(\theta_{t}) \big).
\end{aligned}
\end{equation}
\end{itemize}

\end{proposition}
\proof
See Appendix \ref{app:proposition_sufficient_ICS}. 
\endproof

In Proposition \ref{proposition:sufficient_ICS}, the first statement \textit{(i)} shows three conditions (\ref{eq:payment_phi_eta_1})-(\ref{eq:construct_rho_primal}) to construct the three payment rules in terms of the allocation rule and the threshold function (for $\rho$) through the potential functions and the utility function, such that if the potential functions satisfy the conditions (\ref{eq:potential_function_relation_S})-(\ref{eq:potential_function_relation_ine}) in the second statement \textit{(ii)}, then the mechanism $\{\alpha, \phi,\xi,\rho\}$ is DIC.
Importantly, the construction of $\rho$ in (\ref{eq:construct_rho_primal}) is based on the setting when the agent's optimal stopping rule is a threshold rule with a unique threshold function. The existence of such threshold rule is established under Assumption \ref{assp:single_crossing}.
Therefore, Proposition \ref{proposition:sufficient_ICS} imposes an additional condition for the construction of the potential functions in the third statement \textit{(iii)} to guarantee the monotonicity in Assumption \ref{assp:single_crossing}.
The explicit constructions of the potential functions is constructed from the necessity of the DIC.

Define with a slight abuse of notation, for any $\theta_{t}, \hat{\theta}_{t}\in\Theta_{t}$, $\tau\in\mathbb{T}_{t}$, $t\in\mathbb{T}$,
$$
\bar{\mu}^{\alpha,\phi,\xi,\rho}_{t}(\theta_{t}, \hat{\theta}_{t};\tau) \equiv \mathbb{E}^{\alpha|\theta_{t},\hat{\theta}_{t}}\big[\sum_{s=t}^{\tau-1}L^{\alpha,\phi,\xi,\rho}_{s}(\tilde{\theta}_{s})- \rho(s)  \big]   + \rho(t),
$$
with $\bar{\mu}^{\alpha,\phi,\xi,\rho}_{t}(\theta_{t};\tau)= \bar{\mu}^{\alpha,\phi,\xi,\rho}_{t}(\theta_{t}, \theta_{t};\tau)$,
and, for any $\theta_{t}, \hat{\theta}_{t}\in\Theta_{t}$, $h^{\theta}_{t-1}\in\prod_{s=1}^{t-1}\Theta_{s}$, $\tau\in\mathbb{T}_{t}$, $t\in\mathbb{T}$
\begin{equation}\label{eq:define_gamma}
\begin{split}
    U^{\alpha,\phi,\xi,\rho}_{t}&(\tau, \theta_{t}, \hat{\theta}_{t}|h^{\theta}_{t-1})\equiv\\
    &\begin{cases}
\delta^{t}\big[ u_{1,t}(\theta_{t}, \alpha_{t}(\hat{\theta}_{t}|h^{\theta}_{t-1})) + \xi_{t}(\hat{\theta}_{t}|h^{\theta}_{t-1})\big]+\rho(t), & \text{ if } \tau=t,\\
\delta^{t}\big[ u_{1,t}(\theta_{t}, \alpha_{t}(\hat{\theta}_{t}|h^{\theta}_{t-1})) + \xi_{t}(\hat{\theta}_{t}|h^{\theta}_{t-1})\big]+\bar{\mu}^{\alpha,\phi,\xi,\rho}_{t}(\theta_{t};\tau), & \text{ if } \tau>t,
\end{cases}
\end{split}
\end{equation}
with $U^{\alpha,\phi,\xi,\rho}_{t}(\tau, \theta_{t}|h^{\theta}_{t-1})=U^{\alpha,\phi,\xi,\rho}_{t}(\tau, \theta_{t}, \theta_{t}|h^{\theta}_{t-1})$.

The following lemma takes advantage of the quasilinearity of the payoff function and formulates the partial derivative of $U^{\alpha,\phi,\xi,\rho}_{t}(\tau, \theta_{t}|h^{\theta}_{t-1})$ with respect to $\theta_{t}$.

\begin{lemma}\label{lemma:envelope_conditions}
Suppose that Assumption \ref{assp:full_support} holds.
In any DIC mechanism, $\{\alpha,\phi,\xi,\rho\}$ satisfy, for any $\theta_{t}\in \Theta_{t}$, $h^{\theta}_{t-1}\in\prod_{s=1}^{t-1}\Theta_{s}$, $\tau\in\mathbb{T}_{t}$, $t\in\mathbb{T}$,
\begin{equation}\label{eq:envelope_gamma}
    \begin{aligned}
       \frac{ \partial U^{\alpha,\phi,\xi,\rho}_{t}(\tau,r|h^{\theta}_{t-1}) }{ \partial r }\Big|_{r= \theta_{t}} = & \mathbb{E}^{\alpha|\theta_{t}}\Big[ \sum_{s=t}^{\tau} \delta^{s}\frac{\partial u_{1,s}(r, \alpha_{s}(\tilde{\theta}_{s}|h^{\theta}_{s-1})) }{\partial r }\Big|_{r=\tilde{\theta}_{s}} G_{t,s}(h^{\tilde{\theta}}_{t,s})   \Big],
    \end{aligned}
\end{equation}
where
$$
G_{t,s}(h^{\theta}_{t,s}) =\prod_{k=t+1}^{s} \Big[-\frac{\partial F_{k}(\theta_{k}|\theta_{k-1}, \alpha_{k-1}(\theta_{k-1}|h^{\theta}_{k-2})) }{f_{k}(\theta_{k}|\theta_{k-1}, \alpha_{k-1}(\theta_{k-1}|h^{\theta}_{k-2}))\partial r  }\Big|_{r=\theta_{k-1}}\Big].
$$
\end{lemma}

\proof
See Appendix \ref{app:proof_lemma:envelope_conditions}.
\endproof

We establish the explicit formulations of the potential functions in the following proposition.

\begin{proposition}\label{prop:envelope_necessary}
Fix any allocation rule $\alpha$. 
Suppose that Assumptions \ref{assp:stopping_condition}, \ref{assp:full_support}, and \ref{assp_monotone_transistions} hold.
Suppose additionally that $u_{1,t}(\theta_{t}, a_{t})$ is a non-decreasing function of $\theta_{t}$.
In any DIC mechanism $\{\alpha, \phi, \xi, \rho\}$, with $\phi$, $\xi$, and $\rho$ constructed in (\ref{eq:payment_phi_eta_1})-(\ref{eq:construct_rho_primal})  respectively, 
the potential functions $\beta^{\alpha}_{S,t}$ and $\beta^{\alpha}_{\bar{S},t}$ are constructed in terms of $\alpha$ as follows, for any arbitrary fixed state $\theta_{\epsilon,t}\in \Theta_{t}$, any $\theta_{t}$, $\hat{\theta}_{t}\in\Theta_{t}$, $h^{\theta}_{t-1}\in\prod_{s=1}^{t-1}\Theta_{s}$, $t\in\mathbb{T}$,
\begin{equation}\label{eq:beta_necessary_S}
     \beta^{\alpha}_{\boldsymbol{S},t}(\theta_{t}) = \int^{\theta_{t} }_{\theta_{\epsilon,t}}\gamma^{\alpha}_{t}(t,r|h^{\theta}_{t-1})dr,
\end{equation}
\begin{equation}\label{eq:beta_necessary_SBar}
         \beta^{\alpha}_{\boldsymbol{\bar{S}},t}(\theta_{t}) = \sup_{\tau\in\mathbb{T}_{t}}\Big\{ \int^{\theta_{t}}_{\theta_{\epsilon,t}}\gamma^{\alpha}_{t}(\tau, r|h^{\theta}_{t-1})dr \Big\},
\end{equation}
where $\gamma_{t}^{\alpha}(\tau,\theta_{t}|h^{\theta}_{t-1}) \equiv \frac{\partial U^{\alpha,\phi,\xi,\rho}_{t}(\tau,r|h^{\theta}_{t-1}) }{ \partial r }\Big|_{r= \theta_{t}} $, where $\alpha$ satisfies the conditions (\ref{eq:potential_function_relation_S})-(\ref{eq:potential_function_relation_ine}).
\end{proposition}
\proof

See Appendix \ref{app:prop_formulation_betas}.
\endproof

Proposition \ref{prop:envelope_necessary} provides explicit formulations of the potential functions that depend only on the allocation rule $\alpha$, up to a constant shift (determined by $\theta_{\epsilon,t}$).
As a result, the payment rules $\phi$ and $\xi$ can be pinned down up to a constant given $\alpha$.
Such results lead to the celebrated revenue equivalence theorem, which is important in mechanism design problems in static settings (e.g., \cite{myerson1981optimal, vickrey1961counterspeculation}) as well as in dynamic environments (e.g., \cite{pavan2014dynamic,kruse2018inverse} ). The following proposition summarizes the property of the revenue equivalence of our dynamic mechanism.

\begin{proposition}\label{prop:revenue_equivalence_new}
Fix an allocation rule $\alpha$. 
Suppose that Assumptions \ref{assp:stopping_condition}, \ref{assp:full_support}, and \ref{assp_monotone_transistions} hold.
Suppose additionally that $u_{1,t}(\theta_{t}, a_{t})$ is a non-decreasing function of $\theta_{t}$.
Let $\theta^{a}_{\epsilon,t}$, $\theta^{b}_{\epsilon,t}\in \Theta_{t}$ be any two arbitrary states.
Let $\phi^{a} \equiv \{\phi^{a}_{t}\}$ and $\phi^{b}\equiv \{\phi^{b}_{t}\}$ \big(resp. $\xi^{a}$ and $\xi^{b}$\big) be two payment rules constructed according to (\ref{eq:payment_phi_eta_1}) \big(resp. (\ref{eq:payment_xi_eta})\big) by the potential function $\beta^{\alpha}_{\bar{S},t}$ formulated in (\ref{eq:beta_necessary_SBar}) \big(resp. $\beta^{\alpha}_{S,t}$ formulated in (\ref{eq:beta_necessary_S})\big) with $\theta^{a}_{\epsilon,t}$ and $\theta^{b}_{\epsilon,t}$, respectively, as the lower limit of the integral.
Define $\rho^{a}$ and $\rho^{b}$ in the similar way, for some threshold functions $\eta^{a}$ and $\eta^{b}$, respectively.
Then, there exist constants $\{C_{\tau}\}_{\tau\in\mathbb{T}}$ 
such that, for any $\theta_{t}\in\Theta_{t}$, $\tau\in\mathbb{T}_{t}$, $t\in\mathbb{T}$,
\begin{equation}\label{eq:equivalence_1}
    \begin{aligned}
        \mathbb{E}^{\alpha|\theta_{t}}\Big[ \sum_{s=1}^{\tau-1}\delta^{s}\phi^{a}_{s}(\tilde{\theta}_{s}|&h^{\tilde{\theta}}_{s-1}) + \delta^{\tau}\xi^{a}_{\tau}(\tilde{\theta}_{\tau}|h^{\tilde{\theta}}_{t-1}) +\rho^{a}(\tau) \Big] \\
    &= \mathbb{E}^{\alpha|\theta_{t}}\Big[ \sum_{s=1}^{\tau-1}\delta^{s}\phi^{b}_{s}(\tilde{\theta}_{s}|h^{\tilde{\theta}}_{s-1})+\delta^{\tau}\xi^{b}_{\tau}(\tilde{\theta}_{\tau}|h^{\tilde{\theta}}_{\tau-1}) +\rho^{b}(\tau) \Big] + C_{\tau}.
    \end{aligned}
\end{equation}

\end{proposition}
\proof
See Appendix \ref{app:proof_prop_revenue_equivalence_new}.
\endproof

Proposition \ref{prop:revenue_equivalence_new} implies that different DIC mechanisms with the same allocation rule result in equivalent expected sum of payments paid by the principal up to a constant for any time horizon $\tau\in\mathbb{T}_{t}$ given any $t\in\mathbb{T}$.
The following corollary shows the uniqueness of the state-independent payment rule $\rho$ for any given threshold function $\eta$.

\begin{corollary}\label{prop:uniqueness_rho}
Fix an allocation rule $\alpha$.
Suppose that Assumptions \ref{assp:stopping_condition}, \ref{assp:full_support}, and \ref{assp_monotone_transistions} hold.
Suppose additionally that $u_{1,t}(\theta_{t}, a_{t})$ is a non-decreasing function of $\theta_{t}$.
Let $\eta$ be any threshold function such that $\eta(t)\in\Theta_{t}$, for each $t\in\mathbb{T}$, and let $\theta^{a}_{\epsilon,t}$, $\theta^{b}_{\epsilon,t}\in \Theta_{t}$ be any two arbitrary states.
Construct two state-independent payment rules $\rho^{a}$ and $\rho^{b}$ according to (\ref{eq:construct_rho_primal}) with $\eta$ in which the potential functions are formulated in (\ref{eq:beta_necessary_S}) and (\ref{eq:beta_necessary_SBar}) with $\theta^{a}_{\epsilon,t}$ and $\theta^{b}_{\epsilon,t}$, respectively, as the lower limit of the integral.
Then, there exist constants $\{C^{\rho}_{t}\}_{t\in\mathbb{T}}$ such that, for any $t\in\mathbb{T}$, 
\begin{equation}\label{eq:prop_uniqueness_rho}
    \rho^{a}(t) = \rho^{b}(t) + C^{\rho}_{t},
\end{equation}
with $\rho^{a}(T)=\rho^{b}(T)=C^{\rho}_{t}=0$.
\end{corollary}

If the mechanism does not satisfy the monotonicity condition specified in Assumption \ref{assp:single_crossing}, we cannot guarantee a threshold rule (with a unique threshold function) for the mechanism. Hence, as shown in (\ref{eq:rewrite_stopping_region_v2}), the optimal stopping rule partitions the state space into multiple stopping sub-regions.
Let $\vec{\Theta}_{t}\equiv \{\theta^{k;j}_{t}\}_{k=1, j=\{k, r\} }^{k=n_{t}}$ denote the set of boundaries of stopping region $\vec{\Lambda}^{\alpha,\phi,\xi,\rho}_{1,t}(t|n_{t})$ given in (\ref{eq:rewrite_stopping_region_v2}).
Define an operator 
$$
\theta_{t}	\sqcup \vec{\Theta}_{t} \equiv \begin{cases}
\argmin_{\theta^{k;j}_{t}\in \vec{\Theta}_{t}  }|\theta_{t} - \theta^{k;j}_{t}|, &\text{ if } \theta_{t}\in \vec{\Lambda}^{\alpha, \phi, \xi,\rho}_{1,t}(t|n_{t}),\\
\theta_{t}, & \text{ otherwise. }
\end{cases}
$$
Here, $\theta_{t} \sqcup \vec{\Theta}_{t}$ makes $\theta_{t}$ equal to its closest boundary for any $\theta_{t}\in \Theta_{t}$ which is in any stopping sub-region.  
The following corollary shows a sufficient condition for the DIC mechanism when Assumption \ref{assp:single_crossing} does not hold.

\begin{corollary}\label{corollary:sufficient_multiple_stopping_region}
Fix an allocation rule $\alpha$. Suppose that Assumptions \ref{assp:stopping_condition} and \ref{assp:full_support} hold.
The dynamic mechanism is dynamic incentive-compatible if the payment rules $\phi$ and $\xi$ are constructed in (\ref{eq:payment_phi_eta_1}) and (\ref{eq:payment_xi_eta}), respectively, and $\rho$ is constructed as:
\begin{equation}\label{eq:rho_construction_multiple}
    \begin{aligned}
   \rho(t) = \delta^{-t}\mathbb{E}^{\alpha |\mathbb{E}^{\Xi_{\alpha}}[\tilde{\theta}_{t}]\sqcup \vec{\Theta}_{t} }\Big[\sum_{s=t}^{T-1}\big(&\beta^{\alpha}_{S,s+1}(\tilde{\theta}_{s+1} \sqcup \vec{\Theta}_{s+1}) -\beta^{\alpha}_{S,s}(\tilde{\theta}_{s}\sqcup \vec{\Theta}_{s})  \big) \\
   &- \big( \beta^{\alpha}_{\bar{S},s+1}(\tilde{\theta}_{s+1} \sqcup \vec{\Theta}_{s+1}) -\beta^{\alpha}_{\bar{S},s}(\tilde{\theta}_{s}\sqcup \vec{\Theta}_{s}) \big)   \Big],
\end{aligned}
\end{equation}
where the potential functions $\beta^{\alpha}_{S,t}$ and $\beta^{\alpha}_{\bar{S},t}$ satisfy the conditions (\ref{eq:potential_function_relation_S})-(\ref{eq:potential_function_relation_ine}).
\end{corollary}

Given the stopping region $\vec{\Lambda}_{1,t}^{\alpha,\phi,\xi,\rho}$ with $\vec{\Theta}_{t}$, the proof of Corollary \ref{corollary:sufficient_multiple_stopping_region} directly follows Proposition \ref{proposition:sufficient_ICS}.
The construction of $\rho$ in (\ref{eq:rho_construction_multiple}) is equivalent to (\ref{eq:construct_rho_primal}) if the formulations of the potential functions satisfy the monotonicity condition (\ref{eq:monotonicity_potential_assumption_2}) in Proposition \ref{proposition:sufficient_ICS} (thus Assumption \ref{assp:single_crossing} holds).
The following corollary summarizes that the main results shown in Propositions  \ref{prop:envelope_necessary} and  \ref{prop:revenue_equivalence_new} and Corollary \ref{prop:uniqueness_rho} hold for the general optimal stopping rule.

\begin{corollary}\label{corollary:also_for_multiple_regions}
Fix an allocation rule $\alpha$. Suppose that Assumptions \ref{assp:stopping_condition} and \ref{assp:full_support} hold. In DIC mechanism with $\phi$, $\xi$, and $\rho$ constructed in (\ref{eq:payment_phi_eta_1}), (\ref{eq:payment_xi_eta}), and (\ref{eq:rho_construction_multiple}) respectively, (i) the potential functions are formulated in terms of $\alpha$ in (\ref{eq:beta_necessary_S}) and (\ref{eq:beta_necessary_SBar}); (ii) the results of the revenue equivalence shown in Proposition \ref{prop:revenue_equivalence_new} and Corollary \ref{prop:uniqueness_rho} hold. 
\end{corollary}

However, the sufficient condition given by Corollary \ref{corollary:sufficient_multiple_stopping_region} requires the principal to determine $n_{t}$ and all the boundaries of the $n_{t}$ stopping sub-regions.
If the formulations of the potential functions maintain the monotonicity specified by Assumption \ref{assp:single_crossing}, then the principal's mechanism design problem only needs to deal with a unique threshold function $\eta(t)$.

From the construction of $\rho$ in (\ref{eq:construct_rho_primal}) and the formulations of the potential functions in Proposition \ref{prop:envelope_necessary}, $\rho$ can be completely modeled by $\alpha$ and $\eta$, given the transition kernels.
Specifically, at each $t\in\mathbb{T}$, one $\eta(t)\in\Theta_{t}$ leads to $\rho$ that is constructed by (\ref{eq:construct_rho_primal}) and satisfies (\ref{eq:prop_uniqueness_rho}).
Let $r\equiv\{r_{1},\dots,r_{T}\}\in \mathbb{R}^{T}$ be any sequence of real values.
Define
\begin{equation}\label{eq:control_zoon}
    \begin{aligned}
       \mathring{\mathcal{R}} \equiv& \{r\in \mathbb{R}^{T}:  r_{t} = \delta^{-t}\mathbb{E}^{\alpha|\eta(t)}\Big[\sum_{s=t}^{T-1} \big(\beta^{\alpha}_{S,s+1}(\tilde{\theta}_{s+1}\vee \eta(s+1))-\beta^{\alpha}_{S,s}(\tilde{\theta}_{s}\vee \eta(s)) \big) \\
        &- \big(\beta^{\alpha}_{\bar{S},s+1}(\tilde{\theta}_{s+1}\vee \eta(s+1))-\beta^{\alpha}_{\bar{S},s}(\tilde{\theta}_{s}\vee \eta(s))   \big) \Big], \text{ for all }\theta_{t}\in\Theta_{t}, t\in\mathbb{T}\}.
    \end{aligned}
\end{equation}
Hence, $\mathring{\mathcal{R}}$ is the set of payment sequences that $\rho$ can specify from $t=1$ to $t=T$, given $\alpha$ and $\eta$.

\begin{corollary}\label{corollary:no_omnipotent}
Fix an allocation rule $\alpha$. Suppose that Assumptions \ref{assp:stopping_condition}, \ref{assp:full_support}, and \ref{assp_monotone_transistions} hold.
Suppose additionally that $u_{1,t}(\theta_{t}, a_{t})$ is a non-decreasing function of $\theta_{t}$.
Let $\eta$ be a threshold function.
Suppose that $\rho$ specifies a sequence of payments $r=\{r_{t}, \dots, r_{T}\}$, where $r_{t} = \rho(t)$, $t\in\mathbb{T}$.
Then, the following statements hold.
\begin{itemize}
    \item[(i)]  In DIC mechanisms, $\rho$ cannot punish the agent for stopping; i.e., $r_{t}\geq 0$, for all $t\in\mathbb{T}$.
    \item[(ii)] $\rho$ with $\alpha$ implements the optimal stopping rule (\ref{eq:stopping_rule_formulation}) if and only if $r\in \mathring{\mathcal{R}}$.
    \item[(iii)] Fix any $r\in \mathring{\mathcal{R}}$. Let $r'$ differ from $r$ in arbitrary periods, such that for some $t$, $r'_{t} = r_{t} + \epsilon_{t}$, for some non-zero $\epsilon_{t}\in\mathbb{R}$. 
    Suppose $r'\not\in \mathring{\mathcal{R}}$ due to these $r'_{t}$'s.
    Then, posting $r'$ may fail the incentive compatibility constraints.
    \item[(iv)] Let $\phi$ and $\xi$ be constructed in Proposition \ref{proposition:sufficient_ICS}. 
    If the mechanism does not involve $\rho$,
    then $\{\alpha,\phi,\xi\}$ is dynamic incentive-compatible if and only if there exists $\eta(t)\in \Theta_{t}$ such that, for all $t\in\mathbb{T}$,
    \begin{equation}
        \begin{aligned}
        \mathbb{E}^{\alpha|\eta(t)}\Big[\sum_{s=t}^{T-1}& \big(\beta^{\alpha}_{S,s+1}(\tilde{\theta}_{s+1}\vee \eta(s+1))-\beta^{\alpha}_{S,s}(\tilde{\theta}_{s}\vee \eta(s)) \big) \\
        &- \big(\beta^{\alpha}_{\bar{S},s+1}(\tilde{\theta}_{s+1}\vee \eta(s+1))-\beta^{\alpha}_{\bar{S},s}(\tilde{\theta}_{s}\vee \eta(s))   \big) \Big] = 0,
        \end{aligned}
    \end{equation}
    where $\beta^{\alpha}_{\bar{S},t}$ and $\beta^{\alpha}_{S,t}$ are constructed in terms of $\alpha$ in Proposition \ref{prop:envelope_necessary}.
    %
\end{itemize}
\end{corollary}

The following proposition shows a necessary condition for DIC that follows Lemma \ref{lemma:envelope_conditions}.
\begin{proposition}\label{prop:necessary_alpha}
Suppose that Assumptions \ref{assp:stopping_condition} and \ref{assp:full_support} hold.
In any DIC mechanism $\{\alpha, \phi,$  $\xi,\rho\}$, the following conditions hold: for any arbitrary fixed state $\theta_{\epsilon,t}\in \Theta_{t}$, any $\theta_{t}$, $\hat{\theta}_{t}\in\Theta_{t}$, $h^{\theta}_{t-1}\in\prod_{s=1}^{t-1}\Theta_{s}$, $t\in\mathbb{T}$,
%
\begin{equation}\label{eq:potential_function_relation_S_SufNece}
    \int^{\hat{\theta}_{t}}_{\theta_{t}} \gamma^{\alpha}_{t}(t,r|h^{\theta}_{t-1}) dr \leq \int^{\hat{\theta}_{t}}_{\theta_{t}} \frac{\partial u_{1,t}(x, \alpha_{t}(\hat{\theta}_{t}|h^{\theta}_{t-1}))}{\partial x}\Big|_{x=r} dr,
\end{equation}
\begin{equation}\label{eq:potential_function_relation_nS_SufNece}
    \begin{aligned}
   \sup_{\tau\in\mathbb{T}_{t}}\Big\{ \int^{\hat{\theta}_{t}}_{\theta_{\epsilon,t}}\gamma^{\alpha}_{t}(&\tau, r|h^{\theta}_{t-1})dr\Big\}- \sup_{\tau\in\mathbb{T}_{t}}\Big\{ \int^{\theta_{t}}_{\theta_{\epsilon,t}}\gamma^{\alpha}_{t}(\tau, r|h^{\theta}_{t-1})dr\Big\}\\
   \leq& \sup_{\tau\in\mathbb{T}_{t}}\Big\{\int^{\hat{\theta}_{t}}_{\theta_{t}} \mathbb{E}^{\alpha|r, \hat{\theta}_{t}}\Big[\sum_{s=t}^{\tau}\frac{\partial u_{1,s}(x, \alpha_{s}(\tilde{\theta}_{s}|h^{\hat{\theta}_{t}}_{s-1}))}{\partial x}\Big|_{x=\tilde{\theta}_{s}}G_{t,s}(h^{\tilde{\theta}}_{t,s}) dr  \Big]\Big\},
\end{aligned}
\end{equation}
where $\gamma_{t}^{\alpha}(\tau,\theta_{t}|h^{\theta}_{t-1}) \equiv \frac{\partial U^{\alpha,\phi,\xi,\rho}_{t}(\tau,r|h^{\theta}_{t-1}) }{ \partial r }\Big|_{r= \theta_{t}}$.
%
\end{proposition}
\proof
See Appendix \ref{app:thm_sufficient_necessary_allocation}. 
\endproof

Proposition \ref{prop:necessary_alpha} establishes a first-order necessary condition for DIC. This necessary condition takes advantage of the envelope condition established in Lemma \ref{lemma:envelope_conditions} to characterize the necessity of DIC in terms of the allocation rule. Note that the conditions (\ref{eq:potential_function_relation_S_SufNece}) and (\ref{eq:potential_function_relation_nS_SufNece}) hold for both formulations of $\rho$ in Proposition \ref{proposition:sufficient_ICS} and Corollary \ref{corollary:sufficient_multiple_stopping_region}, respectively.

%

\section{Principal's Optimal Mechanism Design }\label{sec:example_optimal}


At the ex-ante stage, the principal provides a take-it-or-leave-it offer by designing the decision rules $\{\alpha,\phi,\xi,\rho\}$ to maximize her ex-ante expected payoff (\ref{eq:ex_ante_payoff_principal}).
The time horizon of the principal's optimization problem is the \textit{mean first passage time}, $\bar{\tau}$.
Given the agent's optimal stopping rule, $\bar{\tau}$ is determined by $\{\alpha,\phi,\xi,\rho\}$ and the stochastic process $\Xi_{\alpha}$.
Propositions \ref{proposition:sufficient_ICS} and \ref{prop:envelope_necessary} imply that $\phi$ and $\xi$ can be represented by $\alpha$. From Corollary \ref{corollary:no_omnipotent}, $\rho$ can be fully characterized by $\alpha$ and $\eta$.
%
%
Then, we can determine $\bar{\tau}$ by $\alpha$, $\eta$, and $\Xi_{\alpha}$, i.e., $\bar{\tau} = \lambda(\alpha,\eta; \Xi_{\alpha})$.
%
%
Specifically,
\begin{equation}\label{eq:mapping_MFPT}
    \begin{split}
        \bar{\tau} =& \lambda(\alpha,\eta; \Xi_{\alpha})\\
        =& \sum_{t=1}^{T} t\cdot P_{r}\big(\tilde{\theta}_{t}\leq \eta(t)\big)\\
        =& \mathbb{E}^{\Xi_{\alpha}}\Big[\sum^{T}_{t=1}t\cdot F_{t}(\eta(t)| \tilde{\theta}_{t-1}, \alpha_{t-1}(\tilde{\theta}_{t-1}|h^{\tilde{\theta}}_{t-2}))\Big].
    \end{split}
\end{equation}
Since the principal's mechanism design problem is finite-horizon, i.e., $T<\infty$, $\bar{\tau}$ always exists.
%
%

Given $\{\alpha,\phi,\xi,\rho\}$, the agent decides whether to participate at the ex-ante stage, by checking the rational participation (RP) constraint, 
\begin{equation}\label{eq:initial_RPC}
    \text{RP:}\;\;\;J^{\alpha,\phi,\xi,\rho}_{1}(\lambda(\alpha,\eta; \Xi_{\alpha}))\geq 0;
\end{equation}
i.e., the ex-ante expected payoff of the agent is non-negative.
Since the decision making at the ex-ante stage involves no private information, the principal's evaluation of $\bar{\tau}$ coincides with the agent's.
Hence, the principal's problem is
\begin{equation}\label{eq:principal_MD_problem}
    \begin{split}
        \max_{\alpha, \phi, \xi,\rho}\;\;& J^{\alpha, \phi, \xi,\rho}_{0}(\lambda(\alpha,\eta; \Xi_{\alpha}))\\
        \text{s.t., }& \text{ DIC},\\
        & \text{ RP.}
    \end{split}
\end{equation}



The mechanism design problem formulated in (\ref{eq:principal_MD_problem}) can be generally applicable to principal-agent models in different economic environments, in which the agent's private information endogenously changes over time.
The mechanism could be a dynamic contract, in which the principal commits to a $T$-period contract and the agent is allowed to unilaterally terminate the contract without advance notification.
As an illustration, consider that the principal is an employer and the agent is an employee.
The employee's state is his attitude towards work, which dynamically changes over time due to the time-evolving development of employee loyalty.
The employer's mechanism design problem is to seek an optimal way to design and allocate inter-temporal work and salary arrangements to facilitate the ongoing development of employee loyalty and to specify compensation or penalty policy to influence the employee's unilateral termination of the contract and hedge the corresponding risk of losing profit.
Other applications include dynamic subscription policies (e.g., magazine, online courses, and customized products), insurance, and leasing. 
Allowing for the early termination provides additional flexibility for the agent and reduces his risk due to the uncertainty of the dynamic environment, thereby increasing the attractiveness of the principal's offer.

Exactly solving the mechanism design problem (\ref{eq:principal_MD_problem}) needs to satisfy the following criteria: ($C1$) game-theoretic constraints, i.e., DIC and RP, ($C2$) profit-maximizing, assuming truthful reporting, and ($C3$) computational constraints (e.g., scale, size, speed, accuracy of approximation, and complexity).
The classic economic analysis of mechanism design focuses on criteria $C1$ and $C2$, based on some distributional assumptions about the participants (e.g., a known probability distribution of the private information, known dynamics, etc.).
Algorithmic mechanism design studies the mechanism design problem by applying approaches, such as robust analysis and algorithmic approximations, and taking into account the computational constraints as important ones, in addition to the classic game-theoretic analysis.
When it is not possible to exactly solve the mechanism design problem (due to that, e.g., the incentive compatibility constraint is too strong), we need to relax the aforementioned three criteria.

Relaxations can be applied to each of the criteria or the hybrid of them.
Mechanism design problems can be analytically relaxed by, for example, imposing additional assumptions and conditions that reduce the strength of $C1$ or relax the optimality of $C2$.
Computational relaxations, usually referred to as approximations, are considered due to the analytical intractability or because $C3$ is not satisfied.
However, relaxations, especially the computational ones, may inevitably lose the robustness of the incentive compatibility to the agent's strategic reporting; i.e., the agent may gain profits by deviating from truthful reporting.
As highlighted in the algorithmic mechanism design, addressing the tension between the game-theoretic constraints and the computational ones is the central challenge when the economic models are used in practice. 

Nevertheless, the game-theoretic guarantee of strong incentive compatibility is important with the theoretical and the practical interest (\cite{pathak2008leveling,lubin2012approximate,azevedo2019strategy}). For example, strong incentive compatibility leads to the simplicity of the agent's reasoning in the decision making and the dynamic stability such that the agent does not need to modify the reporting strategy in response to changes in the environment and other agents' behaviors in the multi-agent cases. Strong incentive compatibility may also provide normative advice in the sense that the system is fair because agents cannot game the system to their advantage.
Moreover, for empirical analysis, reasonably assuming that reporting is truthful enables the empirical work to adjust mechanism parameters or to reshape economic policies without predicting the corresponding complicated strategic responses (i.e., strategic reporting) of the agents.

\subsection{Relaxation}\label{sec:relaxation_1}

It is beyond the scope of this paper to design efficient computational algorithms and analyze computational relaxations.
We restrict attention to analytical relaxations and provide examples of such relaxations in this section.

The principal's mechanism design problem (\ref{eq:principal_MD_problem}) can be relaxed by applying a first-order approach (see, e.g., \cite{rogerson1985first,tadelis2005lectures}) to modify the objective function, which requires the principal to choose decision rules with which her ex-ante expected payoff is at a stationary point; i.e., $\beta^{\alpha}_{S,t}$ and $\beta^{\alpha}_{\bar{S},t}$ satisfy the conditions in Proposition \ref{prop:envelope_necessary} and $\{\phi,\xi,\rho\}$ are constructed according to (\ref{eq:payment_phi_eta_1})-(\ref{eq:construct_rho_primal}).
%
%
%
Specifically, the original ex-ante expected payoff function (\ref{eq:ex_ante_payoff_principal}) can be rewritten as follows:
\begin{equation}\label{eq:relaxed_principal_object_v0}
    \begin{aligned}
         J^{\alpha,\phi,\xi,\rho}_{0}(\bar{\tau})=& - J^{\alpha,\phi,\xi, \rho}_{1,1}(\bar{\tau},\underline{\theta_{1}})+ \mathbb{E}^{\Xi_{\alpha}}\Big[\sum^{1}_{i=0}\sum^{ \bar{\tau}}_{t=1}\delta^{t}u_{i,t}(\tilde{\theta}_{t}, \alpha_{t}(\tilde{\theta}_{t}|h^{\tilde{\theta}}_{t-1})) \Big]  \\
         -& \mathbb{E}^{\Xi_{\alpha} }\Big[\frac{1- F_{1}(\tilde{\theta}_{1})}{f_{1}(\tilde{\theta}_{1})}\sum_{t=1}^{\bar{\tau}} \delta^{t}\frac{  \partial u_{1,t}(r, \alpha_{t}(\tilde{\theta}_{t}|h^{\tilde{\theta}}_{t-1}))}{ \partial r }\big|_{r = \tilde{\theta}_{t}}G_{1,t}(h^{\tilde{\theta}}_{1,t})\Big].
    \end{aligned}
\end{equation}
%






\begin{lemma}\label{lemma:non-decreasing_J}
Suppose that Assumptions \ref{assp:stopping_condition}, \ref{assp:full_support}, and \ref{assp_monotone_transistions} hold.
Suppose additionally that the utility function $u_{1,t}(\theta_{t}, a_{t})$ is non-decreasing in $\theta_{t}$.  
Then, $J^{\alpha,\phi,\xi,\rho}_{1,t}$ is a non-decreasing function of $\theta_{t}$, for all $t\in\mathbb{T}$.
\end{lemma}
\proof
See Appendix \ref{app:proof_lemma:non-decreasing_J}.
\endproof

From Lemma \ref{lemma:non-decreasing_J},
the condition $J^{\alpha,\phi,\xi, \rho}_{1,1}(\bar{\tau},\theta_{1}) \geq 0$ implies the RP constraint (\ref{eq:initial_RPC}) due to $J^{\alpha, \phi,\xi,\rho}_{1}(\bar{\tau}) = \mathbb{E}^{\Xi_{\alpha}}\Big[J^{\alpha, \phi,\xi,\rho}_{1,1}(\bar{\tau}, \tilde{\theta}_{1})  \Big]$.
Therefore, under the assumption that \\$J^{\alpha,\phi,\xi, \rho}_{1,1}(\bar{\tau},\underline{\theta_{1}})=0$, the principal's mechanism design problem can be solved by the following constrained dynamic programming:
\begin{equation}\label{eq:relaxed_principal_object_v1}
    \begin{aligned}
         \max_{\alpha, \eta}\;\;J^{\alpha}_{0}(\lambda(\alpha,\eta; \Xi_{\alpha}))
         =&\mathbb{E}^{\Xi_{\alpha}}\Bigg[\sum_{t=1}^{\lambda(\alpha,\eta; \Xi_{\alpha})}\delta^{t} \Big[ \sum^{1}_{i=0} u_{i,t}(\tilde{\theta}_{t}, \alpha_{t}(\tilde{\theta}_{t}|h^{\tilde{\theta}}_{t-1})) \\
         &- \frac{1- F_{1}(\tilde{\theta}_{1})}{f_{1}(\tilde{\theta}_{1})} \frac{  \partial u_{1,t}(r, \alpha_{t}(\tilde{\theta}_{t}|h^{\tilde{\theta}}_{t-1}))}{ \partial r }\Big|_{r = \tilde{\theta}_{t}}G_{1,t}(h^{\tilde{\theta}}_{1,t})\Big] \Bigg],\\
         \text{s.t., }& \text{ DIC,}\\
         & J^{\alpha,\phi,\xi, \rho}_{1,1}(\bar{\tau},\underline{\theta_{1}}) = 0.
         %
    \end{aligned}
\end{equation}
The relaxed problem (\ref{eq:relaxed_principal_object_v1}) is independent of the payment rules $\{\phi, \xi, \rho\}$.
%
%
%
\begin{corollary}\label{corollary:if_can_be_solved}
Suppose that Assumptions \ref{assp:stopping_condition}, \ref{assp:full_support}, and \ref{assp_monotone_transistions} hold.
Suppose additionally that the utility function $u_{1,t}(\theta_{t}, a_{t})$ is non-decreasing in $\theta_{t}$.  
If the problem (\ref{eq:relaxed_principal_object_v1}) has solutions, then there exist payment rules $\{\phi^{*}, \xi^{*}, \rho^{*}\}$ such that the mechanism decision rules $\{\alpha^{*}, \phi^{*}, \xi^{*}, \rho^{*}\}$ maximize (\ref{eq:principal_MD_problem}).
The resulting mechanism is DIC and satisfies $J^{\alpha,\phi,\xi,\rho}_{1,1}(\bar{\tau},\underline{\theta}_{t})=0$. 
\end{corollary}



However, the constrained dynamic programming (\ref{eq:relaxed_principal_object_v1}) in general does not have closed-form solutions and computational approximations are necessary.
As mentioned earlier in this section, approximations may reduce the robustness of the DIC, i.e., there exist $\epsilon_{S,t}\geq0$ and $\epsilon_{\bar{S},t}\geq 0$, $t\in\mathbb{T}$, such that, for all $\theta_{t}$, $\hat{\theta}_{t}\in\Theta_{t}$, $h^{\theta}_{t-1}\in\prod_{s=1}^{t-1}\Theta_{s}$,
\begin{equation}\label{eq:O-LA-IC_1_epsilon}
    U_{S,t}^{\alpha,\phi, \xi, \rho}(\theta_{t}|h^{\theta}_{t-1})+ \epsilon_{S,t}\geq U_{S,t}^{\alpha,\phi, \xi, \rho}(\theta_{t}, \hat{\theta}_{t}|h^{\theta}_{t-1}),
\end{equation}
and
\begin{equation}\label{eq:O-LA-IC_2_epsilon}
    U_{\bar{S},t}^{\alpha,\phi, \xi, \rho}(\theta_{t}|h^{\theta}_{t-1})+\epsilon_{\bar{S},t}\geq U_{\bar{S},t}^{\alpha,\phi, \xi, \rho}(\theta_{t}, \hat{\theta}_{t}|h^{\theta}_{t-1}).
\end{equation}
The inequalities (\ref{eq:O-LA-IC_1_epsilon}) and (\ref{eq:O-LA-IC_2_epsilon}) imply that the agent cannot improve his payoff by more than $\epsilon_{S,t}$ (resp. $\epsilon_{\bar{S},t}$) through misreporting his true state if he decides to stop at $t$ (resp. to continue at $t$). We say that such mechanism with $\{\alpha,\phi,\xi,\rho\}$ is $\big\{\epsilon_{S,t}, \epsilon_{\bar{S},t}\big\}$-DIC.
Define the \textit{deviations}, for $t\in\mathbb{T}$,
\begin{equation}\label{eq:relax_distance_stop}
    \begin{aligned}
       d^{\alpha}_{S,t} = \sup_{\theta_{t},\hat{\theta}_{t}\in \Theta_{t}}\Big\{ u_{1,t}(\theta_{t},\alpha_{t}(\hat{\theta}_{t}|h^{\theta}_{t-1})) - u_{1,t}(\hat{\theta}_{t},\alpha_{t}(\hat{\theta}_{t}|h^{\theta}_{t-1}))+
       \beta^{\alpha}_{S,t}(\hat{\theta}_{t}) - \beta^{\alpha}_{S,t}(\theta_{t}) \Big\},
    \end{aligned}
\end{equation}
and
\begin{equation}\label{eq:relax_distance}
    \begin{aligned}
       d^{\alpha}_{\bar{S},t}= \sup_{\theta_{t},\hat{\theta}_{t}\in \Theta_{t}}\Bigg\{\sup_{\tau\in\mathbb{T}_{t}}\Big\{ \pi^{\alpha}_{t}(\theta_{t},\hat{\theta}_{t};\tau) \Big\}-&\sup_{\tau\in\mathbb{T}_{t}}\Big\{ \pi^{\alpha}_{t}(\hat{\theta}_{t};\tau) \Big\} +\beta^{\alpha}_{\bar{S},t}(\hat{\theta}_{t})- \beta^{\alpha}_{\bar{S},t}(\theta_{t})\Bigg\},
    \end{aligned}
\end{equation}
where $\pi^{\alpha}_{t}$ 
is given in (\ref{eq:sub_path_length}). 
%
\begin{proposition}\label{prop:relax_approximation}
Suppose that Assumptions \ref{assp:stopping_condition}, \ref{assp:full_support}, and \ref{assp_monotone_transistions} hold.
Let $\alpha^{\circ}$ be an approximated optimal allocation rule as a solution to (\ref{eq:relaxed_principal_object_v1}). 
Let $\beta^{\alpha^{\circ}}_{\bar{S},t}$ and $\beta^{\alpha^{\circ}}_{S,t}$, respectively, be constructed according to (\ref{eq:beta_necessary_S}) and (\ref{eq:beta_necessary_SBar}).
If the payment rules $\phi^{\circ}$ and $\xi^{\circ}$ are constructed according to (\ref{eq:payment_phi_eta_1}) and (\ref{eq:payment_xi_eta}), respectively, and $\rho^{\circ}$ is constructed according to (\ref{eq:construct_rho_primal}) (with an additional assumption that $u_{1,t}(\theta_{t},a_{t})$ is non-decreasing in $\theta_{t}$) or (\ref{eq:rho_construction_multiple}), then 
\begin{itemize}
    \item[(i)] the mechanism is $\Big\{d^{\alpha^{\circ}}_{S,t}, d^{\alpha^{\circ}}_{\bar{S},t}+\sup_{\tau\in\mathbb{T}_{t}}\Big\{\rho^{\circ}(\tau) \Big\}\Big\}$-DIC, when $d^{\alpha^{\circ }}_{S,t}>0$ and\\ $d^{\alpha^{\circ}}_{\bar{S},t}> -\sup_{\tau\in\mathbb{T}_{t}}(\rho^{\circ}(\tau))$;
    \item[(ii)] the mechanism is DIC, when $d^{\alpha^{\circ}}_{S,t}\leq0$ and $d^{\alpha^{\circ}}_{\bar{S},t}\leq -\sup_{\tau\in\mathbb{T}_{t}}(\rho^{\circ}(\tau))$.
\end{itemize}

\end{proposition}
\proof
See Appendix \ref{app:prop_relax_approximation}.
\endproof

Proposition \ref{prop:relax_approximation} provides an approach to evaluate the worst-case scenario of the agent's strategic misreporting; i.e., the most profitable deviations from truthfulness, when the allocation rule $\alpha^{\circ}$ is an approximate solution of the relaxed problem (\ref{eq:relaxed_principal_object_v1}).
If the payment rules are constructed according to (\ref{eq:payment_phi_eta_1})-(\ref{eq:construct_rho_primal}) in Proposition \ref{proposition:sufficient_ICS}, then the mechanism is $\Big\{d^{\alpha^{\circ}}_{S,t}, d^{\alpha^{\circ}}_{\bar{S},t}+\sup_{\tau\in\mathbb{T}_{t}}\Big\{\rho^{\circ}(\tau) \Big\}  \Big\}$-DIC. 
The evaluation approach also provides a sufficient condition for DIC, i.e., if $d^{\alpha^{\circ}}_{S,t}\leq0$ and $d^{\alpha^{\circ}}_{\bar{S},t}\leq -\sup_{\tau\in\mathbb{T}_{t}}(\rho^{\circ}(\tau))$, then the mechanism is DIC.

\subsection{Case Study}\label{sec:example_relaxed}

We consider a dynamic resource allocation problem to illustrate the verification of the theoretical conditions shown in Propositions \ref{proposition:sufficient_ICS}, \ref{prop:envelope_necessary}, and \ref{prop:necessary_alpha}.
Note that this example is crafted as analytically tractable for the purpose of illustration. 
Consider that the principal allocates resource $a_{t}\in A_{t}=[0,\bar{a}_{t}]$ to the agent based on his state $\theta_{t}$ at each time $t\in\mathbb{T}$.
Suppose that $\Theta_{t}=[\underline{\theta}, \bar{\theta}]\subset \mathbb{R}_{+}$, for all $t\in\mathbb{T}$, the dynamics of the agent's state follow a nonlinear autoregressive model: (see, e.g., \cite{blasques2020nonlinear}): 
for $t\in\mathbb{T}\backslash{\{T\}}$, 
\begin{equation}\label{eq:example_dynamics}
    \theta_{t+1} = \zeta(\theta_{t}) + b_{t+1},
\end{equation}
where $\zeta(\theta_{t})$ is a nonlinear, non-decreasing, and Lipschitz continuous function of $\theta_{t}$ with bounded derivative and $\min_{\theta_{t}\in\Theta_{t}}\zeta(\theta_{t})=\underline{\theta}$ and $\max_{\theta_{t}\in\Theta_{t}}\zeta(\theta_{t})<\bar{\theta}$, and $b_{t}$ is distributed over $B=[0, \bar{\theta} - \zeta(\bar{\theta})]$.
Since $\zeta$ is non-decreasing in $\theta_{t}$, the condition in Assumption \ref{assp_monotone_transistions} is satisfied.
Suppose that $\zeta$ and the distribution of $b_{t}$ are set such that Assumption \ref{assp:full_support} holds.
Suppose additionally that the initial kernel has $F_{1}$ and $f_{1}$ such that $\frac{1-F_{1}(\theta_{1}) }{f_{1}(\theta_{1})}$ is non-increasing in $\theta_{1}$.
According to (\ref{eq:example_dynamics}), $\theta_{k}$, for $k\in\mathbb{T}_{t}$, can be represented in terms of $\theta_{t}$; i.e., there exists some function $\bar{\zeta}_{t,k}:\Theta_{t}\times  B^{k-t}\mapsto \Theta_{k}$, such that $\theta_{k}=\bar{\zeta}_{t,k}(\theta_{t};h^{b}_{t+1,k})$.
The utility functions of the agent and the principal are given as follows, respectively:
\begin{equation}\label{eq:case_utility_agent}
    u_{1,t}(\theta_{t},a_{t}) = \kappa(\theta_{t})a_{t}+c_{1},
\end{equation}
where $\kappa$ is an increasing linear function of $\theta_{t}$ and $c_{1}\in\mathbb{R}_{+}$, and
\begin{equation}\label{eq:case_utility_principal}
    u_{0,1}(\theta_{t}, a_{t}) = c_{2}a^{2}_{t}+c_{3},
\end{equation}
where $c_{2}\in\mathbb{R}_{-}$ and $c_{3}\in \mathbb{R}_{+}$ are two non-zero real numbers.
It is straightforward to see that the agent's utility function $u_{1,t}$ is Lipschitz continuous.
We assume that both the principal's and the agent's utility functions are bounded, i.e., there exist constants $k<\infty, k'<\infty$, $\{g_{t}\}_{t\in \mathbb{T}}$ with $|g_{t}|<\infty$, and $\{g'_{t}\}_{t\in \mathbb{T}}$ with $|g'_{t}|<\infty$, such that $|u_{0,t}(\theta_{t}, a_{t})|\leq k|\theta_{t}|+ g_{t}$ and $|u_{1,t}(\theta_{t}, a_{t})|\leq k'|\theta_{t}|+ g'_{t}$.
%
%
By applying the partial derivatives in Lemma \ref{lemma:envelope_conditions} to the agent's ex-ante and interim expected payoff functions, we obtain that the boundedness condition in Assumption \ref{assp:stopping_condition} is satisfied.
%

Since Assumptions \ref{assp:stopping_condition}, \ref{assp:full_support}, and  \ref{assp_monotone_transistions} hold, Lemma \ref{lemma:non-decreasing_J} yields that the principal's mechanism design problem can be reduced to (\ref{eq:relaxed_principal_object_v1}) which takes the following form:
\begin{equation}\label{eq:example_objective}
    \begin{split}
        \max_{\alpha,\eta}\;J^{\alpha}_{0}(\lambda(\alpha,\eta;\Xi_{\alpha})) =& \mathbb{E}^{\Xi_{\alpha}}\Big[\sum^{\lambda(\alpha,\eta;\Xi_{\alpha})}_{t=1}\delta^{t}\big[\alpha_{t}(\tilde{\theta}_{t})\big(\kappa(\tilde{\theta}_{t})+c_{2}\alpha_{t}(\tilde{\theta}_{t})\big) +c_{1}+c_{3}  \\
        &- \frac{1-F_{1}(\tilde{\theta}_{1})}{f_{1}(\tilde{\theta}_{1})} \frac{d \kappa(x)}{dx}\Big|_{x= \tilde{\theta}_{t}}\alpha_{t}(\tilde{\theta}_{t}) \prod_{s=1}^{t}\frac{d \zeta(x)}{dx}\Big|_{x= \tilde{\theta}_{s}} \big] \Big],\\
        \text{s.t. } & \text{ DIC,}\\
        & J^{\alpha,\phi,\xi, \rho}_{1,1}(\bar{\tau},\underline{\theta }) = 0.
    \end{split}
\end{equation}

Instead of solving the constrained dynamic programming (\ref{eq:example_objective}), we follow a short-cut procedure (see, e.g., \cite{tadelis2005lectures}) by first ignoring the constraints.
Hence, the allocation rule $\alpha^{*}_{t}$ that maximizes (\ref{eq:example_objective}) (without constraints) is
\begin{equation}\label{eq:example_allocation}
    \begin{split}
        &\alpha^{*}_{t}(\theta_{t}|h^{\theta}_{t-1}) \\
        =& \max\Big\{0, \;\; \frac{1}{2c_{2}}\Big[\frac{1-F_{1}(\theta_{1})}{f_{1}(\theta_{1})}\frac{d \kappa(x)}{dx}\Big|_{x=\theta_{t}}\prod^{t}_{s=1}\frac{d \zeta(x)}{dx}\Big|_{x=\theta_{s}} -\kappa(\theta_{t})\Big]\Big\}.
    \end{split}
\end{equation}
Based on the given environment (i.e., with all non-designed variables), $\alpha^{*}_{t}(\theta_{t})$ is non-decreasing in $\theta_{t}$.
%
%
Greedy algorithms can be used to choose $\bar{\tau}^{*}$ such that the pair $\{\alpha^{*}, \tau^{*}\}$ leads to the maximum value $J^{\alpha^{*}}_{0}(\bar{\tau}^{*})$.
Then, the corresponding threshold $\eta$ can be determined or approximated according to (\ref{eq:mapping_MFPT}). 
%

%
From Proposition \ref{prop:envelope_necessary}, the potential functions are constructed as follows: 
%
\begin{equation}\label{eq:example_potential_stop_1}
    \begin{aligned}
   \beta^{\alpha^{*}}_{S,t}(\theta_{t}) = \int^{\theta_{t}}_{\underline{\theta} }\delta^{t}\alpha^{*}_{t}(r|h^{\theta}_{t-1}) \frac{d \kappa(x)}{d x} \Big|_{x=r}dr,
\end{aligned}
\end{equation}
%
and
\begin{equation}\label{eq:example_potential_nonstop_1}
    \begin{aligned}
   \beta^{\alpha^{*}}_{\bar{S},t}(\theta_{t}) = \sup_{\tau\in\mathbb{T}_{t}}\Big\{\int^{\theta_{t}}_{\underline{\theta} }\mathbb{E}^{\alpha^{*}|r}\Big[\sum_{s=t}^{\tau}\delta^{s}  \alpha^{*}_{s}(\tilde{\theta}_{s}|h^{\tilde{\theta}}_{s-1})\frac{d \kappa(x)}{dx}\Big|_{x=\tilde{\theta}_{s}}\prod^{s}_{k=t+1}\frac{d \zeta(x)}{dx}\Big|_{x=\tilde{\theta}_{k-1}}  \Big] dr \Big\}.
\end{aligned}
\end{equation}
For simplicity, let, for $\tau\in\mathbb{T}_{t}$, $t\in\mathbb{T}$,
\begin{equation}\label{eq:example_D}
    D^{\alpha^{*}}_{t,\tau}(\theta_{t}) =\mathbb{E}^{\alpha^{*}|\theta_{t}}\Big[ \sum_{s=t}^{\tau}\delta^{s}  \alpha^{*}_{s}(\tilde{\theta}_{s}|h^{\tilde{\theta}}_{s-1})\frac{d \kappa(x)}{dx}\Big|_{x=\tilde{\theta}_{s}}\prod^{s}_{k=t+1}\frac{d \zeta(x)}{dx}\Big|_{x=\tilde{\theta}_{k-1}}\Big].
\end{equation}
It is straightforward to verify that $D^{\alpha^{*}}_{t,\tau}(\theta_{t})\geq0$, for all $\theta_{t}\in [\underline{\theta}, \bar{\theta}]$, $\tau\in\mathbb{T}_{t}$, $t\in\mathbb{T}$.
Hence, the R.H.S. of (\ref{eq:example_potential_nonstop_1}) attains the maximum when $\tau=T$.
%
Since $u_{1,t}$ is increasing in $\theta_{t}$ and Assumptions \ref{assp:stopping_condition}, \ref{assp:full_support}, and  \ref{assp_monotone_transistions} hold, from Proposition \ref{prop:envelope_necessary}, we can construct the payment rules according to (\ref{eq:payment_phi_eta_1}), (\ref{eq:payment_xi_eta}), and (\ref{eq:construct_rho_primal}) in Proposition \ref{proposition:sufficient_ICS}, respectively, as follows:
%
%
\begin{equation}\label{eq:example_phi}
    \begin{aligned}
       \phi_{t}&(\theta_{t})\\ =&\sup_{\tau\in\mathbb{T}_{t}}\Big\{\int^{\theta_{t}}_{\underline{\theta} }\mathbb{E}^{\alpha^{*}|r}\Big[\sum_{s=t}^{\tau}\delta^{s}  \alpha^{*}_{s}(\tilde{\theta}_{s}|h^{\tilde{\theta}}_{s-1})\frac{d \kappa(x)}{dx}\Big|_{x=\tilde{\theta}_{s}}\prod^{s}_{k=t+1}\frac{d \zeta(x)}{dx}\Big|_{x=\tilde{\theta}_{k-1}}  \Big] dr \Big\}\\
       &- \mathbb{E}^{\alpha^{*}|\theta_{t}}\Big[ \sup_{\tau\in\mathbb{T}_{t}}\Big\{\int^{\tilde{\theta}_{t+1}}_{\underline{\theta} } \sum_{s=t+1}^{\tau}\delta^{s}  \alpha^{*}_{s}(\tilde{\theta}_{s}|h^{\tilde{\theta}}_{s-1})\frac{d \kappa(x)}{dx}\Big|_{x=\tilde{\theta}_{s}}\prod^{s}_{k=t+2}\frac{d \zeta(x)}{dx}\Big|_{x=\tilde{\theta}_{k-1}} dr  \Big\} \Big]\\
       &- u_{1,t}(\theta_{t}, \alpha^{*}_{t}(\theta_{t})),
    \end{aligned}
\end{equation}

\begin{equation}\label{eq:example_xi}
    \begin{aligned}
   \xi_{t}(\theta_{t}) = \int^{\theta_{t}}_{\underline{\theta} }\alpha^{*}_{t}(r|h^{\theta}_{t-1}) \frac{d \kappa(x)}{dx}\Big|_{x=r}dr - u_{1,t}(\theta_{t}, \alpha^{*}_{t}(\theta_{t})),
\end{aligned}
\end{equation}
and
\begin{equation}\label{eq:example_rho}
    \begin{aligned}
    \rho(t) = \sup_{\tau\in\mathbb{T}_{t}}\Big\{\int^{\eta^{*}(t)}_{\underline{\theta} }&\mathbb{E}^{\alpha^{*}|r}\Big[\sum_{s=t+1}^{\tau}\delta^{s}  \alpha^{*}_{s}(\tilde{\theta}_{s}\vee\eta^{*}(s)|h^{\tilde{\theta}}_{s-1})\frac{d \kappa(x)}{dx}\Big|_{x=\tilde{\theta}_{s}\vee \eta^{*}(s)}\\
   &\cdot\prod^{s}_{k=t+1}\frac{d \zeta(x)}{dx}\Big|_{x=\tilde{\theta}_{k-1}\vee \eta^{*}(k-1)}  \Big] dr \Big\}.              
\end{aligned}
\end{equation}

The length function (\ref{eq:path_length_eta_stop}) is given as, for all $t\in\mathbb{T}$,
\begin{equation}\label{eq:example_path_length_stop}
    \begin{split}
        \ell^{\alpha^{*}}_{S,t}(\hat{\theta}_{t}, \theta_{t})=\delta^{t}\big[\kappa(\hat{\theta}_{t})\alpha^{*}_{t}(\hat{\theta}_{t}|h^{\theta}_{t-1}) -  \kappa(\theta_{t})\alpha^{*}(\hat{\theta}_{t}|h^{\theta}_{t-1})  \big].
    \end{split}
\end{equation}
%
%
%
%
From (\ref{eq:example_path_length_stop}), we have
\begin{equation}\label{eq:example_check_stop}
    \begin{aligned}
   \ell^{\alpha^{*}}_{S,t} (\hat{\theta}_{t}, \theta_{t})
   \geq & \delta^{t}\int^{\hat{\theta}_{t}}_{\theta_{t}}  \frac{d \kappa(x)}{dx}\Big|_{x=r}\alpha^{*}_{t}(r|h^{\theta}_{t-1})dr\\
   =& \delta^{t}\int^{\hat{\theta}_{t}}_{\underline{\theta} }  \frac{d \kappa(x)}{dx}\Big|_{x=r}\alpha^{*}_{t}(r|h^{\theta}_{t-1})dr -\delta^{t}\int^{\theta_{t}}_{\underline{\theta} }  \frac{d \kappa(x)}{dx}\Big|_{x=r}\alpha^{*}_{t}(r|h^{\theta}_{t-1})dr\\
   =&\beta^{\alpha^{*}}_{S,t}(\hat{\theta}_{t})-\beta^{\alpha^{*}}_{S,t}(\theta_{t}).
\end{aligned}
\end{equation}
%
Hence, the allocation rule $\alpha^{*}_{t}$ in (\ref{eq:example_allocation}) satisfies condition (\ref{eq:potential_function_relation_S}) in Proposition \ref{proposition:sufficient_ICS} and condition (\ref{eq:potential_function_relation_S_SufNece}) in Proposition \ref{prop:necessary_alpha}.

Next, we check condition (\ref{eq:potential_function_relation_nS_SufNece}) in Proposition \ref{prop:necessary_alpha}.
Let, for $\tau\in\mathbb{T}_{t}$, $t\in\mathbb{T}$,
\begin{equation}
    \begin{aligned}
       M^{\alpha^{*}}_{t}(\theta_{t},\hat{\theta}_{t};\tau) 
       = & \mathbb{E}^{\alpha^{*}|\theta_{t},\hat{\theta}_{t}}\Big[\delta^{t}\kappa(\theta_{t})\alpha^{*}_{t}(\hat{\theta}_{t}|h^{\theta}_{t-1})+\int^{\tilde{\theta}_{t+1}}_{\underline{\theta}}\mathbb{E}^{\alpha|r}\Big[ D^{\alpha^{*}}_{t+1,T}(r)\Big]dr\\ 
       &-\int^{\tilde{\theta}_{\tau}}_{\underline{\theta}}\mathbb{E}^{\alpha|r}\Big[ D^{\alpha^{*}}_{\tau, T}(r)\Big]dr
       +\delta^{\tau}\big[ \int^{\tilde{\theta}_{\tau}}_{\underline{\theta} }D^{\alpha^{*}}_{\tau,\tau}(r)dr \big]   \Big],
    \end{aligned}
\end{equation}
and
\begin{equation}
    \begin{aligned}
       M^{\alpha^{*}}_{t}(\theta_{t}, \hat{\theta}_{t};t)=& \int^{\theta_{t}}_{\underline{\theta}}\delta^{t}\alpha^{*}_{t}(\hat{\theta}_{t}|h^{\theta}_{t-1})\frac{d\kappa(x)}{dx}\Big|_{x=r}dr.
    \end{aligned}
\end{equation}
Then, the length function (\ref{eq:path_length_eta_continue}) is given as, for $\tau\in\mathbb{T}_{t}$, $t\in\mathbb{T}$,
%
\begin{equation}\label{eq:example_path_length_continue}
    \begin{split}
        \ell^{\alpha^{*}}_{\bar{S},t}(\hat{\theta}_{t},\theta_{t};\tau)=M^{\alpha^{*}}_{t}(\hat{\theta}_{t},\hat{\theta}_{t};&\tau)- M^{\alpha^{*}}_{t}(\theta_{t},\hat{\theta}_{t};\tau)\\
        =&\int^{\hat{\theta}_{t}}_{\theta_{t}}\Big[\frac{d}{dx}M^{\alpha^{*}}_{t}(x,\hat{\theta}_{t};\tau)\Big|_{x=r}\Big]dr.
    \end{split}
\end{equation}
Let $\hat{\theta}_{t}\geq \theta_{t}\in [\underline{\theta}, \bar{\theta}]$.
Due to the monotonicity of $\kappa$ and $\alpha^{*}$ and $D^{\alpha^{*}}_{t,s}(\theta_{t})\geq0$, for all $\theta_{t}\in[\underline{\theta}, \bar{\theta}]$, $s\in\mathbb{T}_{t}$, $t\in\mathbb{T}$, we have, from (\ref{eq:example_path_length_continue}), 
%
\begin{equation}\label{eq:example_path_nonstop_2}
    \begin{split}
        \sup_{\tau\in\mathbb{T}_{t}}\Big\{ \ell^{\alpha^{*}}_{\bar{S},t}(\hat{\theta}_{t},&\theta_{t};\tau) \Big\}=\int^{\hat{\theta}_{t}}_{\theta_{t}}\Big[\frac{d}{dx}M^{\alpha^{*}}_{t}(x,\hat{\theta}_{t};T)\Big|_{x=r}\Big]dr.\\
        \geq& \int^{\hat{\theta}_{t}}_{\theta_{t}}\mathbb{E}^{\alpha^{*}|r, \hat{\theta}_{t}}\Big[\sum^{T}_{s=t}\delta^{s} \frac{d\kappa(x)}{dx}\Big|_{x=\tilde{\theta}_{s}}\alpha^{*}_{s}(\tilde{\theta}_{s}|h^{\tilde{\theta}}_{s-1})  \prod^{s}_{k=t+1}\frac{d \zeta(x)}{dx}\Big|_{x=\tilde{\theta}_{k-1}}\Big] dr.
    \end{split}
\end{equation}
From the construction of the potential function $\beta^{\alpha^{*}}_{\bar{S},t}$ in (\ref{eq:example_potential_nonstop_1}), we have
\begin{equation}
    \begin{split}
        &\text{R.H.S. of (\ref{eq:example_path_nonstop_2})}\\
        \geq&\sup_{\tau\in\mathbb{T}_{t}}\Big\{\int^{\hat{\theta}_{t}}_{\underline{\theta} }\mathbb{E}^{\alpha^{*}|r, \hat{\theta}_{t}}\Big[\sum^{\tau}_{s=t}\delta^{s} \frac{d\kappa(x)}{dx}\Big|_{x=\tilde{\theta}_{s}}\alpha^{*}_{s}(\tilde{\theta}_{s}|h^{\tilde{\theta}}_{s-1})  \prod^{s}_{k=t+1}\frac{d \zeta(x)}{dx}\Big|_{x=\tilde{\theta}_{k-1}}\Big] dr\Big\}\\
        &-\sup_{\tau\in\mathbb{T}_{t}}\Big\{\int^{\theta_{t}}_{\underline{\theta} }\mathbb{E}^{\alpha^{*}|r, \hat{\theta}_{t}}\Big[\sum^{\tau}_{s=t}\delta^{s} \frac{d\kappa(x)}{dx}\Big|_{x=\tilde{\theta}_{s}}\alpha^{*}_{s}(\tilde{\theta}_{s}|h^{\tilde{\theta}}_{s-1})  \prod^{s}_{k=t+1}\frac{d \zeta(x)}{dx}\Big|_{x=\tilde{\theta}_{k-1}}\Big] dr\Big\}\\
        =& \beta^{\alpha^{*}}_{\bar{S}, t}(\hat{\theta}_{t}) - \beta^{\alpha^{*}}_{\bar{S}, t}(\theta_{t}). 
    \end{split}
\end{equation}
Let $\hat{\theta}_{t}\leq \theta_{t}\in [\underline{\theta}, \bar{\theta}]$. Then, 
$$
\begin{aligned}
   \sup_{\tau\in\mathbb{T}_{t}}\Big\{\ell^{\alpha^{*}}_{\bar{S}}(\hat{\theta}_{t}, \theta_{t};\tau)\Big\} = & \int^{\hat{\theta}_{t}}_{\theta_{t}}\alpha^{*}_{t}(\hat{\theta}_{t}|h^{\theta}_{t-1})\frac{d\kappa(x)}{dx}\Big|_{x=r}dr. 
\end{aligned}
$$ 
Since $\alpha^{*}_{t}$ is non-decreasing,  $\int^{\hat{\theta}_{t}}_{\theta_{t}}\alpha^{*}_{t}(\hat{\theta}_{t}|h^{\theta}_{t-1})\frac{d\kappa(x)}{dx}\Big|_{x=r}dr\geq \int^{\hat{\theta}_{t}}_{\theta_{t}}\alpha^{*}_{t}(r|h^{\theta}_{t-1})\frac{d\kappa(x)}{dx}\Big|_{x=r}dr$.
From the fact that $0\leq D^{\alpha^{*}}_{t,t}(\theta'_{t})\leq D^{\alpha^{*}}_{t,\tau}(\theta'_{t})$  for any $\theta'_{t}\in [\underline{\theta}, \bar{\theta}]$, $\tau\in\mathbb{T}_{t}$, we have $\int^{\hat{\theta}_{t}}_{{\theta}_{t}}D^{\alpha^{*}}_{t,t}(r)dr\geq \int^{\hat{\theta}_{t}}_{{\theta}_{t}}D^{\alpha^{*}}_{t,\tau}(r)dr$, for all $\tau\in\mathbb{T}_{t}$, $t\in\mathbb{T}$.
Then, it is straightforward to see that condition (\ref{eq:potential_function_relation_nS_SufNece}) is satisfied when $\hat{\theta}_{t}\leq \theta_{t}$, given the construction of $\beta^{\alpha^{*}}_{\bar{S}}$ in (\ref{eq:example_potential_nonstop_1}). 
Hence, the allocation rule $\alpha^{*}_{t}$ in (\ref{eq:example_allocation}) satisfies condition (\ref{eq:potential_function_relation_nS_SufNece}) in Proposition \ref{prop:envelope_necessary}. 

Finally, we use the evaluation approach in Proposition \ref{prop:relax_approximation}.
From (\ref{eq:example_check_stop}), the deviation $d^{\alpha^{*}}_{S,t}=0$.
The deviation $d^{\alpha^{*}}_{\bar{S},t}$ in (\ref{eq:relax_distance}) can be written as
\begin{equation}\label{eq:case_study_deviation_nS}
    \begin{aligned}
   d^{\alpha^{*}}_{\bar{S},t} =&\sup_{\theta_{t},\hat{\theta}_{t}\in\Theta_{t}}\Bigg\{ \sup_{\tau\in\mathbb{T}_{t}}\Big\{M^{\alpha^{*}}_{t}(\theta_{t},\hat{\theta}_{t};\tau) + \rho(\tau) \Big\}- \sup_{\tau\in\mathbb{T}_{t}}\Big\{M^{\alpha^{*}}_{t}(\hat{\theta}_{t},\hat{\theta}_{t};\tau)+ \rho(\tau) \Big\} \\
   &+ \beta^{\alpha^{*}}_{\bar{S},t}(\hat{\theta}_{t}) - \beta^{\alpha^{*}}_{\bar{S},t}(\theta_{t})\Bigg\}\\
=& \sup_{\theta_{t},\hat{\theta}_{t}\in\Theta_{t}}\Big\{\int^{\hat{\theta}_{t}}_{\theta_{t}}\mathbb{E}^{\alpha^{*}|r}\Big[D^{\alpha^{*}}_{t,T}(r) -\delta^{t}\alpha^{*}_{t}(\hat{\theta}_{t}|h^{\theta}_{t-1})\frac{d\kappa(x)}{dx}\Big|_{x=r}dr \Big] \Big\}\\
=& \int^{\bar{\theta}}_{\underline{\theta}}\mathbb{E}^{\alpha^{*}|r}\Big[D^{\alpha^{*}}_{t,T}(r) -\delta^{t}\alpha^{*}_{t}(\bar{\theta}|h^{\theta}_{t-1})\frac{d\kappa(x)}{dx}\Big|_{x=r}dr\Big].
\end{aligned} 
\end{equation}
The following corollary summarizes the results of this case study by using the results in Propositions \ref{proposition:sufficient_ICS}, \ref{prop:envelope_necessary}, and \ref{prop:relax_approximation}.

\begin{corollary}\label{corollary:case_study}
Consider the state dynamics (\ref{eq:example_dynamics}). 
The agent optimizes his interim expected payoff at each period with the utility function (\ref{eq:case_utility_agent}). The principal optimizes her ex-ante expected payoff with the utility function (\ref{eq:case_utility_principal}).
Suppose that $\{\alpha^{*},\eta^{*}\}$ maximizes (\ref{eq:example_objective}) by ignoring the constraints.
Let $\phi$, $\xi$, and $\rho$ be constructed by $\beta^{\alpha^{*}}_{S,t}$ and $\beta^{\alpha^{*}}_{\bar{S},t}$ given in (\ref{eq:example_potential_stop_1}) and (\ref{eq:example_potential_nonstop_1}), respectively.
Then, the mechanism with $\{\alpha^{*}, \phi,\xi,\rho\}$ admits $\big\{0, d^{\alpha^{*}}_{\bar{S},t}+\rho(t)\big\}$-DIC, where $d^{\alpha^{*}}_{\bar{S},t}$ is given in (\ref{eq:case_study_deviation_nS}).
The mechanism is DIC if, for any $t\in\mathbb{T}$,
\begin{equation}\label{eq:example_corollary}
    \begin{aligned}
   \alpha^{*}&(\bar{\theta})\geq \frac{1}{\delta^{t}[\kappa(\bar{\theta})-\kappa(\underline{\theta})]}\Bigg[ \rho(t) + \int^{\bar{\theta}}_{\underline{\theta}} \mathbb{E}^{\alpha^{*}}\Big[D^{\alpha^{*}}_{t,T}(r)\Big]dr \Bigg],
\end{aligned}
\end{equation}
where $D^{\alpha^{*}}_{t,T}$ and $\rho(t)$ are given in (\ref{eq:example_D}) and (\ref{eq:example_rho}), respectively.
%
%
%
\end{corollary}

As stated in Corollary \ref{corollary:case_study}, the mechanism with decision rules that solves (\ref{eq:example_objective}) is in general not DIC. Condition (\ref{eq:example_corollary}) is directly from the statement $(ii)$ of Proposition \ref{prop:relax_approximation}.
The mechanism is DIC if the non-designed components of the model (i.e., $\kappa$, $\bar{\theta}$, $\underline{\theta}$, and $\xi$) satisfy condition (\ref{eq:example_corollary}).

\section{Conclusion}\label{sec:conclusion}

This work focuses on the theoretical characterizations of the incentive compatibility of a finite-horizon dynamic mechanism design problem, in which the agent has the right to stop the mechanism at any period.
We have studied an optimal stopping time problem under this dynamic environment for the agent to optimally select the time of stopping.
A state-independent payment rule has been introduced that delivers a payment only at the realized stopping time.
This payment rule enables the principal to directly influence the realization of the agent's stopping time.
We have also shown that under certain conditions, the optimal stopping problem can be fully represented by a threshold rule.
Dynamic incentive compatibility has been defined in terms of the Bellman equations.
A one-shot deviation principle has been established to address the complexity from the dynamic nature of the environment and the coupling of the agent's reporting choices and stopping decisions.
By relying on a set of formulations characterized by the non-monetary allocation rule and the potential functions, we have constructed the payment rules to obtain the sufficiency argument of the dynamic incentive compatibility.
The quasilinear payoff formulation enables us to
derive a necessary condition for the dynamic incentive compatibility from the envelope theorem which determines the explicit formulation of the potential functions.
These settings naturally lead to the revenue equivalence.

Our analysis provides a new design paradigm for optimal direct mechanism in general quasilinear dynamic environments when the dynamic incentive compatibility takes into account not only the agent's reporting behaviors but also his stopping decisions.
From the necessary and the sufficient conditions, we can design the state-independent terminal payment rule by the allocation rule and the threshold function.
As a result, the expected first-passage time (i.e., the expected time-horizon of the principal's ex-ante expected payoff) seen at the ex-ante stage is fully determined by the allocation rule and the threshold function given the transition kernels of the state.
We have described the principal's optimal mechanism design by applying a relaxation approach that reformulates the principal's optimal mechanism design by making the principal make decisions at a stationary point. 
Hence, the principal's optimization problem can be handled by finding the optimal allocation rule and the optimal threshold function.
A regular condition has been provided as a design principle for the state-independent payment rule.
An evaluation approach has been provided to evaluate the loss of the robustness of the dynamic incentive compatibility due to relaxations or computational approximations.
In a case study, we have shown an example of a relaxed mechanism design and used the evaluation approach to obtain an approximate dynamic incentive-compatible mechanism.

The extension to multiple-agent environments would be a natural next step. 
In environments with multiple agents, allowing the early exit of each agent leads to a dynamic population over time, which is state-dependent.
The relationships between population and individual payoffs could complicate the analysis of the incentive compatibility.
One nontrivial extension could introduce the arrivals of new agents whose incentives of participation are characterized by both dynamic rational participation constraints as well as the history of the mechanism.
Involving renegotiation after the agent realizes the stopping time is another direction of nontrivial extensions, especially when the agent can predict and plan how he will renegotiate with the principal. In the agent's decision makings, he can leverage this prediction into his joint decision of reporting and stopping at each period. Due to the dynamics of the agent's state, the (planned) renegotiation also changes over time. As a result, the principal's mechanism design needs to address the effect of the predicted renegotiation on the characterization of the dynamic incentive compatibility.

\section*{Appendix}
\appendix

For notational compactness, we suppress the notations $h_{t-1}^{\theta}$, $h_{t-1}^{\hat{\theta}}$, and $h^{a}_{t-1}$ from $\alpha, \phi, \xi$, and $\sigma$.

\section{Proof of Proposition \ref{prop:one_shot_deviation_ICS} }\label{app:proof_lemma_one_shot_deviation_ICS}

%
The proof of the \textit{only if} part directly follows from the optimality of truthful reporting and here we only provide the proof of the \textit{if} part.
Suppose, on the contrary, the truthful reporting strategy $\sigma^{*}$ satisfies (\ref{eq:ICS_bellman_1SD}) and (\ref{eq:ICS_T_1SD}) but not (\ref{eq:ICS_bellman}) and (\ref{eq:ICS_T}).
Then there exists a reporting strategy $\sigma'$ and a state $\theta_{t}$, at period $t\in\mathbb{T}$,  such that $V^{\alpha,\phi,\xi,\rho}_{t}(\theta_{t};\sigma')> V^{\alpha,\phi,\xi,\rho}_{t}(\theta_{t};\sigma^{*})$. 
Suppose that the optimal stopping rule with $\sigma^{*}$ calls for stopping and the agent decides to continue by using $\sigma'$, i.e.,
$$
J^{\alpha,\phi,\xi,\rho}_{1,t}(t,\theta_{t})< \mathbb{E}^{\Xi_{\alpha;\sigma'}[h^{\theta}_{t}]}\big[V^{\alpha,\phi,\xi,\rho}_{t+1}(\tilde{\theta}_{t+1}) \big].
$$
Hence, there exists some $\varepsilon>0$ such that 
\begin{equation}\label{eq:app_1_1}
    \mathbb{E}^{\Xi_{\alpha;\sigma'}[h^{\theta}_{t}]}\big[V^{\alpha,\phi,\xi,\rho}_{t+1}(\tilde{\theta}_{t+1}) \big] \geq J^{\alpha,\phi,\xi,\rho}_{1,t}(t,\theta_{t}) + 2\varepsilon.
\end{equation}
Let $\sigma''$ be the reporting strategy such that if $\sigma''$ and $\sigma'$ have the same reporting strategies from period $t$ to $t+k$, for some $k\geq 0$, then
\begin{equation}\label{eq:app_1_2}
    \mathbb{E}^{\Xi_{\alpha;\sigma''}[h^{\theta}_{t}]}\big[V^{\alpha,\phi,\xi,\rho}_{t+1}(\tilde{\theta}_{t+1}) \big] \geq \mathbb{E}^{\Xi_{\alpha;\sigma'}[h^{\theta}_{t}]}\big[V^{\alpha,\phi,\xi,\rho}_{t+1}(\tilde{\theta}_{t+1}) \big] - \varepsilon.
\end{equation}
From (\ref{eq:app_1_1}) and (\ref{eq:app_1_2}), we have
\begin{equation}\label{eq:app_1_3}
    \mathbb{E}^{\Xi_{\alpha;\sigma''}[h^{\theta}_{t}]}\big[V^{\alpha,\phi,\xi,\rho}_{t+1}(\tilde{\theta}_{t+1}) \big] \geq J^{\alpha,\phi,\xi,\rho}_{1,t}(t,\theta_{t}) + \varepsilon.
\end{equation}
Here, (\ref{eq:app_1_3}) implies that any deviation(s) for the periods from $t$ to $t+k$ (reporting truthfully for all other periods) can improve the value $V^{\alpha,\phi,\xi,\rho}_{t}$.

Let $\hat{\sigma}^{s}$ denote the reporting strategy that differs only at period $s$ from $\sigma^{*}$ and $\hat{\sigma}^{s}_{s} = \sigma''_{s}$, for $s\in[t,t+k]$. 
Then, we have
\begin{equation}\label{eq:app_1_4}
    \mathbb{E}^{\Xi_{\alpha;\hat{\sigma}^{t+k-1}}[h^{\theta}_{t}]}\big[V^{\alpha,\phi,\xi,\rho}_{t+1}(\tilde{\theta}_{t+1}) \big] > J^{\alpha,\phi,\xi,\rho}_{1,t}(t,\theta_{t}).
\end{equation}
Now, we look at period $t+k-1$. Because $\sigma^{*}$ satisfies (\ref{eq:ICS_bellman_1SD}) and (\ref{eq:ICS_T_1SD}), we have, for all $\theta_{t+k-1}\in\Theta_{t+k-1}$,
\begin{equation}\label{eq:app_1_5}
    \begin{aligned}
        &\mathbb{E}^{\Xi_{\alpha;\hat{\sigma}^{t+k-2}}[h^{\theta}_{t+k-1}]}\big[V^{\alpha,\phi,\xi,\rho}_{t+k-1}(\theta_{t+k-1}) \big] = V^{\alpha,\phi,\xi,\rho}_{t+k-1}(\theta_{t+k-1})\\
        \geq& \max\Big(J^{\alpha,\phi,\xi,\rho}_{1,t+k-1}(t+k-1,\theta_{t+k-1}, \hat{\sigma}^{t+k-1}_{t+k-1}(\theta_{t+k-1})|h^{\theta}_{t+k-2})
        , \mathbb{E}^{\Xi_{\alpha;\hat{\sigma}^{t+k-1}}[h^{\theta}_{t+k-1}]}\big[V^{\alpha,\phi,\xi,\rho}_{t+k}(\tilde{\theta}_{t+k}) \big]  \Big)\\
        =& V^{\alpha,\phi,\xi,\rho}_{t+k-1}(\theta_{t+k-1};\hat{\sigma}^{t+k-1}_{t+k-1}).
    \end{aligned}
\end{equation}
Then, 
\begin{equation}\label{eq:app_1_6}
    \begin{aligned}
        \mathbb{E}^{\Xi_{\alpha;\hat{\sigma}^{t+k-2}}[h^{\theta}_{t}]}\big[V^{\alpha,\phi,\xi,\rho}_{t+1}(\tilde{\theta}_{t+1}) \big] \geq  \mathbb{E}^{\Xi_{\alpha;\hat{\sigma}^{t+k-1}}[h^{\theta}_{t}]}\big[V^{\alpha,\phi,\xi,\rho}_{t+1}(\tilde{\theta}_{t+1}) \big].
    \end{aligned}
\end{equation}
From (\ref{eq:app_1_4}) and (\ref{eq:app_1_6}), we have
$$
\mathbb{E}^{\Xi_{\alpha;\hat{\sigma}^{t+k-2}}[h^{\theta}_{t}]}\big[V^{\alpha,\phi,\xi,\rho}_{t+1}(\tilde{\theta}_{t+1}) \big] > J^{\alpha,\phi,\xi,\rho}_{1,t}(t,\theta_{t}).
$$
Backward induction yields
$$
\mathbb{E}^{\Xi_{\alpha;\hat{\sigma}^{t}}[h^{\theta}_{t}]}\big[V^{\alpha,\phi,\xi,\rho}_{t+1}(\tilde{\theta}_{t+1}) \big] > J^{\alpha,\phi,\xi,\rho}_{1,t}(t,\theta_{t}),
$$
which contradicts the fact that $\sigma^{*}$ satisfies (\ref{eq:ICS_bellman_1SD}) and (\ref{eq:ICS_T_1SD}). 

Following the similar analysis, we can prove the cases when the optimal stopping rule with truthful $\sigma^{*}$ (1) calls for stopping and the agent decides to stop, (2) calls for continuing and the agent decides to continue, and (3) calls for continuing and the agent decides to stop.

\qed

\section{Proof of Lemma \ref{lemma:marginal_value_represt} }\label{app:proof_lemma:marginal_value_represt}

%
We prove (\ref{eq:ex_ante_payoff_marginal}) here. The proof of (\ref{eq:interim_payoff_marginal}) can be done analogously.
For any $\tau\in\mathbb{T}$, the agent's ex-ante expected payoff (\ref{eq:ex_ante_payoff_agent}) can be written as
$$
\begin{aligned}
    J^{\alpha,\phi,\xi,\rho}_{1}(\tau)
    = \mathbb{E}^{\Xi_{\alpha}}\Bigg[ \sum^{\tau-1}_{t=1}\Big[ \delta^{t+1}\big[& u_{1,t+1}(\tilde{\theta}_{t+1}, \alpha_{t+1}( \tilde{\theta}_{t+1} )) +\xi_{t+1}(\tilde{\theta}_{t+1})\big]+\rho(t+1) \\
    &+ \delta^{t}\big[ \phi_{t}(\tilde{\theta}_{t}) -\xi_{t}(\tilde{\theta}_{t}) \big] - \rho(t)\Big]   + \delta\big[u_{1,1}(\tilde{\theta}_{1},\alpha_{1}(\tilde{\theta}_{1}))+\xi_{1}(\tilde{\theta}_{1})] +\rho(1) \Bigg].
\end{aligned}
$$
From law of total expectation, we have
$$
\begin{aligned}
    J^{\alpha,\phi,\xi,\rho}_{1}(\tau)
    = \mathbb{E}^{\Xi_{\alpha}}\Bigg[ \sum^{\tau-1}_{t=1}\Big[ \mathbb{E}^{\alpha;\theta_{t}}&\Big[ \delta^{t+1}\big[ u_{1,t+1}(\tilde{\theta}_{t+1}, \alpha_{t+1}( \tilde{\theta}_{t+1} )) +\xi_{t+1}(\tilde{\theta}_{t+1})\big]+\rho(t+1) \\
    &+ \delta^{t}\big[ \phi_{t}(\tilde{\theta}_{t}) -\xi_{t}(\tilde{\theta}_{t}) \big]\Big] - \rho(t) \Big]   + \delta\big[u_{1,1}(\tilde{\theta}_{1},\alpha_{1}(\tilde{\theta}_{1}))+\xi_{1}(\tilde{\theta}_{1})] +\rho(1) \Bigg]\\
    =& \mathbb{E}^{\Xi_{\alpha}}\Big[ \sum_{s=1}^{\tau-1} L^{\alpha, \phi, \xi, \rho }_{s}(\tilde{\theta}_{s}) -\rho(s)\Big] + J^{\alpha, \phi, \xi, \rho}_{1}(1).
\end{aligned}
$$

\qed

\section{Proof of Lemma \ref{prop:threshold_rule}}\label{app:prop_threshold_rule}

Let $\hat{\sigma}[t]$ be the one-shot deviation strategy that reports $\hat{\theta}_{t}$ for the true state $\theta_{t}$ at $t$.
Let $\Omega^{*}[\hat{\sigma}[t]]$ be the optimal stopping time rule defined in (\ref{eq:stopping_rule_formulation}) with the stopping region $\Lambda^{\alpha, \phi, \xi, \rho }_{1,t}(t;\hat{\sigma}[t])$ given in (\ref{eq:original_stopping_zone}) (equivalently, (\ref{eq:rewrite_stopping_region})).
Suppose that at period $t$ the agent observes a state $\theta_{t}\in \Lambda^{\alpha, \phi, \xi, \rho }_{1,t}(t;\hat{\sigma}[t])$.
Hence, the agent stops at $t$ optimally.
Then, we obtain, for every $\theta'_{t}\leq \theta_{t}$,
$$
\rho(t) \geq \bar{\mu}^{\alpha, \phi, \xi, \rho}_{t}(\theta_{t}, \hat{\theta}_{t})  \geq \bar{\mu}^{\alpha, \phi, \xi, \rho}_{t}(\theta'_{t},\hat{\theta}_{t}),
$$
where the inequality is due to Lemma \ref{lemma:nonotone_continuing_value}.
Therefore, $\theta'_{t}\in \Lambda^{\alpha, \phi, \xi, \rho }_{1,t}(t;\hat{\sigma}[t])$ for every $\theta'_{t}\leq \theta_{t}$, which implies that $\Lambda^{\alpha, \phi, \xi, \rho }_{1,t}(t;\hat{\sigma}[t])$ is an interval left-bounded by $\underline{\theta}_{t}$.
Since $L^{\alpha, \phi, \xi, \rho }_{t}$ is continuous, $\Lambda^{\alpha, \phi, \xi, \rho }_{1,t}(t;\hat{\sigma}[t])$ is closed.
Hence, according to Assumption \ref{assp:full_support}, there exists some $\eta(t)\in \Theta_{t}$ such that $\Lambda^{\alpha, \phi, \xi, \rho }_{1,t}(t;\hat{\sigma}[t]) = [\underline{\theta}_{t}, \eta(t)]$. 
\qed

\section{Proof of Lemma \ref{lemma:threshold_rule}}\label{app:lemma_threshold_rule}

%
Let $\Omega[\hat{\sigma}[t]]|\eta$ and $\Omega[\hat{\sigma}[t]]|\eta'$ denote the optimal stopping rule with threshold functions $\eta$ and $\eta'$, respectively.
Let $\tau_{\eta}$ and $\tau_{\eta'}$ denote the expected realized stopping time from $\Omega[\hat{\sigma}[t]]|\eta$ and $\Omega[\hat{\sigma}[t]]|\eta'$, respectively.
Without loss of generality, suppose $\eta(t)<\eta'(t)$ for some $t\in \mathbb{T}$. 
Here, we obtain the probability of $\tau_{\eta}=t$ as:
$$
\begin{aligned}
    P_{r}(\tau_{\eta}=t) = P_{r}(\theta_{t}\leq\eta(t), \tau_{\eta}>t-1)
    =& \mathbb{E}^{\Xi_{\alpha}}\Bigg[\mathbb{E}^{\alpha|\theta_{t-1}}\Big[\mathbf{1}_{\{\tilde{\theta}_{t}\leq \eta(t)\}}   \Big] \mathbf{1}_{\{\tau_{\eta}>t-1\}}  \Bigg].
\end{aligned}
$$
We can obtain $P_{r}(\tau_{\eta'}=t)$ in a similar way. Then, 
\begin{equation}\label{eq:Appendix_4_1}
    \begin{aligned}
    P_{r}(\tau_{\eta'}=t)-P_{r}(\tau_{\eta}=t)
    =& \mathbb{E}^{\Xi_{\alpha}}\Bigg[\mathbb{E}^{\alpha|\theta_{t-1}}\Big[\mathbf{1}_{\{\tilde{\theta}_{t}\leq \eta'(t)\}}   \Big] \mathbf{1}_{\{\tau_{\eta'}>t-1\}}  \Bigg]-\mathbb{E}^{\Xi_{\alpha}}\Bigg[\mathbb{E}^{\alpha|\theta_{t-1}}\Big[\mathbf{1}_{\{\tilde{\theta}_{t}\leq \eta(t)\}}   \Big] \mathbf{1}_{\{\tau_{\eta}>t-1\}}  \Bigg]\\
    =& \mathbb{E}^{\Xi_{\alpha}}\Bigg[ \mathbb{E}^{\alpha|\theta_{t-1}}\Big[\mathbf{1}_{\{\eta(t)\leq \tilde{\theta}_{t} \leq \eta'(t)  \} }   \Big] \mathbf{1}_{\tau_{\eta}>t-1} \Bigg].
\end{aligned}
\end{equation}
Since $\tau_{\eta}=\tau_{\eta'}$, the probabilities $P_{r}(\tau_{\eta'}=t)$ and $P_{r}(\tau_{\eta}=t)$ are equal, i.e., (\ref{eq:Appendix_4_1}) equals $0$.
However, from Assumption \ref{assp:full_support}, we know $\mathbb{E}^{\alpha|\theta_{t-1}}\Big[\mathbf{1}_{\{\eta(t)\leq \tilde{\theta}_{t} \leq \eta'(t)  \} }   \Big]>0$ and $P_{r}(\tau_{\eta}>t-1)>0$, which implies that  
$$
\mathbb{E}^{\Xi_{\alpha}}\Bigg[ \mathbb{E}^{\alpha|\theta_{t-1}}\Big[\mathbf{1}_{\{\eta(t)\leq \tilde{\theta}_{t} \leq \eta'(t)  \} }   \Big] \mathbf{1}_{\tau_{\eta}>t-1} \Bigg]>0.
$$
This contradiction implies that $\eta$ is unique.

\qed

\section{Proof of Proposition \ref{proposition:sufficient_ICS}}\label{app:proposition_sufficient_ICS}

%
From the construction of $\xi$ in (\ref{eq:payment_xi_eta}), we have
\begin{equation}\label{eq:app_sufficiency_1}
    \begin{aligned}
   \xi_{t}(\hat{\theta}_{t}) - \xi_{t}(\theta_{t}) =& \delta^{-t}\beta^{\alpha}_{S,t}(\hat{\theta}_{t}) -\delta^{-t}\beta^{\alpha}_{S,t}(\theta_{t}) + u_{1,t}(\theta_{t}, \alpha_{t}(\theta_{t}))  - u_{1,t}(\hat{\theta}_{t}, \alpha_{t}(\hat{\theta}_{t}))\\
   =& \delta^{-t}\beta^{\alpha}_{S,t}(\hat{\theta}_{t}) -\delta^{-t}\beta^{\alpha}_{S,t}(\theta_{t}) -(u_{1,t}(\hat{\theta}_{t}, \alpha_{t}(\hat{\theta}_{t}))- u_{1,t}(\theta_{t}, \alpha_{t}(\hat{\theta}_{t})))\\
   &+ u_{1,t}(\theta_{t}, \alpha_{t}(\theta_{t}))  - u_{1,t}(\theta_{t}, \alpha_{t}(\hat{\theta}_{t}))).
\end{aligned}
\end{equation}
From the definition of $\ell^{\alpha}_{S,t}$ in (\ref{eq:path_length_eta_stop}) and the condition (\ref{eq:potential_function_relation_S}), 
$$
\begin{aligned}
   \text{R.H.S. of (\ref{eq:app_sufficiency_1}) } =&\delta^{-t}\beta^{\alpha}_{S,t}(\hat{\theta}_{t}) -\delta^{-t}\beta^{\alpha}_{S,t}(\theta_{t}) + u_{1,t}(\theta_{t}, \alpha_{t}(\theta_{t})) -\ell^{\alpha}_{S,t}(\hat{\theta}_{t}, \theta_{t})  \\
   &+ u_{1,t}(\theta_{t}, \alpha_{t}(\theta_{t}))  - u_{1,t}(\theta_{t}, \alpha_{t}(\hat{\theta}_{t})))\\
   \leq& u_{1,t}(\theta_{t}, \alpha_{t}(\theta_{t}))  - u_{1,t}(\theta_{t}, \alpha_{t}(\hat{\theta}_{t}))),
\end{aligned}
$$
which implies
\begin{equation}\label{eq:app_sufficient_stop_terminal}
    u_{1,t}(\theta_{t}, \alpha_{t}(\theta_{t}))  +\xi_{t}(\theta_{t}) \geq u_{1,t}(\theta_{t}, \alpha_{t}(\hat{\theta}_{t}))) + \xi_{t}(\hat{\theta}_{t}),
\end{equation}
i.e., 
$$
U^{\alpha,\phi,\xi,\rho}_{S, t}(\theta_{t}|h^{\theta}_{t-1})\geq U^{\alpha,\phi,\xi,\rho}_{S, t}(\theta_{t},\hat{\theta}_{t}|h^{\theta}_{t-1}).
$$

From the construction of $\phi$ in (\ref{eq:payment_phi_eta_1}), we have, for any $\tau\in\mathbb{T}_{t+1}$,
\begin{equation}\label{eq:suff_nonstop_app}
    \begin{split}
        \phi_{t}(\hat{\theta}_{t}) - \phi_{t}(\theta_{t}) 
       = & \beta^{\alpha}_{\bar{S},t}(\hat{\theta}_{t}) - \mathbb{E}^{\alpha|\hat{\theta}_{t}}\Big[\beta_{\bar{S},t+1}(\tilde{\theta}_{t+1})\Big] - u_{1,t}(\hat{\theta}_{t}, \alpha_{t}(\hat{\theta}_{t}))\\
       -&\beta^{\alpha}_{\bar{S},t}(\theta_{t}) + \mathbb{E}^{\alpha|\theta_{t}}\Big[\beta_{\bar{S},t+1}(\tilde{\theta}_{t+1})\Big] + u_{1,t}(\theta_{t}, \alpha_{t}(\theta_{t}))\\
       =& \beta^{\alpha}_{\bar{S},t}(\hat{\theta}_{t}) - \beta^{\alpha}_{\bar{S},t}(\theta_{t}) + \mathbb{E}^{\alpha|\theta_{t}}\Big[ \sum_{s=t}^{T}\delta^{s} u_{1,s}(\tilde{\theta}_{s},\alpha_{s}(\tilde{\theta}_{s}))+ \sum_{s=t+1}^{T-1}\delta^{s}\phi_{s}(\tilde{\theta}_{s}) +\delta^{\tau}\xi_{T}(\tilde{\theta}_{T})\Big]\\
       -&\mathbb{E}^{\alpha|\hat{\theta}_{t}}\Big[ \sum_{s=t}^{\tau-1}\delta^{s} u_{1,s}(\tilde{\theta}_{s},\alpha_{s}(\tilde{\theta}_{s}))+ \sum_{s=t+1}^{\tau-1}\delta^{s}\phi_{s}(\tilde{\theta}_{s}) +\beta^{\alpha}_{\bar{S},\tau}(\tilde{\theta}_{\tau}) \Big].
    \end{split}
\end{equation}
From the condition (\ref{eq:potential_function_relation_nS}), we have

\begin{equation}\label{eq:suff_nonstop_app2}
    \begin{split}
        \text{R.H.S. of (\ref{eq:suff_nonstop_app})} \leq & \inf_{\tau\in\mathbb{T}_{t}}\Big\{ \mathbb{E}^{\alpha|\hat{\theta}_{t}} \Big[ \sum_{s=t}^{\tau }\delta^{s}\big[u_{1,s}(\tilde{\theta}_{s}, \alpha_{s}(\tilde{\theta}_{s}))\big] +\sum_{s=t+1}^{\tau-1}\delta^{s}\phi_{s}(\tilde{\theta}_{s})  +\delta^{\tau }\xi_{\tau }(\tilde{\theta}_{\tau })   \Big]\\
        &-\mathbb{E}^{\alpha|\theta_{t},\hat{\theta}_{t}}\Big[ \sum_{s=t}^{\tau}\delta^{s}\big[u_{1,s}(\tilde{\theta}_{s}, \alpha_{s}(\tilde{\theta}_{s}))\big] +\sum_{s=t+1}^{\tau-1}\delta^{s}\phi_{s}(\tilde{\theta}_{s})  +\delta^{\tau}\xi_{\tau}(\tilde{\theta}_{\tau })   \Big] \Big\}-\sup_{\tau\in\mathbb{T}_{t}}\rho(\tau)\\
        &+ \mathbb{E}^{\alpha|\theta_{t}}\Big[ \sum_{s=t}^{T}\delta^{s} u_{1,s}(\tilde{\theta}_{s},\alpha_{s}(\tilde{\theta}_{s}))+ \sum_{s=t+1}^{T-1}\delta^{s}\phi_{s}(\tilde{\theta}_{s}) +\delta^{T}\xi_{T}(\tilde{\theta}_{T})\Big]\\
       &-\mathbb{E}^{\alpha|\hat{\theta}_{t}}\Big[ \sum_{s=t}^{\tau-1}\delta^{s}u_{1,s}(\tilde{\theta}_{s},\alpha_{s}(\tilde{\theta}_{s}))+ \sum_{s=t+1}^{\tau-1}\delta^{s}\phi_{s}(\tilde{\theta}_{s}) +\beta^{\alpha}_{\bar{S},\tau}(\tilde{\theta}_{\tau}) \Big].
    \end{split}
\end{equation}
From the condition (\ref{eq:potential_function_relation_ine}), $\beta^{\alpha}_{\bar{S},t}(\theta_{t})\geq \beta^{\alpha}_{S,t}(\theta_{t})$, for all $\theta_{t}\in\Theta_{t}$, $t\in\mathbb{T}$. Hence, $\beta^{\alpha}_{\bar{S},t}(\theta_{t})\geq \xi_{t}(\theta_{t})$ $+ u_{1,t}(\theta_{t},\alpha_{t}(\theta_{t}))$, for all $\theta_{t}\in \mathbb{T}$, $t\in\mathbb{T}$.
Then, 
\begin{equation}\label{eq:app_haha1_1}
    \begin{split}
        \inf_{\tau\in\mathbb{T}_{t}}\Big\{ \mathbb{E}^{\alpha|\hat{\theta}_{t}} &\Big[ \sum_{s=t}^{\tau }\delta^{s}\big[u_{1,s}(\tilde{\theta}_{s}, \alpha_{s}(\tilde{\theta}_{s}))\big] +\sum_{s=t+1}^{\tau-1}\delta^{s}\phi_{s}(\tilde{\theta}_{s})  +\delta^{\tau }\xi_{\tau }(\tilde{\theta}_{\tau })   \Big]\\
        &-\mathbb{E}^{\alpha|\hat{\theta}_{t}}\Big[ \sum_{s=t}^{\tau-1}\delta^{s}u_{1,s}(\tilde{\theta}_{s},\alpha_{s}(\tilde{\theta}_{s}))+ \sum_{s=t+1}^{\tau-1}\delta^{s}\phi_{s}(\tilde{\theta}_{s}) +\beta^{\alpha}_{\bar{S},\tau}(\tilde{\theta}_{\tau})  \Big]\Big\}\\
        \leq\inf_{\tau\in\mathbb{T}_{t}}\Big\{ \mathbb{E}^{\alpha|\hat{\theta}_{t}} &\Big[ \sum_{s=t}^{\tau }\delta^{s}\big[u_{1,s}(\tilde{\theta}_{s}, \alpha_{s}(\tilde{\theta}_{s}))\big] +\sum_{s=t+1}^{\tau-1}\delta^{s}\phi_{s}(\tilde{\theta}_{s})  +\delta^{\tau }\xi_{\tau }(\tilde{\theta}_{\tau })   \Big]\\
        &-\mathbb{E}^{\alpha|\hat{\theta}_{t}}\Big[ \sum_{s=t}^{\tau}\delta^{s}u_{1,s}(\tilde{\theta}_{s},\alpha_{s}(\tilde{\theta}_{s}))+ \sum_{s=t+1}^{\tau-1}\delta^{s}\phi_{s}(\tilde{\theta}_{s}) +\delta^{\tau}\xi_{\tau}(\tilde{\theta}_{\tau}) \Big]\Big\}\\
        =& 0.
    \end{split}
\end{equation}
Hence, from (\ref{eq:app_haha1_1}), we have
\begin{equation}\label{eq:app_94_1}
    \begin{aligned}
   \text{R.H.S. of (\ref{eq:suff_nonstop_app2}) }
       \leq&\mathbb{E}^{\alpha|\theta_{t}}\Big[\sum_{s=t}^{T}\delta^{s} u_{1,s}(\tilde{\theta}_{s},\alpha_{s}(\tilde{\theta}_{s}))+ \sum_{s=t+1}^{T-1}\delta^{s}\phi_{s}(\tilde{\theta}_{s}) +\delta^{T}\xi_{T}(\tilde{\theta}_{T})\Big]\\ &+\inf_{\tau\in\mathbb{T}_{t}}\Big\{-\mathbb{E}^{\alpha|\theta_{t},\hat{\theta}_{t}}\Big[ \sum_{s=t}^{\tau}\delta^{s}\big[u_{1,s}(\tilde{\theta}_{s}, \alpha_{s}(\tilde{\theta}_{s}))\big] +\sum_{s=t+1}^{\tau-1}\delta^{s}\phi_{s}(\tilde{\theta}_{s})  +\delta^{\tau}\xi_{\tau}(\tilde{\theta}_{\tau })    \Big] \Big\}-\sup_{\tau\in\mathbb{T}_{t}}\rho(\tau).
\end{aligned}
\end{equation}
From the construction of $\rho$ in (\ref{eq:construct_rho_primal}) and Lemma \ref{lemma:monotonicity_of_L}, we have, for some $\tau'\in\mathbb{T}_{t}$
\begin{equation}\label{eq:app_redut}
    \begin{aligned}
       \text{R.H.S. of (\ref{eq:app_94_1})} \leq&\mathbb{E}^{\alpha|\theta_{t}}\Big[\sum_{s=t}^{\tau'}\delta^{s} u_{1,s}(\tilde{\theta}_{s},\alpha_{s}(\tilde{\theta}_{s}))+ \sum_{s=t+1}^{\tau'-1}\delta^{s}\phi_{s}(\tilde{\theta}_{s}) +\delta^{\tau'}\xi_{\tau'}(\tilde{\theta}_{\tau'}) +\rho(\tau')\Big]\\ &+\inf_{\tau\in\mathbb{T}_{t}}\Big\{-\mathbb{E}^{\alpha|\theta_{t},\hat{\theta}_{t}}\Big[ \sum_{s=t}^{\tau}\delta^{s}\big[u_{1,s}(\tilde{\theta}_{s}, \alpha_{s}(\tilde{\theta}_{s}))\big] +\sum_{s=t+1}^{\tau-1}\delta^{s}\phi_{s}(\tilde{\theta}_{s})  +\delta^{\tau}\xi_{\tau}(\tilde{\theta}_{\tau }) \Big] \Big\}\\
       &+\inf_{\tau\in\mathbb{T}_{t}}\Big\{-\rho(\tau)\Big\},
    \end{aligned}
\end{equation}
which can be further bounded as
\begin{equation*}
    \begin{aligned}
       \text{R.H.S. of (\ref{eq:app_redut})} \leq& \sup_{\tau\in\mathbb{T}_{t}}\Big\{\mathbb{E}^{\alpha|\theta_{t}}\Big[\sum_{s=t}^{\tau}\delta^{s} u_{1,s}(\tilde{\theta}_{s},\alpha_{s}(\tilde{\theta}_{s}))+ \sum_{s=t+1}^{\tau-1}\delta^{s}\phi_{s}(\tilde{\theta}_{s}) +\delta^{\tau}\xi_{\tau}(\tilde{\theta}_{\tau}) +\rho(\tau)\Big]\Big\}\\ +&\inf_{\tau\in\mathbb{T}_{t}}\Big\{-\mathbb{E}^{\alpha|\theta_{t},\hat{\theta}_{t}}\Big[ \sum_{s=t}^{\tau}\delta^{s}\big[u_{1,s}(\tilde{\theta}_{s}, \alpha_{s}(\tilde{\theta}_{s}))\big] +\sum_{s=t+1}^{\tau-1}\delta^{s}\phi_{s}(\tilde{\theta}_{s})  +\delta^{\tau}\xi_{\tau}(\tilde{\theta}_{\tau }) +\rho(\tau)\Big] \Big\}\\
       =& \sup_{\tau\in\mathbb{T}_{t}}\Big\{\mathbb{E}^{\alpha|\theta_{t}}\Big[\sum_{s=t}^{\tau}\delta^{s} u_{1,s}(\tilde{\theta}_{s},\alpha_{s}(\tilde{\theta}_{s}))+ \sum_{s=t+1}^{\tau-1}\delta^{s}\phi_{s}(\tilde{\theta}_{s}) +\delta^{\tau}\xi_{\tau}(\tilde{\theta}_{\tau}) +\rho(\tau)\Big]\Big\}\\ -&\sup_{\tau\in\mathbb{T}_{t}}\Big\{\mathbb{E}^{\alpha|\theta_{t},\hat{\theta}_{t}}\Big[ \sum_{s=t}^{\tau}\delta^{s}\big[u_{1,s}(\tilde{\theta}_{s}, \alpha_{s}(\tilde{\theta}_{s}))\big] +\sum_{s=t+1}^{\tau-1}\delta^{s}\phi_{s}(\tilde{\theta}_{s})  +\delta^{\tau}\xi_{\tau}(\tilde{\theta}_{\tau }) +\rho(\tau)\Big] \Big\}.
    \end{aligned}
\end{equation*}

Hence, from the definition of $U^{\alpha,\phi,\xi,\rho}_{\bar{S},t}$, we have
$$
U^{\alpha,\phi,\xi,\rho}_{\bar{S},t}(\theta_{t}|h^{\theta}_{t-1})\geq U^{\alpha,\phi,\xi,\rho}_{\bar{S},t}(\theta_{t}, \hat{\theta}_{t}|h^{\theta}_{t-1}).
$$

Therefore, we can conclude that the mechanism is DIC.

\qed

\section{Proof of Lemma \ref{lemma:envelope_conditions}}\label{app:proof_lemma:envelope_conditions}

%
Let $\tilde{m}_{t}$ be uniformly distributed over $(0,1)$.
Given the kernel $K_{t}$, define the inverse of $F_{t}(\cdot|\theta_{t-1}, a_{t-1})$ as follows:
$$
\begin{aligned}
   F^{-1}_{t}(m_{t}|\theta_{t-1}, a_{t-1})= \inf\{\theta_{t}\in\Theta_{t}: F_{t}(\theta_{t}|\theta_{t-1}, a_{t-1})\geq m_{t}\}.
\end{aligned}
$$
Let $\theta_{t}\in \Theta_{t}$ and $\theta_{t+1}\in \Theta_{t+1}$ be any two realized states at two adjacent periods, for any $t\in \mathbb{T}\backslash\{T\}$.
Then, we have
$$
\begin{aligned}
   \frac{\partial  \theta_{t+1}}{\partial  r}\Big|_{r= \theta_{t}} = &\frac{\partial  F^{-1}_{t+1}(m_{t+1}|r, a_{t}) }{\partial r  }\Big|_{r=\theta_{t}}= \frac{-\partial F_{t+1}(\theta_{t+1}| r, a_{t})}{f_{t+1}(\theta_{t+1}|\theta_{t}, a_{t}) \partial r}\Big|_{r=\theta_{t}}.
\end{aligned}
$$
Then, for any sequence of realized states $\{\theta_{t},\theta_{t+1},\dots, \theta_{t+k}\}$, for some $k>1$, we have
$$
\begin{aligned}
   \frac{\partial  \theta_{t+k}}{\partial  r}\Big|_{r= \theta_{t}} =& \prod_{s=t+1}^{t+k}\frac{\partial  F^{-1}_{s}(m_{s}|r, a_{s-1}) }{\partial r  }\Big|_{r=\theta_{s-1}}= \prod_{s=t+1}^{t+k}\Big[\frac{-\partial F_{s}(\theta_{s}| r, a_{s-1})}{f_{s}(\theta_{s}|\theta_{s-1}, a_{s-1}) \partial r}\Big|_{r=\theta_{s-1}}\Big].
\end{aligned}
$$
In any DIC mechanism, truthful reporting strategy is optimal. Then, the envelope theorem yields the following:
$$
\begin{aligned}
   \frac{\partial U^{\alpha,\phi,\xi,\rho}_{t}(\tau, r|h^{\theta}_{t-1})}{\partial r}\Big|_{r= \theta_{t}} =& \mathbb{E}^{\alpha|\theta_{t}}\Big[\sum_{s=t}^{\tau}\frac{\partial u_{1,s}(r, \alpha_{s}(\tilde{\theta}_{s}))}{\partial r }\Big|_{r=\tilde{\theta}_{s}} \cdot \frac{\partial   \tilde{\theta}_{s}}{ \partial l  }\Big|_{l=\theta_{t}}\Big]\\
   =&\mathbb{E}^{\alpha|\theta_{t}}\Bigg[\sum_{s=t}^{\tau}\frac{\partial u_{1,s}(r, \alpha_{s}(\tilde{\theta}_{s}))}{\partial r }\Big|_{r=\tilde{\theta}_{s}} \cdot \prod_{k=t+1}^{s}\Big[\frac{-\partial F_{k}(\theta_{k}| r, a_{k-1})}{f_{k}(\theta_{k}|\theta_{k-1}, a_{k-1}) \partial r}\Big|_{r=\theta_{k-1}}\Big]\Bigg]\\
   =& \mathbb{E}^{\alpha|\theta_{t}}\Big[\sum_{s=t}^{\tau}\frac{\partial u_{1,s}(r, \alpha_{s}(\tilde{\theta}_{s}))}{\partial r }\Big|_{r=\tilde{\theta}_{s}}  \cdot G_{t,s}(h^{\theta}_{t,s})\Big].
\end{aligned}
$$

\qed

\section{Proof of Proposition \ref{prop:envelope_necessary} }\label{app:prop_formulation_betas}

%
Since $u_{1,t}(\theta_{t}, a_{t})$ is a non-decreasing function of $\theta_{t}$, then
$
\frac{\partial  u_{1,t}(r, a_{t})}{\partial r}\Big|_{r = \theta_{t}} \geq 0
$
, for all $t\in\mathbb{T}$.
From Assumption \ref{assp_monotone_transistions}, we have $\frac{\partial F_{t+1}(\theta_{t+1}|r, a^{t})}{\partial r}\Big|_{r=\theta_{t}} \leq 0$. Therefore, from Lemma \ref{eq:envelope_gamma}, the term $\gamma^{\alpha}_{t}(\tau, \theta_{t}|h^{\theta}_{t-1})$ is non-negative.

From the definition of $\chi^{\alpha, \phi, \xi }_{1,t}(\theta_{t})$ in (\ref{eq:difference_Z}), we have
$$
\begin{aligned}
   &\chi^{\alpha, \phi, \xi }_{1,t}(\theta_{t}) = Z^{\alpha, \phi, \xi}_{1,t}(t+1, \theta_{t}|h^{\theta }_{t-1}) - Z^{\alpha, \phi, \xi}_{1,t}(t, \theta_{t}|h^{\theta }_{t-1})\\
   =& \mathbb{E}^{\alpha|\theta_{t}}\Big[\sum_{s=t}^{t+1}\delta^{s} u_{1,s}(\tilde{\theta}_{s}, \alpha_{s}(\tilde{\theta}_{s})) + \delta^{t+1}\xi_{t+1}(\tilde{\theta}_{t+1}) + \delta^{t}\phi_{t}(\theta_{t})\Big] - [\delta^{t}u_{1,t}(\theta_{t}, \alpha_{t}(\theta_{t}))+ \delta^{t}\xi_{t}(\theta_{t})].
\end{aligned}
$$
Substituting the constructions of $\phi$ and $\xi$ given by (\ref{eq:payment_phi_eta_1}) and (\ref{eq:payment_xi_eta}), respectively, yields
\begin{equation}\label{eq:app_Xi_1}
    \begin{aligned}
   \chi^{\alpha, \phi, \xi }_{1,t}(\theta_{t}) =& \mathbb{E}^{\alpha|\theta_{t}}\Big[\beta_{\bar{S},t}(\theta_{t}) -\beta_{\bar{S},t+1}(\tilde{\theta}_{t+1}) \Big] + \mathbb{E}^{\alpha|\theta_{t}}\Big[\beta_{S,t+1}(\tilde{\theta}_{t+1}) - \beta_{S,t}(\theta_{t}) \Big].
\end{aligned}
\end{equation}

Given the formulations of $\beta_{S,t}$ and $\beta_{\bar{S},t}$ in (\ref{eq:beta_necessary_S}) and (\ref{eq:beta_necessary_SBar}), respectively, we have
\begin{equation}\label{eq:app_Xi_2_1}
    \begin{aligned}
       \text{R.H.S. of (\ref{eq:app_Xi_1}) } =& \mathbb{E}^{\alpha|\theta_{t}}\Big[\sup_{\tau\in\mathbb{T}_{t}}\Big\{\int^{\theta_{t}}_{\theta_{\epsilon,t}} \gamma^{\alpha}_{t}(\tau, r|h^{\theta}_{t-1})dr\Big\} -\sup_{\tau\in\mathbb{T}_{t+1}}\Big\{\int^{\tilde{\theta}_{t+1} }_{\theta_{\epsilon,t+1}} \gamma^{\alpha}_{t+1}(\tau, r|h^{\theta}_{t})dr\Big\} \Big] \\
       &+ \mathbb{E}^{\alpha|\theta_{t}}\Big[\int^{\tilde{\theta}_{t+1}}_{\theta_{\epsilon,t+1}} \gamma^{\alpha}_{t+1}(t+1, r|h^{\theta}_{t})dr  - \int^{\theta_{t}}_{\theta_{\epsilon,t}} \gamma^{\alpha}_{t}(t, r|h^{\theta}_{t-1})dr\Big].
    \end{aligned}
\end{equation}

Since $\gamma^{\alpha}_{t}$ is non-negative for all $t\in \mathbb{T}$, then 
\begin{equation}\label{eq:app_Xi_2}
    \begin{aligned}
       \text{R.H.S. of (\ref{eq:app_Xi_2_1}) } =& \mathbb{E}^{\alpha|\theta_{t}}\Big[\int^{\theta_{t}}_{\theta_{\epsilon,t}} \gamma^{\alpha}_{t}(T, r|h^{\theta}_{t-1})dr -\int^{\tilde{\theta}_{t+1} }_{\theta_{\epsilon,t+1}} \gamma^{\alpha}_{t+1}(T, r|h^{\theta}_{t})dr \Big]\\
       &+ \mathbb{E}^{\alpha|\theta_{t}}\Big[\int^{\tilde{\theta}_{t+1}}_{\theta_{\epsilon,t+1}} \gamma^{\alpha}_{t+1}(t+1, r|h^{\theta}_{t})dr  - \int^{\theta_{t}}_{\theta_{\epsilon,t}} \gamma^{\alpha}_{t}(t, r|h^{\theta}_{t-1})dr\Big].
    \end{aligned}
\end{equation}

Taking partial derivative of $\chi^{\alpha, \phi,\zeta}_{1,t}$ given in (\ref{eq:app_Xi_2}) with respect to $\theta_{t}$ gives

$$
\begin{aligned}
   \frac{\partial \chi^{\alpha, \phi, \xi }_{1,t}(r)  }{\partial r}\Big|_{r=\theta_{t}} =& \mathbb{E}^{\alpha|\theta_{t}}\Big[ \gamma^{\alpha}_{t}(T,\theta_{t}|h^{\theta}_{t-1}) -  \gamma^{\alpha}_{t+1}(T, \tilde{\theta}_{t+1}|h^{\theta}_{t})G_{t,t+1}(\tilde{\theta}_{t+1})  \Big] \\
   &+ \mathbb{E}^{\alpha|\theta_{t}}\Big[\gamma^{\alpha}_{t+1}(t+1, \tilde{\theta}_{t+1}|h^{\theta}_{t})G_{t,t+1}(\tilde{\theta}_{t+1}) - \gamma^{\alpha}_{t}(t, \theta_{t}|h^{\theta}_{t-1})  \Big].
\end{aligned}
$$

From Lemma \ref{lemma:envelope_conditions}, we have
$$
\begin{aligned}
   \mathbb{E}^{\alpha|\theta_{t}}\Big[ \gamma^{\alpha}_{t}(T,\theta_{t}|h^{\theta}_{t-1})dr& - \gamma^{\alpha}_{t}(t, \theta_{t}|h^{\theta}_{t-1}) \Big]\\
   =& \max\Big\{ \mathbb{E}^{\alpha|\theta_{t}}\Big[ \sum_{s=t+1}^{T} \delta^{s}\frac{\partial u_{1,s}(r, \alpha_{s}(\tilde{\theta}_{s})) }{\partial r }\Big|_{r=\tilde{\theta}_{s}} G_{t,s}(h^{\tilde{\theta}}_{t,s})   \Big], \;\; 0 \Big\}\\
   =& \mathbb{E}^{\alpha|\theta_{t}}\Big[ \sum_{s=t+1}^{T} \delta^{s}\frac{\partial u_{1,s}(r, \alpha_{s}(\tilde{\theta}_{s})) }{\partial r }\Big|_{r=\tilde{\theta}_{s}} G_{t,s}(h^{\tilde{\theta}}_{t,s})   \Big],
\end{aligned}
$$
where the second equality is from the fact that $\gamma^{\alpha}_{t}$ is non-negative; and

$$
\begin{aligned}
   \mathbb{E}^{\alpha|\theta_{t}}\Big[  \gamma^{\alpha}_{t+1}(T, \tilde{\theta}_{t+1}|h^{\theta}_{t}) G_{t,s}(h^{\tilde{\theta}}_{t,s}) \Big]=&\mathbb{E}^{\alpha|\theta_{t}}\Big[ \sum_{s=t+1}^{T} \delta^{s}\frac{\partial u_{1,s}(r, \alpha_{s}(\tilde{\theta}_{s})) }{\partial r }\Big|_{r=\tilde{\theta}_{s}} G_{t+1,s}(h^{\tilde{\theta}}_{t+1,s})G_{t,s}(h^{\tilde{\theta}}_{t,s})  \Big]\\
   =& \mathbb{E}^{\alpha|\theta_{t}}\Big[ \sum_{s=t+1}^{T} \delta^{s}\frac{\partial u_{1,s}(r, \alpha_{s}(\tilde{\theta}_{s})) }{\partial r }\Big|_{r=\tilde{\theta}_{s}} G_{t,s}(h^{\tilde{\theta}}_{t,s})   \Big].
\end{aligned}
$$
Also, Assumption \ref{assp_monotone_transistions} implies that $G_{t,t+1}(\theta_{t+1})\geq 0$ for all $\theta_{t+1}\in \Theta_{t+1}$.
Hence, we have
$$
\begin{aligned}
   \frac{\partial \chi^{\alpha, \phi, \xi }_{1,t}(r)  }{\partial r}\Big|_{r=\theta_{t}} =& \mathbb{E}^{\alpha|\theta_{t}}\Big[\gamma^{\alpha}_{t+1}(t+1, \tilde{\theta}_{t+1}|h^{\theta}_{t})G_{t,t+1}(\tilde{\theta}_{t+1}) \Big]
   \geq 0.
\end{aligned}
$$
Therefore, the constructions of potential functions given in (\ref{eq:beta_necessary_S}) and (\ref{eq:beta_necessary_SBar}) satisfy the monotonicity condition specified by Assumption \ref{assp:single_crossing}, i.e., the statement \textit{(iii)} in Proposition \ref{proposition:sufficient_ICS} is satisfied.

\qed

\section{Proof of Proposition \ref{prop:necessary_alpha} }\label{app:thm_sufficient_necessary_allocation}

%
Fix an arbitrary $\hat{\theta}_{\epsilon,t}\in \Theta_{\epsilon}$.
We discuss the following two cases:

\subsection*{1. \texorpdfstring{$\boldsymbol{\theta_{t}\in \Lambda_{t}(t)}:$}{Lg}}
%
%


Let $\hat{\theta}_{t} \in \Lambda_{t}(t)$.
Without loss of generality, suppose $\hat{\theta}_{t}\leq \theta_{t}$. Let $\theta$, $\theta^{1}$, $\theta^{2}\in \Bar{\Theta}_{t} \equiv [\Hat{\theta}_{t}, \theta_{t}]$. 
Since the mechanism is DIC, there exists $\xi$ such that
\begin{equation}\label{eq:app_E_implement}
    \delta^{t}\big[ u_{1,t}(\theta_{t }, \alpha_{t}(\theta_{t})) + \xi_{t}(\theta_{t})\big] +\rho(t) \geq \delta^{t}\big[ u_{1,t}(\theta_{t }, \alpha_{t}(\Hat{\theta}_{t})) + \xi_{t}(\Hat{\theta}_{t})\big]+\rho(t).
\end{equation}
Define 
\begin{equation}\label{eq:app_E_max}
    B_{t}(\theta) \equiv \max_{x\in \Bar{\Theta}_{t} }\delta^{t}\Big[  u_{1,t}(\theta, \alpha_{t}(x)) +\xi_{t}(x) \Big].
\end{equation}
DIC implies that 
$$
\theta \in \argmax_{x\in\Bar{\Theta}_{t}}\delta^{t}\Big[  u_{1,t}(\theta, \alpha_{t}(x)) +\xi_{t}(x) \Big].
$$
Then, we obtain 
$$
\begin{aligned}
    |B_{t}(\theta^{2})- B_{t}(\theta^{1}) | \leq & \max_{x\in\Bar{\Theta}_{t}} \delta^{t}\big|  u_{1,t}( \theta^{2}, \alpha_{t}(x)) - u_{1,t}( \theta^{1}, \alpha_{t}(x)) \big|\\
    =& \max_{x\in\Bar{\Theta}_{t}}\delta^{t}\Big|\int_{\theta^{1}}^{\theta^{2}} \frac{\partial u_{1,t}(y, \alpha_{t}(x))  }{ \partial y } \big|_{y = \theta}  d \theta\Big|\\
    =& \max_{x\in\Bar{\Theta}_{t}}\delta^{t}\Big|\beta^{\alpha}_{S,t}(\theta^{2}) - \beta^{\alpha}_{S,t}(\theta^{1})  \Big|.
\end{aligned}
$$
By Assumption \ref{assp:stopping_condition}, we have that $B_{t}$ is Lipschitz continuous. Thus, $B_{t}$ is differentiable almost everywhere. Therefore, we have
$$
B_{t}(\theta_{t}) - B_{t}(\Hat{\theta}_{t}) = \int^{\theta_{t}}_{\hat{\theta }_{t}} \frac{d B_{t}(y)}{d y} \big|_{y = \theta} d\theta. 
$$
Applying envelope theorem to $B_{t}$ yields
$$
\begin{aligned}
    \frac{d B_{t}(y)}{d y} \big|_{y = \theta} = & \frac{\partial }{\partial x }\big[ \delta^{t} u_{1,t}(x, \alpha_{t}(\theta)) + \xi_{t}(\theta)   \big] \Big|_{x= \theta}\\
    =& \frac{\partial }{\partial x } \delta^{t} u_{1,t}(x, \alpha_{t}(\theta) )\Big|_{x= \theta}\\
    =& \gamma^{\alpha }_{t}  (t, \theta |h^{\theta}_{t-1}).
\end{aligned}
$$

Therefore, we have 
$$
\begin{aligned}
   \beta^{\alpha}_{S,t}(\theta_{t}) - \beta^{\alpha}_{S,t}(\hat{\theta}_{t}) =& B_{t}(\theta_{t}) - B_{t}(\Hat{\theta}_{t}) \\
   =& \delta^{t}\big[ u_{1,t}(\theta_{t}, \alpha_{t}(\theta_{t})) + \xi_{t}(\theta_{t})\big] - \delta^{t}\big[u_{1,t}(\hat{\theta}_{t}, \alpha_{t}(\hat{\theta}_{t})) + \xi_{t}(\hat{\theta}_{t})\big]
\end{aligned}
$$
From the definition of $\ell^{\alpha}_{S,t}(\theta_{t},\hat{\theta}_{t})$, we have
$$
\begin{aligned}
    \ell^{\alpha}_{S,t}(\theta_{t}, \Hat{\theta}_{t}) 
    =& \delta^{t}u_{1,t}(\theta_{t}, \alpha_{t}(\theta_{t})) - \delta^{t} u_{1,t}(\Hat{\theta}_{t}, \alpha_{t}(\theta_{t}))\\
    =& \delta^{t}u_{1,t}(\theta_{t}, \alpha_{t}(\theta_{t})) - \delta^{t} u_{1,t}(\hat{\theta}_{t}, \alpha_{t}(\hat{\theta}_{t})  ) +\delta^{t} u_{1,t}(\hat{\theta}_{t}, \alpha_{t}(\hat{\theta}_{t})  ) - \delta^{t} u_{1,t}(\Hat{\theta}_{t}, \alpha_{t}(\theta_{t})) \\
    \geq & \delta^{t}\big[u_{1,t}(\theta_{t}, \alpha_{t}(\theta_{t})) -  u_{1,t}(\hat{\theta}_{t}, \alpha_{t}(\hat{\theta}_{t})  ) + \xi_{t}(\theta_{t}) - \xi_{t}(\Hat{\theta}_{t})\big]\\
    =& \beta^{\alpha}_{S,t}(\theta_{t}) - \beta^{\alpha}_{S,t}(\hat{\theta}_{t}).
\end{aligned}
$$
\subsection*{2. $\theta_{t}\not\in \Lambda_{t}(t):$}



Similar to the case when $\theta_{t}\in \Lambda_{t}(t)$, DIC implies the existence of $\phi$ and $\xi_{s}$ such that
$$
\begin{aligned}
    \delta^{t}\big[ u_{1,t}(\theta_{t}, \alpha_{t}(\theta_{t}))+ &\xi_{t}(\theta_{t}) \big]+\bar{\mu}^{\alpha,\phi, \xi,\rho}_{t}(\theta_{t}) \\
    \geq& \delta^{t}\big[ u_{1,t}(\theta_{t}, \alpha_{t}(\hat{\theta}_{t}))+ \xi_{t}(\hat{\theta}_{t}) \big]+\bar{\mu}^{\alpha,\phi, \xi,\rho}_{t}(\theta_{t}, \hat{\theta}_{t}).
\end{aligned}
$$
Define 
\begin{equation}\label{eq:app_E_max_B}
    B'_{t}(\theta) \equiv \max_{x\in \Bar{\Theta}_{t} }\Big[ \delta^{t}\big[ u_{1,t}(\theta, \alpha_{t}(x)) + \xi_{t}(x)\big]+\bar{\mu}^{\alpha,\phi, \xi,\rho}_{t}(\theta, x) \Big].
\end{equation}
Since the mechanism is DIC, we have
$$
\theta \in \argmax_{x\in\Bar{\Theta}_{t}} \Big[ \delta^{t}\big[ u_{1,t}(\theta, \alpha_{t}(x)) + \xi_{t}(x)\big]+\bar{\mu}^{\alpha,\phi, \xi,\rho}_{t}(\theta, x) \Big].
$$
%
Envelope theorem yields the following:
\begin{equation}\label{eq:app_theorem_1}
    \begin{aligned}
      \beta^{\alpha}_{\bar{S},t}(\Hat{\theta}_{t})-\beta^{\alpha}_{\bar{S},t}(\theta_{t}) =& B'_{t}(\Hat{\theta}_{t})-B'_{t}(\theta_{t})\\
      =& \delta^{t}\big[ u_{1,t}(\hat{\theta}_{t}, \alpha_{t}(\hat{\theta}_{t}))+\xi_{t}(\hat{\theta}_{t}) \big]+\bar{\mu}^{\alpha,\phi, \xi,\rho}_{t}(\hat{\theta}_{t})\\
     &- \delta^{t}\big[ u_{1,t}(\theta_{t}, \alpha_{t}(\theta_{t}))+ \xi_{t}(\theta_{t}) \big]-\bar{\mu}^{\alpha,\phi, \xi,\rho}_{t}(\theta_{t}).
\end{aligned} 
\end{equation}
%
%
From the definition of $\bar{\mu}^{\alpha,\phi, \xi,\rho}_{t}$, (\ref{eq:app_theorem_1}) can be extended as follows:
\begin{equation}\label{eq:app_theorem_2}
    \begin{aligned}
       &\beta^{\alpha}_{\bar{S},t}(\Hat{\theta}_{t})-\beta^{\alpha}_{\bar{S},t}(\theta_{t}) \\
       =&   \sup_{\tau\in\mathbb{T}_{t+1}}\Bigg\{\mathbb{E}^{\alpha|\hat{\theta}_{t}}\Big[\sum_{s=t}^{\tau}\delta^{s}\big[u_{1,s}(\tilde{\theta}_{s}, \alpha_{s}(\tilde{\theta}_{s})) +\sum_{s=t+1}^{\tau-1} \phi_{s}(\tilde{\theta}_{s}) +\xi_{\tau}(\tilde{\theta}_{\tau}) \big] +\rho(\tau)\Big] \Bigg\}\\
       &- \sup_{\tau\in\mathbb{T}_{t+1}}\Bigg\{\mathbb{E}^{\alpha|\theta_{t}}\Big[\sum_{s=t}^{\tau}\delta^{s}\big[u_{1,s}(\tilde{\theta}_{s}, \alpha_{s}(\tilde{\theta}_{s}))+\sum_{s=t+1}^{\tau-1} \phi_{s}(\tilde{\theta}_{s}) +\xi_{\tau}(\tilde{\theta}_{\tau}) \big] +\rho(\tau)\Big] \Bigg\} + \phi_{t}(\hat{\theta}_{t}) -\phi_{t}(\theta_{t}).
    \end{aligned}
\end{equation}
Since the mechanism is DIC, we have

\begin{equation*}
    \begin{aligned}
       \text{R.H.S. of (\ref{eq:app_theorem_2})  } \leq & \sup_{\tau\in\mathbb{T}_{t+1}}\Bigg\{\mathbb{E}^{\alpha|\hat{\theta}_{t}}\Big[\sum_{s=t}^{\tau}\delta^{s}\big[u_{1,s}(\tilde{\theta}_{s}, \alpha_{s}(\tilde{\theta}_{s})) +\sum_{s=t+1}^{\tau-1} \phi_{s}(\tilde{\theta}_{s}) +\xi_{\tau}(\tilde{\theta}_{\tau}) \big] +\rho(\tau)\Big] \Bigg\}\\
       &- \sup_{\tau\in\mathbb{T}_{t+1}}\Bigg\{\mathbb{E}^{\alpha|\theta_{t}}\Big[\sum_{s=t}^{\tau}\delta^{s}\big[u_{1,s}(\tilde{\theta}_{s}, \alpha_{s}(\tilde{\theta}_{s}))+\sum_{s=t+1}^{\tau-1} \phi_{s}(\tilde{\theta}_{s}) +\xi_{\tau}(\tilde{\theta}_{\tau}) \big] +\rho(\tau)\Big] \Bigg\} \\
       &+ \sup_{\tau\in\mathbb{T}_{t+1}}\Bigg\{\mathbb{E}^{\alpha|\theta_{t}}\Big[\sum_{s=t}^{\tau}\delta^{s}\big[u_{1,s}(\tilde{\theta}_{s}, \alpha_{s}(\tilde{\theta}_{s}))+\sum_{s=t+1}^{\tau-1} \phi_{s}(\tilde{\theta}_{s}) +\xi_{\tau}(\tilde{\theta}_{\tau}) \big] +\rho(\tau)\Big] \Bigg\} \\
       &-\sup_{\tau\in\mathbb{T}_{t+1}}\Bigg\{\mathbb{E}^{\alpha|\theta_{t},\hat{\theta}_{t}}\Big[\sum_{s=t}^{\tau}\delta^{s}\big[u_{1,s}(\tilde{\theta}_{s}, \alpha_{s}(\tilde{\theta}_{s}))+\sum_{s=t+1}^{\tau-1} \phi_{s}(\tilde{\theta}_{s}) +\xi_{\tau}(\tilde{\theta}_{\tau}) \big] +\rho(\tau)\Big] \Bigg\},
    \end{aligned}
\end{equation*}
which is equal to
\begin{equation*}
    \begin{aligned}
       &\sup_{\tau\in\mathbb{T}_{t+1}}\Bigg\{\mathbb{E}^{\alpha|\hat{\theta}_{t}}\Big[\sum_{s=t}^{\tau}\delta^{s}\big[u_{1,s}(\tilde{\theta}_{s}, \alpha_{s}(\tilde{\theta}_{s})) +\sum_{s=t+1}^{\tau-1} \phi_{s}(\tilde{\theta}_{s}) +\xi_{\tau}(\tilde{\theta}_{\tau}) \big] +\rho(\tau)\Big] \Bigg\}\\
       &-\sup_{\tau\in\mathbb{T}_{t+1}}\Bigg\{\mathbb{E}^{\alpha|\theta_{t},\hat{\theta}_{t}}\Big[\sum_{s=t}^{\tau}\delta^{s}\big[u_{1,s}(\tilde{\theta}_{s}, \alpha_{s}(\tilde{\theta}_{s}))+\sum_{s=t+1}^{\tau-1} \phi_{s}(\tilde{\theta}_{s}) +\xi_{\tau}(\tilde{\theta}_{\tau}) \big] +\rho(\tau)\Big] \Bigg\} \\
       =& \sup_{\tau\in\mathbb{T}_{t+1}}\Big\{\pi^{\alpha}_{t}(\hat{\theta}_{t}) \Big\} - \sup_{\tau\in\mathbb{T}_{t+1}}\Big\{\pi^{\alpha}_{t}(\theta_{t},\hat{\theta}_{t}) \Big\}.
    \end{aligned}
\end{equation*}

Hence, the condition (\ref{eq:potential_function_relation_nS_SufNece}) is satisfied.

\qed

\section{Proof of Proposition \ref{prop:revenue_equivalence_new} }\label{app:proof_prop_revenue_equivalence_new}

Let $\theta^{a}_{\epsilon,t}$, $\theta^{b}_{\epsilon,t}\in\Theta$ associate with $\beta^{\alpha;a}_{\bar{S}, t}$ and $\beta^{\alpha;b}_{\bar{S}, t}$, respectively.
For any period $t\in\mathbb{T}\backslash\{1\}$, time horizon $\tau\in\mathbb{T}$,
\begin{equation*}\label{eq:appE_1}
    \begin{aligned}
        &J{}^{\alpha,\phi^{a},\xi,\rho}_{1,t}(\tau,\theta_{t}|h^{\theta}_{t-1}) - J^{\alpha,\phi^{a},\xi,\rho}_{1,t-1}(\tau,\theta_{t-1}|h^{\theta}_{t-2}) \\
        =& J^{\alpha,\phi^{a},\xi,\rho}_{1,t}(\tau,\theta_{t}|h^{\theta}_{t-1}) - \mathbb{E}^{F_{t}(\theta_{t-1},a_{t-1})}\Big[J^{\alpha,\phi^{a},\xi,\rho}_{1,t-1} (\tau, \tilde{\theta}_{t}|h^{\theta}_{t-2})      \Big]\\
        =& \mathbb{E}^{F_{t}(\theta_{t-1},a_{t-1})}\Big[ J^{\alpha,\phi^{a},\xi,\rho}_{1,t}(\tau,\theta_{t}|h^{\theta}_{t-1}) -  J^{\alpha,\phi^{a},\xi,\rho}_{1,t-1} (\tau, \tilde{\theta}_{t}|h^{\theta}_{t-1})  \Big]\\
        =& \mathbb{E}^{F_{t}(\theta_{t-1},a_{t-1})}\Big[ \int^{\theta_{t}}_{\theta^{a}_{\epsilon,t}} \gamma^{\alpha}_{t-1}(\tau,r|h^{\theta}_{t-2})dr -  \int^{\tilde{\theta}_{t}}_{\theta^{a}_{\epsilon,t}} \gamma^{\alpha}_{t-1}(\tau,r|h^{\theta}_{t-2})dr \Big]\\
        =& \mathbb{E}^{F_{t}(\theta_{t-1},a_{t-1})}\Big[ \int^{\theta_{t}}_{\tilde{\theta}_{t}} \gamma^{\alpha}_{t-1}(\tau,r|h^{\theta}_{t-2})dr \Big]\\
        =& \mathbb{E}^{F_{t}(\theta_{t-1},a_{t-1})}\Big[ \int^{\theta_{t}}_{\theta^{b}_{\epsilon,t}} \gamma^{\alpha}_{t-1}(\tau,r|h^{\theta}_{t-2})dr -  \int^{\tilde{\theta}_{t}}_{\theta^{b}_{\epsilon,t}} \gamma^{\alpha}_{t-1}(\tau,r|h^{\theta}_{t-2})dr \Big]\\
        =&  J^{\alpha,\phi^{b},\xi,\rho}_{1,t}(\tau,\theta_{t}|h^{\theta}_{t-1}) - J^{\alpha,\phi^{b},\xi,\rho}_{1,t-1}(\tau,\theta_{t-1}|h^{\theta}_{t-2}). 
    \end{aligned}
\end{equation*}
Hence, we have
\begin{equation}
    \begin{aligned}
        \mathbb{E}^{\alpha|\theta_{t}}\Big[&\sum_{s=1}^{\tau-1}\delta^{s}\phi^{a}_{s}(\tilde{\theta}_{s})+ \delta^{\tau}\xi^{a}_{\tau}(\tilde{\theta}_{\tau}) +\rho^{a}(\tau) \Big] - \mathbb{E}^{\alpha|\theta_{t}}\Big[\sum_{s=1}^{\tau-1}\delta^{s}\phi^{b}_{s}(\tilde{\theta}_{s}) + \delta^{\tau}\xi^{b}_{\tau}(\tilde{\theta}_{\tau}) +\rho^{b}(\tau)\Big]\\
        =& \mathbb{E}^{\alpha|\theta_{t-1}}\Big[\sum_{s=1}^{\tau-1}\delta^{s}\phi^{a}_{s}(\tilde{\theta}_{s})+ \delta^{\tau}\xi^{a}_{\tau}(\tilde{\theta}_{\tau}) +\rho^{a}(\tau) \Big] - \mathbb{E}^{\alpha|\theta_{t-1}}\Big[\sum_{s=1}^{\tau-1}\delta^{s}\phi^{b}_{s}(\tilde{\theta}_{s}) + \delta^{\tau}\xi^{b}_{\tau}(\tilde{\theta}_{\tau}) +\rho^{b}(\tau)\Big].
    \end{aligned}
\end{equation}
Induction gives the following
\begin{equation}
    \begin{aligned}
        \mathbb{E}^{\alpha|\theta_{t}}\Big[&\sum_{s=1}^{\tau}\delta^{s}\phi^{a}_{s}(\tilde{\theta}_{s})+ \delta^{\tau}\xi^{a}_{\tau}(\tilde{\theta}_{\tau}) +\rho^{a}(\tau) \Big] - \mathbb{E}^{\alpha|\theta_{t}}\Big[\sum_{s=1}^{\tau}\delta^{s}\phi^{b}_{s}(\tilde{\theta}_{s}) + \delta^{\tau}\xi^{b}_{\tau}(\tilde{\theta}_{\tau}) +\rho^{b}(\tau)\Big]\\
        =& \mathbb{E}^{\Xi_{\alpha}}\Big[\sum_{s=1}^{\tau}\delta^{s}\phi^{a}_{s}(\tilde{\theta}_{s})+ \delta^{\tau}\xi^{a}_{\tau}(\tilde{\theta}_{\tau}) +\rho^{a}(\tau) \Big] - \mathbb{E}^{\Xi_{\alpha}}\Big[\sum_{s=1}^{\tau}\delta^{s}\phi^{b}_{s}(\tilde{\theta}_{s}) + \delta^{\tau}\xi^{b}_{\tau}(\tilde{\theta}_{\tau}) +\rho^{b}(\tau)\Big]\\
        =&C_{\tau}.
    \end{aligned}
\end{equation}

\qed

\section{Proof of Lemma \ref{lemma:non-decreasing_J} }\label{app:proof_lemma:non-decreasing_J}

%
It is straightforward to see that
$$
\begin{aligned}
   \frac{ \partial J^{\alpha,\phi,\xi,\rho}_{1,t}(\tau,r|h^{\theta}_{t-1})  }{ \partial r  } \Big|_{r= \theta_{t}} =& \frac{ \partial U^{\alpha,\phi,\xi,\rho}_{t}(\tau,r|h^{\theta}_{t-1})  }{\partial r   }\Big|_{r=\theta_{t}}\\
   =& \mathbb{E}^{\alpha|\theta_{t}}\Big[ \sum_{s=t}^{\tau} \frac{\partial u_{1,s}(r, \alpha_{s}(\tilde{\theta}_{s})) }{\partial r }\Big|_{r=\tilde{\theta}_{s}} G_{t,s}(h^{\tilde{\theta}}_{t,s})   \Big].
\end{aligned}
$$
From Assumption \ref{assp_monotone_transistions}, we have $G_{t,s}(h^{\tilde{\theta}}_{t,s}) \geq 0$. 
Since $u_{1,t}$ is a non-decreasing function of $\theta_{t}$, $\frac{ \partial J^{\alpha,\phi,\xi,\rho}_{1,t}(\tau,r|h^{\theta}_{t-1})  }{ \partial r  }$  $\Big|_{r= \theta_{t}}\geq 0$. 
Therefore, $J^{\alpha,\phi,\xi,\rho}_{1,t}(\tau,r|h^{\theta}_{t-1})$ is a non-decreasing function of $\theta_{t}$, for all $t\in\mathbb{T}$.

\qed

\section{Proof of Proposition \ref{prop:relax_approximation}}\label{app:prop_relax_approximation}

%
From the construction of $\phi$ in (\ref{eq:payment_phi_eta_1}), we have, for any $\tau',\tau''\in\mathbb{T}_{t+1}$,
\begin{equation}\label{eq:suff_nonstop_app_new}
    \begin{split}
        \phi_{t}(\hat{\theta}_{t}) - \phi_{t}(\theta_{t}) 
       = & \beta^{\alpha}_{\bar{S},t}(\hat{\theta}_{t}) - \mathbb{E}^{\alpha|\hat{\theta}_{t}}\Big[\beta_{\bar{S},t+1}(\tilde{\theta}_{t+1})\Big] - u_{1,t}(\hat{\theta}_{t}, \alpha_{t}(\hat{\theta}_{t}))\\
       &-\beta^{\alpha}_{\bar{S},t}(\theta_{t}) + \mathbb{E}^{\alpha|\theta_{t}}\Big[\beta_{\bar{S},t+1}(\tilde{\theta}_{t+1})\Big] + u_{1,t}(\theta_{t}, \alpha_{t}(\theta_{t}))\\
       =& \beta^{\alpha}_{\bar{S},t}(\hat{\theta}_{t}) - \beta^{\alpha}_{\bar{S},t}(\theta_{t}) + \mathbb{E}^{\alpha|\theta_{t}}\Big[ \sum_{s=t}^{\tau''-1}\delta^{s}u_{1,s}(\tilde{\theta}_{s},\alpha_{s}(\tilde{\theta}_{s}))+ \sum_{s=t+1}^{\tau''-1}\delta^{s}\phi_{s}(\tilde{\theta}_{s}) +\beta^{\alpha}_{\bar{S},\tau''}(\tilde{\theta}_{\tau''})  \Big]\\
       &-\mathbb{E}^{\alpha|\hat{\theta}_{t}}\Big[ \sum_{s=t}^{\tau'-1}\delta^{s}u_{1,s}(\tilde{\theta}_{s},\alpha_{s}(\tilde{\theta}_{s}))+ \sum_{s=t+1}^{\tau'-1}\delta^{s}\phi_{s}(\tilde{\theta}_{s}) +\beta^{\alpha}_{\bar{S},\tau'}(\tilde{\theta}_{\tau'}) \Big].
    \end{split}
\end{equation}
%
From the definition of $d^{\alpha}_{\bar{S}, t}$ in (\ref{eq:relax_distance}), we have for any $\tau'$, $\tau''\in\mathbb{T}_{t}$,

\begin{equation}\label{eq:suff_nonstop_app2_new}
    \begin{split}
        &\text{R.H.S. of (\ref{eq:suff_nonstop_app_new})} 
        \leq  \sup_{\tau\in\mathbb{T}_{t}}\Big\{ \mathbb{E}^{\alpha|\hat{\theta}_{t}} \Big[ \sum_{s=t}^{\tau }\delta^{s}\big[u_{1,s}(\tilde{\theta}_{s}, \alpha_{s}(\tilde{\theta}_{s}))\big] +\sum_{s=t+1}^{\tau-1}\delta^{s}\phi_{s}(\tilde{\theta}_{s})  +\delta^{\tau }\xi_{\tau }(\tilde{\theta}_{\tau })   \Big] +\rho(\tau)\Big\}\\
        &-\sup_{\tau\in\mathbb{T}_{t}}\Big\{\mathbb{E}^{\alpha|\theta_{t},\hat{\theta}_{t}}\Big[ \sum_{s=t}^{\tau}\delta^{s}\big[u_{1,s}(\tilde{\theta}_{s}, \alpha_{s}(\tilde{\theta}_{s}))\big] +\sum_{s=t+1}^{\tau-1}\delta^{s}\phi_{s}(\tilde{\theta}_{s})  +\delta^{\tau}\xi_{\tau}(\tilde{\theta}_{\tau })   \Big] +\rho(\tau) \Big\}+d^{\alpha}_{\bar{S}, t}\\
        &+ \mathbb{E}^{\alpha|\theta_{t}}\Big[ \sum_{s=t}^{\tau''-1}\delta^{s}u_{1,s}(\tilde{\theta}_{s},\alpha_{s}(\tilde{\theta}_{s}))+ \sum_{s=t+1}^{\tau''-1}\delta^{s}\phi_{s}(\tilde{\theta}_{s}) +\beta^{\alpha}_{\bar{S},\tau''}(\tilde{\theta}_{\tau''})  \Big]\\
       &-\mathbb{E}^{\alpha|\hat{\theta}_{t}}\Big[ \sum_{s=t}^{\tau'-1}\delta^{s}u_{1,s}(\tilde{\theta}_{s},\alpha_{s}(\tilde{\theta}_{s}))+ \sum_{s=t+1}^{\tau'-1}\delta^{s}\phi_{s}(\tilde{\theta}_{s}) +\beta^{\alpha}_{\bar{S},\tau'}(\tilde{\theta}_{\tau'}) \Big].
    \end{split}
\end{equation}

From the condition (\ref{eq:potential_function_relation_ine}), $\beta^{\alpha}_{\bar{S},t}(\theta_{t})\geq \beta^{\alpha}_{S,t}(\theta_{t})$, for all $\theta_{t}\in\Theta_{t}$, $t\in\mathbb{T}$. Hence, $\beta^{\alpha}_{\bar{S},t}(\theta_{t})\geq $ $ \xi_{t}(\theta_{t}) + $  $ u_{1,t}(\theta_{t},$  $\alpha_{t}(\theta_{t}))$, for all $\theta_{t}\in \mathbb{T}$, $t\in\mathbb{T}$.
Then, 
\begin{equation}\label{eq:app_haha1}
    \begin{split}
        \sup_{\tau\in\mathbb{T}_{t}}\Big\{ \mathbb{E}^{\alpha|\hat{\theta}_{t}} &\Big[ \sum_{s=t}^{\tau }\delta^{s}\big[u_{1,s}(\tilde{\theta}_{s}, \alpha_{s}(\tilde{\theta}_{s}))\big] +\sum_{s=t+1}^{\tau-1}\delta^{s}\phi_{s}(\tilde{\theta}_{s})  +\delta^{\tau }\xi_{\tau }(\tilde{\theta}_{\tau })   +\rho(\tau)\Big]\\
        &-\mathbb{E}^{\alpha|\hat{\theta}_{t}}\Big[ \sum_{s=t}^{\tau-1}\delta^{s}u_{1,s}(\tilde{\theta}_{s},\alpha_{s}(\tilde{\theta}_{s}))+ \sum_{s=t+1}^{\tau-1}\delta^{s}\phi_{s}(\tilde{\theta}_{s}) +\beta^{\alpha}_{\bar{S},\tau}(\tilde{\theta}_{\tau})  \Big]\Big\}\\
        \leq\sup_{\tau\in\mathbb{T}_{t}}\Big\{ \mathbb{E}^{\alpha|\hat{\theta}_{t}} &\Big[ \sum_{s=t}^{\tau }\delta^{s}\big[u_{1,s}(\tilde{\theta}_{s}, \alpha_{s}(\tilde{\theta}_{s}))\big] +\sum_{s=t+1}^{\tau-1}\delta^{s}\phi_{s}(\tilde{\theta}_{s})  +\delta^{\tau }\xi_{\tau }(\tilde{\theta}_{\tau })   +\rho(\tau)\Big]\\
        &-\mathbb{E}^{\alpha|\hat{\theta}_{t}}\Big[ \sum_{s=t}^{\tau}\delta^{s}u_{1,s}(\tilde{\theta}_{s},\alpha_{s}(\tilde{\theta}_{s}))+ \sum_{s=t+1}^{\tau-1}\delta^{s}\phi_{s}(\tilde{\theta}_{s}) +\delta^{\tau}\xi_{\tau}(\tilde{\theta}_{\tau})  \Big]\Big\}\\
        =\sup_{\tau\in\mathbb{T}_{t}}\Big\{ \rho(\tau)&\Big\}.
    \end{split}
\end{equation}
Hence, from (\ref{eq:app_haha1}), we have, for any $\tau'\in \mathbb{T}_{t}$
\begin{equation}\label{eq:app_94}
    \begin{aligned}
   \text{R.H.S. of (\ref{eq:suff_nonstop_app2_new}) } 
       %
       \leq&\mathbb{E}^{\alpha|\theta_{t}}\Big[ \sum_{s=t}^{\tau'-1}\delta^{s}u_{1,s}(\tilde{\theta}_{s},\alpha_{s}(\tilde{\theta}_{s}))+ \sum_{s=t+1}^{\tau'-1}\delta^{s}\phi_{s}(\tilde{\theta}_{s}) +\beta^{\alpha}_{\bar{S},\tau'}(\tilde{\theta}_{\tau'})  \Big] + \sup_{\tau\in\mathbb{T}_{t}}\Big\{\rho(\tau) \Big\}+d^{\alpha}_{\bar{S}, t}\\ &-\sup_{\tau\in\mathbb{T}_{t}}\Big\{\mathbb{E}^{\alpha|\theta_{t},\hat{\theta}_{t}}\Big[ \sum_{s=t}^{\tau}\delta^{s}\big[u_{1,s}(\tilde{\theta}_{s}, \alpha_{s}(\tilde{\theta}_{s}))\big] +\sum_{s=t+1}^{\tau-1}\delta^{s}\phi_{s}(\tilde{\theta}_{s})  +\delta^{\tau}\xi_{\tau}(\tilde{\theta}_{\tau })    \Big] +\rho(\tau) \Big\}.
\end{aligned}
\end{equation}
From the construction of $\rho$ in (\ref{eq:construct_rho_primal}) and Lemma \ref{lemma:monotonicity_of_L}, we have, for some $\tau'\in\mathbb{T}_{t}$

\begin{equation*}
    \begin{aligned}
       \text{R.H.S. of (\ref{eq:app_94})} \leq&\mathbb{E}^{\alpha|\theta_{t}}\Big[\sum_{s=t}^{\tau'}\delta^{s} u_{1,s}(\tilde{\theta}_{s},\alpha_{s}(\tilde{\theta}_{s}))+ \sum_{s=t+1}^{\tau'-1}\delta^{s}\phi_{s}(\tilde{\theta}_{s}) +\delta^{\tau'}\xi_{\tau'}(\tilde{\theta}_{\tau'}) +\rho(\tau')\Big]+d^{\alpha}_{\bar{S}, t}+\sup_{\tau\in\mathbb{T}_{t}}\Big\{\rho(\tau) \Big\}\\ &-\sup_{\tau\in\mathbb{T}_{t}}\Big\{\mathbb{E}^{\alpha|\theta_{t},\hat{\theta}_{t}}\Big[ \sum_{s=t}^{\tau}\delta^{s}\big[u_{1,s}(\tilde{\theta}_{s}, \alpha_{s}(\tilde{\theta}_{s}))\big] +\sum_{s=t+1}^{\tau-1}\delta^{s}\phi_{s}(\tilde{\theta}_{s})  +\delta^{\tau}\xi_{\tau}(\tilde{\theta}_{\tau }) +\rho(\tau)\Big] \Big\}\\
       =& \sup_{\tau\in\mathbb{T}_{t}}\Big\{\mathbb{E}^{\alpha|\theta_{t}}\Big[\sum_{s=t}^{\tau}\delta^{s} u_{1,s}(\tilde{\theta}_{s},\alpha_{s}(\tilde{\theta}_{s}))+ \sum_{s=t+1}^{\tau-1}\delta^{s}\phi_{s}(\tilde{\theta}_{s}) +\delta^{\tau}\xi_{\tau}(\tilde{\theta}_{\tau}) +\rho(\tau)\Big]\Big\}\\ &-\sup_{\tau\in\mathbb{T}_{t}}\Big\{\mathbb{E}^{\alpha|\theta_{t},\hat{\theta}_{t}}\Big[ \sum_{s=t}^{\tau}\delta^{s}\big[u_{1,s}(\tilde{\theta}_{s}, \alpha_{s}(\tilde{\theta}_{s}))\big] +\sum_{s=t+1}^{\tau-1}\delta^{s}\phi_{s}(\tilde{\theta}_{s})  +\delta^{\tau}\xi_{\tau}(\tilde{\theta}_{\tau }) +\rho(\tau)\Big] \Big\}\\
       &+d^{\alpha}_{\bar{S}, t}+\sup_{\tau\in\mathbb{T}_{t}}\Big\{\rho(\tau) \Big\},
    \end{aligned}
\end{equation*}
which is equal to
\begin{equation*}
    \begin{aligned}
       & \sup_{\tau\in\mathbb{T}_{t}}\Big\{\mathbb{E}^{\alpha|\theta_{t}}\Big[\sum_{s=t}^{\tau}\delta^{s} u_{1,s}(\tilde{\theta}_{s},\alpha_{s}(\tilde{\theta}_{s}))+ \sum_{s=t+1}^{\tau-1}\delta^{s}\phi_{s}(\tilde{\theta}_{s}) +\delta^{\tau}\xi_{\tau}(\tilde{\theta}_{\tau}) +\rho(\tau)\Big]\Big\}\\ &-\sup_{\tau\in\mathbb{T}_{t}}\Big\{\mathbb{E}^{\alpha|\theta_{t},\hat{\theta}_{t}}\Big[ \sum_{s=t}^{\tau}\delta^{s}\big[u_{1,s}(\tilde{\theta}_{s}, \alpha_{s}(\tilde{\theta}_{s}))\big] +\sum_{s=t+1}^{\tau-1}\delta^{s}\phi_{s}(\tilde{\theta}_{s})  +\delta^{\tau}\xi_{\tau}(\tilde{\theta}_{\tau }) +\rho(\tau)\Big] \Big\}\\
       &+d^{\alpha}_{\bar{S}, t}+\sup_{\tau\in\mathbb{T}_{t}}\Big\{\rho(\tau) \Big\}.
    \end{aligned}
\end{equation*}
Hence, from the definition of $U^{\alpha,\phi,\xi,\rho}_{\bar{S},t}$, we have
$$
U^{\alpha,\phi,\xi,\rho}_{\bar{S},t}(\theta_{t}|h^{\theta}_{t-1})+d^{\alpha}_{\bar{S}, t}+\sup_{\tau\in\mathbb{T}_{t}}\Big\{\rho(\tau) \Big\}\geq U^{\alpha,\phi,\xi,\rho}_{\bar{S},t}(\theta_{t}, \hat{\theta}_{t}|h^{\theta}_{t-1}).
$$
Following the similar way, we can prove the following
$$
U_{S,t}^{\alpha,\phi, \xi, \rho}(\theta_{t}|h^{\theta}_{t-1})+ d^{\alpha}_{S,t}\geq U_{S,t}^{\alpha,\phi, \xi, \rho}(\theta_{t}, \hat{\theta}_{t}|h^{\theta}_{t-1}),
$$
where $d^{\alpha}_{S,t}$ is defined in (\ref{eq:relax_distance_stop}).

Therefore, we can conclude that the mechanism is $\Big\{d^{\alpha^{\circ}}_{S,t}, d^{\alpha^{\circ}}_{\bar{S},t}+\sup_{\tau\in\mathbb{T}_{t}}\Big\{\rho^{\circ}(\tau) \Big\}\Big\}$-DIC if $d^{\alpha}_{S,t}>0$ and $\Big[d^{\alpha}_{t}+ \sup_{\tau\in\mathbb{T}_{t}}$ $\Big\{\rho(\tau) \Big\}>0$; otherwise, the mechanism is DIC.
We can prove the case when the payment rule $\rho$ is constructed according to (\ref{eq:rho_construction_multiple}) by following the similar way.

\qed


\bibliography{mybibfile}


\end{document}